\documentclass[moor]{informs1}              

\usepackage{latexsym}
\usepackage{epsfig}
\usepackage{subfig}
\usepackage{amsmath}
\usepackage{amssymb}
\usepackage{graphicx}
\usepackage{color}
\usepackage{multirow}
\usepackage{endnotes}
\usepackage{enumerate}
\usepackage{bm}
\usepackage{booktabs}
\usepackage[compact]{titlesec}
\usepackage{lscape}
\usepackage[normalem]{ulem}
\usepackage{enumitem}
\usepackage{caption}
\newcommand{\ot}[1]{\tilde{#1}}

\newcommand{\ru}[1]{\big\lceil #1 \big\rceil}
\newcommand{\rd}[1]{\lfloor #1 \rfloor}
\newcommand{\ep}{\epsilon}
\newcommand{\al}{\alpha}

\newcommand{\E}{\mathbf{E}}
\newcommand{\pr}{\mathbf{P}}

\newcommand{\dr}{\mathrm{d}}
\newcommand{\e}{\mathrm{e}}

\newcommand{\D}{\mathbb{D}}

\newcommand{\I}{\mathbb{I}}
\newcommand{\N}{\mathbb{N}}
\newcommand{\R}{\mathbb{R}}

\newcommand{\mA}{\mathcal{A}}
\newcommand{\mB}{\mathcal{B}}

\newcommand{\mF}{\mathcal{F}}

\newcommand{\mI}{\mathcal{I}}
\newcommand{\mJ}{\mathcal{J}}
\newcommand{\mK}{\mathcal{K}}
\newcommand{\mN}{\mathcal{N}}

\newcommand{\mX}{\mathcal{X}}

\newcommand{\eProof}{\hspace*{.1in}\hfill  $\blacksquare$}
\newcommand{\bProof}{\noindent\textbf{Proof: }}

\usepackage{natbib}
 \NatBibNumeric
 \def\bibfont{\small}%
 \bibpunct[, ]{[}{]}{,}{n}{}{,}%

\usepackage[colorlinks=true,breaklinks=true,bookmarks=true,urlcolor=blue,
     citecolor=blue,linkcolor=blue,bookmarksopen=false,draft=false,hypertexnames=false]{hyperref}

\def\EMAIL#1{\href{mailto:#1}{#1}}

\TheoremsNumberedThrough     

\EquationsNumberedThrough    

\begin{document}


\RUNAUTHOR{\"Ozkan}

\RUNTITLE{Control of Fork-Join Processing Networks}

\TITLE{Control of Fork-Join Processing Networks with Multiple Job Types and Parallel Shared Resources}

\ARTICLEAUTHORS{
\AUTHOR{Erhun \"Ozkan}
\AFF{College of Administrative Sciences and Economics, Ko\c c University, Istanbul, Turkey, \EMAIL{erhozkan@ku.edu.tr}}
} 

\ABSTRACT{%
A fork-join processing network is a queueing network in which tasks associated with a job can be processed simultaneously. Fork-join processing networks are prevalent in computer systems, healthcare, manufacturing, project management, justice system, etc. Unlike the conventional queueing networks, fork-join processing networks have synchronization constraints that arise due to the parallel processing of tasks and can cause significant job delays. We study scheduling control in fork-join processing networks with multiple job types and parallel shared resources. Jobs arriving in the system fork into arbitrary number of tasks, then those tasks are processed in parallel, and then they join and leave the network. There are shared resources processing multiple job types. We study the scheduling problem for those shared resources (that is, which type of job to prioritize at any given time) and propose an asymptotically optimal scheduling policy in diffusion scale.
}


\KEYWORDS{Fork-join processing network; scheduling control; asymptotic optimality; diffusion scale}
\MSCCLASS{60K25, 90B22, 90B36, 93E20, 60F17}
\ORMSCLASS{Primary: Stochastic model applications, queues, optimization ; secondary: probability, diffusion models, limit theorems }
\HISTORY{}

\maketitle

%

\section{Introduction}\label{s_intro}

A fork-join processing network is a queueing network in which tasks associated with a job can be processed simultaneously. Fork-join networks are prevalent in computer systems (see \citet{tho14}, \citet{zen18}), healthcare (see \citet{arm15}, \citet{car18}), manufacturing (see \citet{dal92}), project management (see \citet{adl95}), justice system (see \citet{lar93}), etc.

\begin{figure}[htb]
\begin{center}
\includegraphics[width=0.8\textwidth]{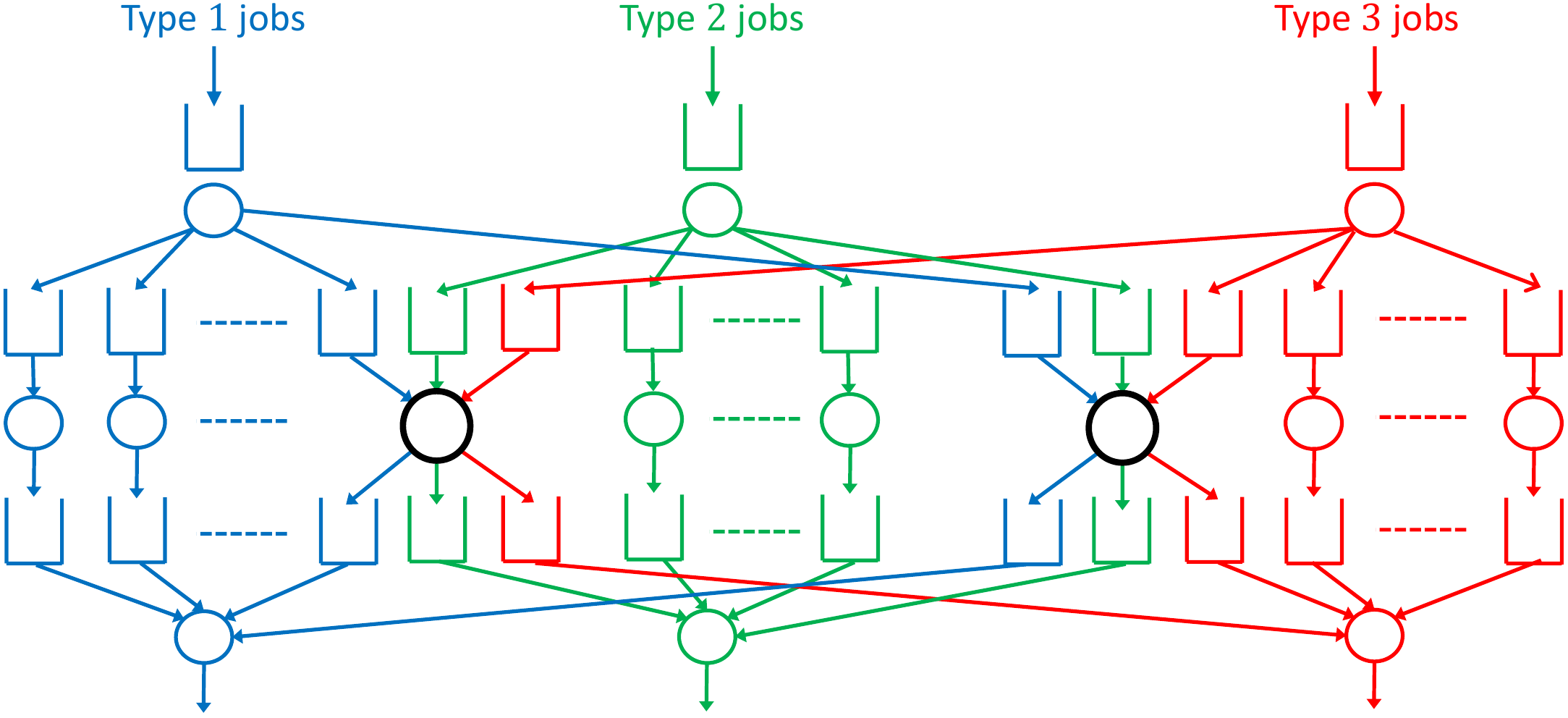}
\caption{(Color online) A fork-join processing network with 3 job types and 2 shared servers. Each job type is forked into arbitrary but finite number of tasks. Circles denote servers and bins denote buffers. There are four different types of servers: fork, join, dedicated, and shared servers. Dedicated servers process a single job type and the shared servers process multiple job types.}\label{fj_network_0}
\end{center}
\end{figure}

We study scheduling decisions in fork-join networks with multiple customer classes that share multiple processing resources. Our main motivation is patient-flow process in emergency departments (EDs, see Figure 1 in \citet{car18}). After triage, a patient may need to have some lab tests (e.g., blood, urine), radiology exams (e.g., CT scan, X-ray, ultra sound), etc. Some of those tests and exams can be taken simultaneously. For example, while his/her blood sample is analyzed, a patient can have a CT scan. A patient cannot be discharged until all of the test results are ready. Therefore, the patient-flow diagram can be illustrated as the fork-join processing network depicted in Figure \ref{fj_network_0}, in which job (patient) types represent condition severity of the patients and the resources (servers) represent the labs or facilities where the tests and exams are taken. Resources such as CT scanners have a large impact on patient waiting time (see \citet{hub11}) because they are very expensive and so hospitals generally own at most a few of them. This motivates us to study the problem of how to schedule resources that are used by multiple different job types.

\begin{figure}[htb]
\begin{center}
\includegraphics[width=0.23\textwidth]{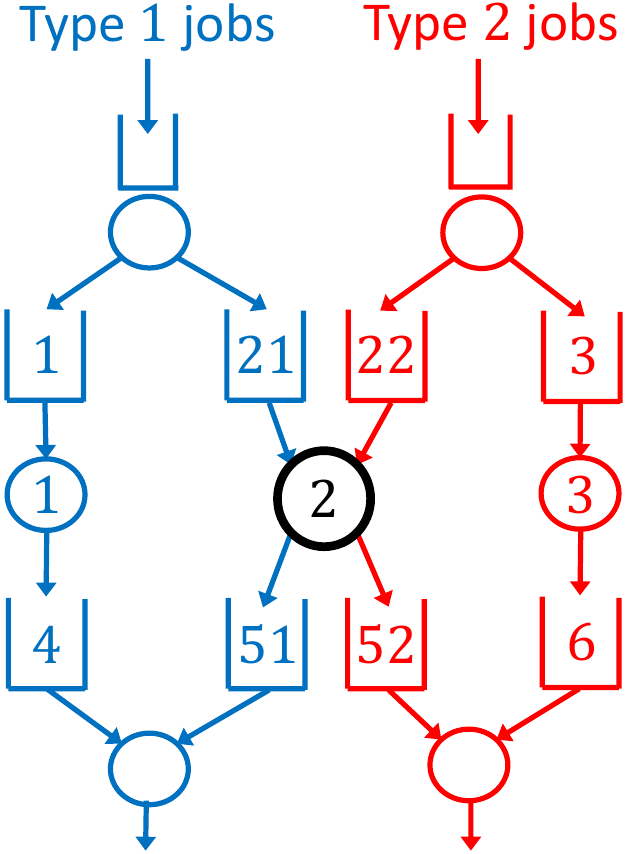}
\caption{(Color online) A fork-join processing network with 2 job types and a single shared server.}\label{fj_network_1}
\end{center}
\end{figure}

The parallel processing of tasks gives rise to synchronization constraints which can cause job delays. Although delays in fork-join networks can be approximated under the first-in-first-out (FIFO) scheduling discipline (see \citet{ngu93,ngu94}), FIFO scheduling rule does not necessarily minimize delay (see \citet{ata12} and \citet{ozk19}). To see this, let us consider the simple fork-join network in Figure \ref{fj_network_1}. There are synchronization constraints because type 1 (2) jobs cannot be joined until there is at least one job in both buffers 4 and 51 (52 and 6). Server 2 processes both job types, but can only serve one job at a time. The control decision is to decide which job type server 2 should prioritize. Suppose that $h_1\mu_{21} \geq h_2\mu_{22}$, where $h_1$ ($h_2$) denotes the holding cost per a type $1$ ($2$) job per unit time and $\mu_{21}$ ($\mu_{22}$) denotes the service rate of server 2 for type 1 (2) jobs. According to the $c\mu$ rule, server 2 should always give priority to type $1$ jobs. However, if there are multiple jobs waiting in buffers 51 and 6 and no jobs waiting in buffers 4 and 52, it may be better to have server 2 work on a type 2 job instead of a type 1 job. This is because server 1 and 2 block the join operations of the type 1 and 2 jobs, respectively. Therefore, static scheduling rules such as FIFO or $c\mu$ rule can perform poorly in the fork-join network in Figure \ref{fj_network_1}.

Deriving an exact optimal control policy is very challenging even for the simple network in Figure \ref{fj_network_1}. A potential approach is to use Markov Decision Process (MDP) techniques under the assumption that the interarrival and service times are exponentially distributed. However, because the associated state is the number of jobs in each buffer (that is, a 10-dimensional state space), curse of dimensionality arises. Therefore, a more efficient solution approach is to derive asymptotically optimal control policies in the conventional heavy-traffic regime as done by \citet{ozk19}. They prove asymptotic optimality of a continuous-review and state-dependent control policy in diffusion scale under the assumption that server 2 is in heavy traffic, that is, its processing capacity is barely enough to process all incoming jobs. Otherwise, the scheduling control in server 2 has negligible impact on job delays. However, it is not clear how to extend the results of \citet{ozk19} to more general fork-join networks (see Section \ref{s_lit} for details). Furthermore, there are only a few studies in the literature that consider control of fork-join networks. For example, similar to \citet{ozk19}, \citet{ata12} consider the control of a very specific fork-join network. \citet{ped14a,ped14b,ped17} consider fork-join networks with very general topological structure but they focus on throughput optimality and ignore delay minimization. Consequently, control of fork-join networks is a relatively unexplored research area, even though those networks are prevalent in many application domains.

Our main contribution is an asymptotically optimal control policy in diffusion scale for the fork-join network in Figure \ref{fj_network_0} extended with arbitrary but finite number of job types and arbitrary but finite number of shared servers. The objective is minimizing the expected total discounted holding cost. We assume that all of the shared servers are in conventional heavy-traffic regime. Otherwise, if a shared server is in light traffic, that is, if its processing capacity is more than enough to process the incoming jobs, then the scheduling decisions in that shared server will not be very important because any work-conserving policy will perform well, when considered in the heavy-traffic regime. We also assume that the join servers are in light traffic so that the synchronization constraints are the main reason for delays in the join operations. We do not have any assumption on the processing capacities of the fork servers or the dedicated servers, that is, a fork server or a dedicated server can be in either heavy or light traffic.

The proposed policy is a continuous-review, state-dependent, and non-preemptive policy under which a linear program (LP) is solved at discrete time epochs. The parameters of the LP are holding cost rates of the job types, service rates for the job types in the shared servers, the numbers of jobs waiting in front of the dedicated servers, and the weighted total number of jobs waiting in front of each shared server. The decision variables of the LP are the numbers of each job type that should wait in front of the shared servers. After the LP is solved, the system controller compares the numbers of jobs waiting in front of a shared server with an optimal LP solution. If those numbers are different, then the shared server processes the jobs until the numbers of jobs waiting in front of it becomes sufficiently close to the optimal LP solution. At that point, the system controller resolves the LP and follows the same procedure. Therefore, under the proposed policy, the numbers of jobs in front of the shared servers always track an optimal LP solution. If the LP has multiple optimal solutions at a time epoch, then we need to choose an optimal solution which does not deviate a lot from the previous optimal LP solutions. We accomplish this goal by solving a quadratic program (QP) which finds an optimal LP solution with the desired property. The QP is convex and so is solvable in polynomial time.

The proposed policy does not require the knowledge of arrival rates of the job types. This is important in practice because estimating the arrival rates accurately can be difficult in many applications. For example, in the ED case, arrival rates of the patients can change dramatically over time. Many studies in the queueing literature such as \citet{bel01}, \citet{ata05}, \citet{dai08}, and \citet{ozk19} prove asymptotic optimality of preemptive control policies due to their mathematical simplicity. In contrast, our proposed policy is non-preemptive which has a practical appeal.

We use Harrison's classical scheme in the paper (see \citet{har97} and \citet{har00}): We formulate a diffusion control problem (DCP), next solve the DCP and interpret a control policy from the solution, and then prove the asymptotic optimality of the proposed policy. The main technical challenge for our paper is that the resulting DCP is multidimensional. Specifically, the dimension of the resulting DCP is equal to the number of shared servers in the system and so the resulting workload process is multidimensional. Although there are many studies considering one-dimensional workload process (see for example \citet{bel01,man04,sto04,ata05,dai08,ozk19}), studies considering multidimensional workload process is rare (see for example \citet{pes16}). This is because solving a multidimensional DCP and proving the asymptotic optimality of the control policy interpreted from the DCP solution are generally very challenging. We overcome this challenge by utilizing the special structure of the fork-join network that we consider. Specifically, because the shared servers are parallel to each other and the join servers are in light-traffic, the network effect is limited in the resulting DCP and we are able to prove that under any work-conserving policy, the multidimensional workload process weakly converges to the same limit. Consequently, the DCP is time-decomposable, and we can convert it to an LP and solve it numerically at discrete time epochs. Finally, by tracking the optimal LP solutions in the shared servers and utilizing a Lipschitz continuity result associated with the optimal LP solutions, we prove asymptotic optimality of the proposed policy.

We present a literature review in Section \ref{s_lit} and some notation in Section \ref{s_not}. Then, we present the model description and the objective in Section \ref{s_model}. We present the asymptotic framework in Section \ref{s_asym} and derive an asymptotic lower bound on the performance of any admissible policies in Section \ref{s_alb}. We present the formal definition of the proposed policy and prove its asymptotic optimality in Section \ref{s_policy}. Finally, we present some modeling extensions in Section \ref{s_ext}. All of the proofs are presented in either the appendix or the electronic companion.

\subsection{Literature Review}\label{s_lit}

Although there are many studies focusing on performance evaluation of the fork-join networks (see \citet{ngu93,ngu94}, \citet{tho14} and references therein, \citet{lu16a,lu16b,lu17}), there are only a few studies focusing on control of fork-join networks (see \citet{ata12}, \citet{ped14a,ped14b,ped17}, \citet{ozk19}). \citet{ata12} consider the control of a specific fork-join network with probabilistic feedback mechanism. Their motivation is also patient flow process in EDs and the feedback represents cases in which a patient should retake a radiology exam or have a lab test again. In contrast, there is no feedback in the network that we consider. \citet{ped14a,ped14b,ped17} consider the control of fork-join networks with very general topological structure in discrete time. Their focus is throughput optimality instead of delay minimization. However, in the fork-join network that we consider, any work-conserving control policy maximizes the throughput, but average job waiting time can differ significantly among the work-conserving policies  (see the numerical experiments in Section E of the E-companion of \citet{ozk19}). Hence, we focus on delay (or in general holding cost) minimization.

\citet{ozk19} consider the control of the fork-join network in Figure \ref{fj_network_1}. They also use Harrison's classical scheme in their paper. Because there is a single shared server in their network, the resulting DCP is one-dimensional. They find a closed-form solution to the DCP and prove weak convergence of the queue length processes to the closed-form DCP solution. However, it is not clear how to extend their results to fork-join networks with more than one shared servers. First, finding a closed-form solution to multidimensional DCPs is very challenging, if not impossible. Although \citet{ozk19} are able to derive a closed-form solution to a specific two-dimensional DCP, the policy that they interpret is complicated enough such that it is not clear how to extend their asymptotic optimality proof to that case. Under the policy that they interpret, the shared servers change the job types that they prioritize frequently depending on the system state, which complicates proving weak convergence of the individual queue length processes to the closed-form DCP solution. In contrast, we convert the DCP into an LP, solve the LP numerically in discrete-time epochs, and use a simple policy which keeps the queue lengths close to the optimal LP solutions. Consequently, we are able to prove the asymptotic optimality of our proposed policy for networks with arbitrary number of job types and shared servers. 


There are also studies focusing on throughput scalability of fork-join networks (see \citet{zen18} and references therein). \citet{zen18} call a network throughput scalable if throughput does not decrease to zero as the network size grows to infinity. They provide necessary and sufficient conditions on the throughput scalability of fork-join networks with general topological structure.

\subsection{Notation}\label{s_not}

The set of nonnegative and strictly positive integers are denoted by $\N$ and $\N_+$, respectively. For all $n\in\N_+$, $\R^n$ denotes the $n$-dimensional Euclidean space and $\R_+^n$ denotes the nonnegative orthant in $\R^n$. For any $x,y\in \R$, $x\vee y:= \max\{x,y\}$, $x\wedge y:= \min\{x,y\}$, and $(x)^+:=x\vee 0$. For any $\bm{x}:=(x_1,x_2,\ldots,x_n)\in\R^n$ and $\bm{y}:=(y_1,y_2,\ldots,y_n)\in\R^n$, we let $|\bm{x}-\bm{y}|_{\infty}:=\max_{i\in\{1,2,\ldots,n\}}|x_i-y_i|$. For any $x\in \R$, $\rd{x}$ ($\lceil x \rceil$) denotes the greatest (smallest) integer which is smaller (greater) than or equal to $x$. For any given set $\mX$, $|\mX|$ denotes the cardinality of $\mX$.

For all $n\in\N_+$, $\D^n$ denotes the set of functions $f:\R_+\rightarrow\R^n$ that are right continuous with left limits. We let $\textbf{0},e\in \D$ be such that $\textbf{0}(t)=0$ and $e(t)=t$ for all $t\in\R_+$. For $x,y\in \D$, $x\vee y$, $x\wedge y$, and $(x)^+$ are functions in $\D$ such that $(x\vee y)(t):=x(t)\vee y(t)$, $(x\wedge y)(t):= x(t)\wedge y(t)$, and $(x)^+(t):=(x(t))^+$ for all $t\in\R_+$. For any $x\in \D$, we define the mappings $\Psi, \Phi:\D\rightarrow\D$ such that for all $t\in\R_+$,
\begin{equation}\label{eq_rm}
\Psi(x)(t) := \sup_{0\leq s\leq t} (-x(s))^+,\quad\quad \Phi(x)(t) := x(t)+\Psi(x)(t),
\end{equation}
where $\Phi$ is the one-sided and one-dimensional reflection map (see Chapter 13.5 of \citet{whi02}). For $x\in \D$ and $t\in \R_+$, we let $\Vert x \Vert_t := \sup_{0\leq s \leq t} |x(s)|$. We consider $\D^n$ endowed with the usual Skorokhod $J_1$ topology (see Chapter 3 of \citet{bil99}). Let $\mB(\D^n)$ denote the Borel $\sigma$-algebra on $\D^n$ associated with Skorokhod $J_1$ topology. For stochastic processes $\{W^r$, $r\in\N_+\}$ and $W$ whose sample paths are in $\D^n$ for some $n\in\N_+$, ``$W^r \Rightarrow W$'' means that the probability measures induced by $\{W^r$, $r\in\N_+\}$ on $(\D^n,\mB(\D^n))$ converge weakly to the one induced by $W$ on $(\D^n,\mB(\D^n))$ as $r\rightarrow \infty$. All of the convergence results hold as $r\rightarrow\infty$.

Let $\mN=\{1,2,\ldots,n\}$ and $x_i\in\D$ for all $i\in \mN$. Then $(x_i,i\in \mN)$ denotes the process $(x_1,x_2,\ldots,x_n)$ in $\D^n$. We abbreviate the phrase ``uniformly on compact intervals'' by ``u.o.c.'' and ``almost surely'' by ``a.s.''. We let $\xrightarrow{a.s.}$ denote almost sure convergence. We repeatedly use the fact that convergence in the $J_1$ metric is equivalent to u.o.c. convergence when the limit process is continuous (see page 124 in \citet{bil99}). Let $\{x^r,r\in\N\}$ be a sequence in $\D$ and $x\in\D$. Then $x^r\rightarrow x$ u.o.c., if $\Vert x^r-x\Vert_{t}\rightarrow 0$ for all $t\in\R_+$. We let ``$\circ$'' denote the composition map and $\I$ denote the indicator function. We assume that all the random variables and stochastic processes are defined in the same complete probability space $(\Omega, \mF, \pr)$, $\E$ denotes expectation under $\pr$, and $\pr(A, B):=\pr(A\cap B)$.

\section{Model Description}\label{s_model}

There are $J\in\N_+$ different job types arriving in the network and we let $\mJ:=\{1,2,\ldots,J\}$ denote the set of job types. For all $j\in\mJ$, each incoming type $j$ job is first forked into arbitrary but finite number of jobs. Some of those forked jobs are processed in some of the shared servers and the remaining ones are processed in the dedicated servers associated with type $j$ jobs. Then all of those forked jobs are joined together and leave the system. We assume that the fork and join operations are done instantaneously to simplify the notation. Later, we will relax this assumption in Section \ref{s_ext_ni}. Consequently, there are two different server types in the network: dedicated and shared servers. Dedicated servers process only a single job type. In contrast, shared servers process at least two job types. Each server can process at most a single job at a time. 

There are $I$ different shared servers and we let $\mI:=\{1,2,\ldots,I\}$ denote the set of shared servers. For all $i\in\mI$ and $j\in\mJ$, if type $j$ jobs are processed in shared server $i$, then we let $P_{ij}:=1$; otherwise, $P_{ij}:=0$. We let $\mI_j:=\{i\in\mI: P_{ij}=1\}$ for all $j\in\mJ$ and $\mJ_i:=\{j\in\mJ: P_{ij}=1\}$ for all $i\in\mI$. Thus, $\mI_j$ is the set of shared servers that process type $j$ jobs and $\mJ_i$ is the set of job types that are processed in the shared server $i$. We assume that $|\mI_j|\geq 1$ for all $j\in\mJ$, implying that each job type is processed in at least one shared server (otherwise there is no scheduling decision for that job type). We also assume that $|\mJ_i|\geq 2$ for all $i\in\mI$, implying that each shared server processes at least two job types (otherwise that server is not a shared server by definition).

For all $j\in\mJ$, each incoming type $j$ job is first forked into $K_j + |\mI_j|$ job types where $K_j\in\N$ denote the number of dedicated servers that process type $j$ jobs. We let $\mK_j$ denote the set of dedicated servers associated with the type $j$ jobs. If $K_j =0$, then $\mK_j=\emptyset$. By definition, $\mK_j\cap\mK_l=\emptyset$ for all $j,l\in\mJ$ such that $j\neq l$. The join operation of a type $j$ job happens when all of the $K_j + |\mI_j|$ forked jobs are processed in the associated dedicated and shared servers. 

There are $\sum_{j\in\mJ} 2(K_j + |\mI_j| )$ buffers in the network such that each buffer has infinite capacity, half of the buffers are in the upper layer, and the remaining half are in the lower layer. In the upper layer, there exists a buffer in front of each dedicated server. Moreover, for all $j\in\mJ$ and $i\in\mI_j$, there exists a buffer in front of the shared server $i$ in which type $j$ jobs wait for service. In the lower layer, there exists a buffer after each dedicated server in which jobs processed in the dedicated server wait for the join operation. Furthermore, for all $j\in\mJ$ and $i\in\mI_j$, there exists a buffer after the shared server $i$ in which type $j$ jobs processed in the shared server $i$ wait for the join operation.

\subsection{Stochastic Primitives}\label{s_primitives}

\textbf{\textit{External arrivals}} We associate the external arrival times of type $j\in\mJ$ jobs with strictly positive and independent and identically distributed (i.i.d.) sequence of random variables $\{\bar{u}_{jn},n\in\N_+\}$ and the constant $\lambda_j>0$. For all $j\in\mJ$ and $n\in\N_+$, $\E[\bar{u}_{jn}]=1$, the variance of $\bar{u}_{jn}$ is $\beta_j^2$, and $u_{jn}:=\bar{u}_{jn}/\lambda_j$ denotes the inter-arrival time between the $(n-1)$st and $n$th type $j$ job. Then, for all $j\in\mJ$, $\{u_{jn},n\in\N_+\}$ is an i.i.d. sequence of random variables with mean $1/\lambda_j$ and squared coefficient of variance $\beta_j^2$. For all $j\in\mJ$, $n\in\N_+$, and $t\in\R_+$, we let $U_j(0):=0$ and
\begin{equation*}
U_j(n):= \sum_{l=1}^n u_{jl},\hspace{1cm} A_j(t):= \sup\left\{n\in\N: U_j(n)\leq t\right\}.
\end{equation*}
Then, $A_j$ is a renewal process such that $A_j(t)$ is the number of external type $j$ job arrivals up to time $t\in\R_+$.

\textbf{\textit{Service processes in the dedicated servers}} For all $j\in\mJ$ and $k\in\mK_j$, let $\{v_{kn},n\in\N_+\}$ be a strictly positive and i.i.d. sequence of random variables with mean $1/\mu_{k}$ and squared coefficient of variance $\sigma_{k}^2$. We let $v_{kn}$ denote the service time of the $n$th type $j$ job in the dedicated server $k$ for all $j\in\mJ$, $k\in\mK_j$, and $n\in\N_+$. For all $j\in\mJ$, $k\in\mK_j$, $n\in\N_+$, and $t\in\R_+$, let $V_{k}(0):=0$ and
\begin{equation*}
V_{k}(n):= \sum_{l=1}^n v_{kl},\hspace{1cm} S_{k}(t):= \sup\left\{n\in\N: V_{k}(n)\leq t\right\}.
\end{equation*}
Then, $S_{k}$ is a renewal process such that $S_{k}(t)$ is the number of service completions in the dedicated server $k\in\mK_j$ up to time $t\in\R_+$ given that the dedicated server never idles during $[0,t]$.

\textbf{\textit{Service processes in the shared servers}} For all $j\in\mJ$ and $i\in\mI_j$, let $\{v_{ijn},n\in\N_+\}$ be a strictly positive and i.i.d. sequence of random variables with mean $1/\mu_{ij}$ and squared coefficient of variance $\sigma_{ij}^2$. We let $v_{ijn}$ denote the service time of the $n$th type $j$ job in the shared server $i$ for all $j\in\mJ$, $i\in\mI_j$, and $n\in\N_+$. For all $j\in\mJ$, $i\in\mI_j$, $n\in\N_+$, and $t\in\R_+$, let $V_{ij}(0):=0$ and
\begin{equation*}
V_{ij}(n):= \sum_{l=1}^n v_{ijl},\hspace{1cm} S_{ij}(t):= \sup\left\{n\in\N: V_{ij}(n)\leq t\right\}.
\end{equation*}

For all $j\in\mJ$, $k\in\mK_j$, and $i\in\mI_j$, we assume that the sequences $\{\bar{u}_{jn},n\in\N_+\}$, $\{v_{kn},n\in\N_+\}$, and $\{v_{ijn},n\in\N_+\}$ are mutually independent of each other and of all other stochastic primitives.

\subsection{Network Dynamics and Scheduling Control}\label{s_dynamics}

For all $j\in\mJ$, $k\in\mK_j$, $i\in\mI_j$ and $t\in\R_+$, we let $T_{k}(t)$ denote the cumulative amount of time that the dedicated server $k$ works on type $j$ jobs during $[0,t]$ and $T_{ij}(t)$ denote the cumulative amount of time that the shared server $i$ works on type $j$ jobs during $[0,t]$. The scheduling control is defined by the process $(T_{ij},j\in\mJ,i\in\mI_j)$. For all $t\in\R_+$, we let
\begin{subequations}\label{eq_idle}
\begin{align}
&I_{k}(t):= t - T_{k}(t)\quad\text{for all $j\in\mJ$ and $k\in\mK_j$},\label{eq_idle_1}\\
&I_i(t):= t - \sum_{j\in\mJ_i} T_{ij}(t)\quad\text{for all $i\in\mI$},\label{eq_idle_2} 
\end{align}
\end{subequations}
denote the cumulative idle time of the dedicated server $k\in\mK_j$ and the shared server $i\in\mI$ up to time $t$, respectively.

For all $j\in\mJ$, $k\in\mK_j$, and $t\in\R_+$, we let $Q_k(t)$ denote the number of type $j$ jobs waiting in front of the dedicated server $k$ at time $t$, including the job that is in service; and we let $Q_k^{(1)}(t)$ denote the number of type $j$ jobs waiting after the dedicated server $k$ for the join operation at time $t$. For all $j\in\mJ$, $i\in\mI_j$, and $t\in\R_+$, we let $Q_{ij}(t)$ denote the number of type $j$ jobs waiting to be served by the shared server $i$ at time $t$, including the job that is in service; and we let $Q_{ij}^{(1)}(t)$ denote the number of type $j$ jobs waiting after the shared server $i$ for the join operation at time $t$. Then, for all $j\in\mJ$, $k\in\mK_j$, $i\in\mI_j$, and $t\in\R_+$,
\begin{subequations}\label{eq_queue}
\begin{align}
& Q_k(t)= Q_k(0)+A_j(t)-S_{k}(T_{k}(t)),\label{eq_queue_1}\\
& Q_{ij}(t)= Q_{ij}(0)+A_j(t)-S_{ij}(T_{ij}(t)),\label{eq_queue_2}\\
& Q_k(t),\; Q_{ij}(t),\; Q_k^{(1)}(t),\; Q_{ij}^{(1)}(t) \geq 0,\label{eq_queue_3}
\end{align}
\end{subequations}
where $S_{k}(T_{k}(t))$ and $S_{ij}(T_{ij}(t))$ denote the cumulative number of type $j$ jobs processed in the dedicated server $k$ and in the shared server $i$ up to time $t$, respectively.

Let $N_j(t)$ denote the number of type $j\in\mJ$ jobs in the system at time $t\in\R_+$ by counting a job that is forked into multiple jobs as a single job. Then, for all $j\in\mJ$ and $t\in\R_+$, we have
\begin{subequations}\label{eq_sync}
\begin{align}
&N_j(t) = Q_k(t) + Q_k^{(1)}(t) = Q_{ij}(t) + Q_{ij}^{(1)}(t)\quad\text{for all $k\in\mK_j$ and $i\in\mI_j$},\label{eq_sync_1}\\
&\min_{k\in\mK_j} Q_k^{(1)}(t) \wedge \min_{i\in\mI_j} Q_{ij}^{(1)}(t) = 0,\label{eq_sync_2}
\end{align}
\end{subequations}
where \eqref{eq_sync_2} is because the join operations happen instantaneously.

For all $j\in\mJ$, $k\in\mK_j$, $i\in\mI_j$, and $t\in\R_+$, we have
\begin{subequations}\label{eq_hl}
\begin{align}
&V_k(S_k(T_k(t))) \leq T_k(t) < V_k(S_k(T_k(t))+1),\label{eq_hl_1}\\
&V_{ij}(S_{ij}(T_{ij}(t))) \leq T_{ij}(t) < V_{ij}(S_{ij}(T_{ij}(t))+1),\label{eq_hl_2}
\end{align}
\end{subequations}
which implies that we consider only head-of-the-line (HL) policies, where jobs are processed in FIFO order within each buffer. Notice that a forked job associated with a specific job cannot join a forked job originating in another job under the HL policies.

For all $j\in\mJ$, $k\in\mK_j$, $i\in\mI$, and $t\in\R_+$, we have
\begin{subequations}\label{eq_wc}
\begin{align}
&\text{$I_k$ is nondecreasing, $I_k(0)=0$, and $I_k(t)$ increases if and only if $Q_k(t)=0$} ,\label{eq_wc_1}\\
&\text{$I_i$ is nondecreasing, $I_i(0)=0$, and $I_i(t)$ increases if and only if $\max_{j\in\mJ_i} Q_{ij}(t)=0$} ,\label{eq_wc_2}
\end{align}
\end{subequations}
which implies that all of the servers work in a work-conserving fashion. We assume that holding cost rate per job per unit time does not change when a job is served in a dedicated or shared server. Therefore, work-conserving policies are more efficient than non-work-conserving policies.

\begin{definition} (Admissible policies)
A scheduling policy $\pi:=(T_{ij},j\in\mJ,i\in\mI_j )$ is admissible if the processes $(T_{k},T_{ij},j\in\mJ,k\in\mK_j,i\in\mI_j)$, $(I_{k},I_{i}, j\in\mJ, k\in\mK_j,i\in\mI)$, and $(Q_{k},Q_{ij},Q_{k}^{(1)},Q_{ij}^{(1)}, j\in\mJ, k\in\mK_j,i\in\mI_j)$ satisfy \eqref{eq_idle}, \eqref{eq_queue}, \eqref{eq_sync}, \eqref{eq_hl}, \eqref{eq_wc}; and for all $j\in\mJ$ and $i\in\mI_j$, we have
\begin{subequations}\label{eq_admissible}
\begin{align}
&\text{$T_{ij}(t)$ is $\mF$-measurable (that is, $T_{ij}(t)\in\mF$) for all $t\in\R_+$},\label{eq_admissible_1}\\
&\text{$T_{ij}$ is continuous and nondecreasing and $T_{ij}(0)=0$}.\label{eq_admissible_2}
\end{align}
\end{subequations}
\end{definition}

Condition \eqref{eq_admissible_1} implies that the set of admissible policies includes even the ones that can anticipate the future.

\subsection{Objective}\label{s_objective}

Our objective is to minimize the expected total discounted holding cost. Let $h_j\in\R_+$ denote the holding cost rate per a type $j$ job per unit time for all $j\in\mJ$. We assume that $\max_{j\in\mJ} h_j>0$. Let $\delta > 0$ be the discount parameter and $\Pi$ denote the set of admissible policies. Then, we want to find
\begin{equation}\label{eq_obj_1}
\argmin_{\pi\in\Pi} \sum_{j\in\mJ} h_j  \E\left[ \int_0^\infty \e^{-\delta t} N_j^\pi(t) \dr t\right].
\end{equation}
We will first focus on the following objective: For any given $\ep>0$ and $t \in\R_+$, we want to find
\begin{equation}\label{eq_obj_2}
\argmin_{\pi\in\Pi} \;\pr\Bigg(\sum_{j\in\mJ} h_j N_j^\pi(t) > \ep  \Bigg).
\end{equation}
Then, we will focus on the objective \eqref{eq_obj_1}. Observe that any admissible policy that minimizes the objective \eqref{eq_obj_2} for all $\ep>0$ and $t \in\R_+$ also minimizes the objective \eqref{eq_obj_1}.

\section{Asymptotic Framework}\label{s_asym}

Deriving an optimal control policy for the fork-join network described in Section \ref{s_model} is very challenging. A potential approach is to use MDP techniques under the assumption that the inter-arrival and service times are exponentially distributed. However, because the associated state is the number of jobs in each buffer, curse of dimensionality arises. Therefore, a more efficient solution approach is to derive asymptotically optimal control policies in the conventional heavy-traffic regime in diffusion scale. Specifically, we assume that all of the shared servers are in heavy traffic. We do not have any assumption on the processing capacities of the dedicated servers, that is, a dedicated server can be in either heavy or light traffic. 

First, we introduce a sequence of fork-join networks and present our main assumptions in Section \ref{s_sequence}. Then, we present fluid and diffusion scaled processes and two convergence results that hold under any work-conserving policy in Section \ref{s_fluid}.

\subsection{A Sequence of Fork-Join Networks}\label{s_sequence}

We consider a sequence of fork-join networks indexed by $r\in\N_+$. Each fork-join network has the same structure with the original network defined in Section \ref{s_model} except that the constant $\lambda_j$ depends on $r$ for all $j\in\mJ$. Specifically, in the $r$th system, we associate the inter-arrival times of type $j\in\mJ$ jobs with the sequence of random variables $\{\bar{u}_{jn},n\in\N_+\}$, defined in Section \ref{s_primitives}, and the constant $\lambda_j^r>0$. For all $j\in\mJ$, $r\in\N_+$, and $n\in\N_+$, we let $u_{jn}^r:=\bar{u}_{jn}/\lambda_j^r$ denote the inter-arrival time between the $(n-1)$st and $n$th type $j$ job in the $r$th system. Then, in the $r$th system, arrival rate of type $j$ jobs is $\lambda_j^r$, whereas the squared coefficient of variation of the inter-arrival times is $\beta_j^2$, which is equal to the one in the original system. From this point forward, we will use the superscript $r$ to show the dependence of the stochastic processes to the $r$th fork-join network.

Next we present two main assumptions. The first one is the exponential moment assumption for the inter-arrival and service times.
\begin{assumption}\label{a_moment} (Moment)
There exists an $\bar{\al}>0$ such that for all $\al\in(-\bar{\al},\bar{\al})$,
\begin{equation*}
\E\left[ \e^{\al\bar{u}_{j1}}\right]<\infty,\quad\E\left[ \e^{\al v_{k1}}\right]<\infty,\quad\E\left[ \e^{\al v_{ij1}}\right]<\infty,\quad\forall j\in\mJ,k\in\mK_j,i\in\mI_j.
\end{equation*}
\end{assumption}
Exponential moment assumption is common in the queueing literature, see for example \citet{har98,bel01,mag03,mey03,ozk19}.

The second assumption sets up the asymptotic regime.
\begin{assumption}\label{a_ht} (Asymptotic Regime)
\begin{enumerate}
\item $\lambda_j^r \rightarrow \lambda_j$  for all $j\in\mJ$.

\item $\sum_{j\in\mJ_i} \lambda_j/ \mu_{ij} =1$ for all $i\in\mI$.

\item $r\left( \left(\sum_{j\in\mJ_i} \lambda_j^r/ \mu_{ij} \right)-1\right) \rightarrow \theta_i\in\R$ for all $i\in\mI$.

\item $r(\lambda_j^r - \mu_k)\rightarrow \theta_k\in\R\cup\{-\infty\}$ for all $j\in\mJ$ and $k\in\mK_j$.
\end{enumerate}
\end{assumption}

If a shared server is in light traffic, any admissible policy will perform well in that shared server and so the control will become trivial. Therefore, we assume that all shared servers are in heavy traffic in Parts 2 and 3 of Assumption \ref{a_ht}. Part 4 of Assumption \ref{a_ht} states that the dedicated servers can be in either light or heavy traffic. On the one hand, if $\theta_k=-\infty$ for some $j\in\mJ$ and $k\in\mK_j$, then the dedicated server $k$ is in light traffic. On the other hand, if $\theta_k\in\R$, then the dedicated server $k$ is in heavy traffic. For all $j\in\mJ$, we let $\mK_j^L:=\{k\in\mK_j: \theta_k=-\infty\}$ and $\mK_j^H:=\{k\in\mK_j: \theta_k\in\R\}$. Then, $\mK_j^L$ ($\mK_j^H$) denotes the set of dedicated servers associated with type $j$ jobs which are in light (heavy) traffic and $\mK_j^L\cup \mK_j^H = \mK_j$ for all $j\in\mJ$. Because $\lambda_j>0$ for all $j\in\mJ$ and $|\mJ_i|\geq 2$ for all $i\in\mI$, we have $\mu_{ij}>\lambda_j$ for all $i\in\mI$ and $j\in\mJ_i$ by Assumption \ref{a_ht} Part 2.

For simplicity, we assume that the system is initially empty, that is, $Q_k^r(0)=Q_{ij}^r(0)=Q_k^{(1),r}(0)=Q_{ij}^{(1),r}(0)=0$ for all $j\in\mJ$, $k\in\mK_j$, $i\in\mI_j$, and $r\in\N_+$. We relax this assumption in Section \ref{s_ext_i}.

\subsection{Fluid and Diffusion Scaled Processes}\label{s_fluid}

For all $j\in\mJ$, $k\in\mK_j$, $t\in\R_+$, and $r\in\N_+$, the fluid scaled processes are defined as 
\begin{subequations}\label{eq_fluid}
\begin{align}
& \bar{A}_j^r(t):=r^{-2}A_j^r(r^2t)&& \bar{N}_j^r(t):=r^{-2}N_j^r(r^2t) \label{eq_fluid_1},\\
& \bar{S}_k^r(t):=r^{-2}S_k(r^2t)&& \bar{S}_{ij}^r(t):=r^{-2}S_{ij}(r^2t)\quad\forall i\in\mI_j ,\label{eq_fluin_2}\\
&\bar{T}_k^r(t):=r^{-2}T_k^r(r^2t) && \bar{T}_{ij}^r(t):=r^{-2}T_{ij}^r(r^2t)\quad\forall i\in\mI_j,\label{eq_fluin_3}\\
&\bar{I}_k^r(t):=r^{-2}I_k^r(r^2t) && \bar{I}_{i}^r(t):=r^{-2}I_{i}^r(r^2t)\quad\forall i\in\mI,\label{eq_fluid_4}\\
&\bar{Q}_k^r(t):=r^{-2}Q_k^r(r^2t) && \bar{Q}_{ij}^r(t):=r^{-2}Q_{ij}^r(r^2t)\quad\forall i\in\mI_j,\label{eq_fluid_5}\\
&\bar{Q}_k^{(1),r}(t):=r^{-2}Q_k^{(1),r}(r^2t) && \bar{Q}_{ij}^{(1),r}(t):=r^{-2}Q_{ij}^{(1),r}(r^2t)\quad\forall i\in\mI_j.\label{eq_fluid_6}
\end{align}
\end{subequations}

For all $j\in\mJ$, $k\in\mK_j$, $t\in\R_+$, and $r\in\N_+$, the diffusion scaled processes are defined as
\begin{subequations}\label{eq_diffusion}
\begin{align}
& \hat{A}_j^r(t):=r\left(\bar{A}_j^r(t)-\lambda_j^r t\right)&& \hat{N}_j^r(t):=r\bar{N}_j^r(t) \label{eq_diffusion_1},\\
& \hat{S}_k^r(t):=r\left(\bar{S}_k^r(t)-\mu_kt\right)&& \hat{S}_{ij}^r(t):=r\left(\bar{S}_{ij}^r(t)-\mu_{ij}t\right)\quad\forall i\in\mI_j ,\label{eq_diffusion_2}\\
& \hat{T}_k^r(t):=r\bar{T}_k^r(t) && \hat{T}_{ij}^r(t):=r\bar{T}_{ij}^r(t)\quad\forall i\in\mI_j ,\label{eq_diffusion_3}\\
&\hat{I}_k^r(t):=r \bar{I}_k^r(t) && \hat{I}_{i}^r(t):=r \bar{I}_{i}^r(t)\quad\forall i\in\mI,\label{eq_diffusion_4}\\
&\hat{Q}_k^r(t):=r \bar{Q}_k^r(t) && \hat{Q}_{ij}^r(t):=r\bar{Q}_{ij}^r(t)\quad\forall i\in\mI_j,\label{eq_diffusion_5}\\
&\hat{Q}_k^{(1),r}(t):=r\bar{Q}_k^{(1),r}(t) && \hat{Q}_{ij}^{(1),r}(t):=r\bar{Q}_{ij}^{(1),r}(t)\quad\forall i\in\mI_j.\label{eq_diffusion_6}
\end{align}
\end{subequations}

For all $i\in\mI$, $t\in\R_+$, and $r\in\N_+$, we define the workload process in the shared server $i$ as
\begin{equation}\label{eq_workload}
W_i^r(t):=\sum_{j\in\mJ_i} \frac{Q_{ij}^r(t)}{\mu_{ij}}.
\end{equation}
Then, $W_i^r(t)$ is the expected time that the shared server $i$ should spend in order to process all of the jobs in front of it given that no more jobs arrive in the system. We let $\bar{W}_i^r(t) := r^{-2}W_i^r(r^2t)$ and $\hat{W}_i^r(t):=r\bar{W}_i^r(t)$ denote the fluid and diffusion scaled workload in the shared server $i$, respectively, for all $i\in\mI$, $t\in\R_+$, and $r\in\N_+$.

Next, we present a convergence result for the fluid scaled processes.

\begin{proposition}\label{p_fluid}
Let $\pi=\{\pi^r,r\in\N_+\}$ be an arbitrary sequence of admissible policies. Then,
\begin{align*}
&\left( \bar{Q}_k^r, \bar{Q}_k^{(1),\pi,r},\bar{T}_{k}^r, j\in\mJ,k\in\mK_j;\;\; \bar{Q}_{ij}^{\pi,r}, \bar{Q}_{ij}^{(1),\pi,r},\bar{T}_{ij}^{\pi,r}, j\in\mJ,i\in\mI_j;\;\; \bar{W}_i^{\pi,r},i\in\mI \right)\\
&\hspace{2cm}\xrightarrow{a.s.} \left( \bar{Q}_k, \bar{Q}_k^{(1)},\bar{T}_{k},j\in\mJ, k\in\mK_j;\;\; \bar{Q}_{ij}, \bar{Q}_{ij}^{(1)},\bar{T}_{ij}, j\in\mJ,i\in\mI_j;\;\; \bar{W}_i,i\in\mI \right)\quad\text{u.o.c.},
\end{align*}
where $\bar{Q}_k = \bar{Q}_{ij}=\bar{Q}_k^{(1)}=\bar{Q}_{ij}^{(1)}=\bm{0}$ for all $j\in\mJ$, $k\in\mK_j$, $i\in\mI_j$; $\bar{W}_i=\bm{0}$ for all $i\in\mI$; and $\bar{T}_k(t)=(\lambda_j/\mu_k)t$ and $\bar{T}_{ij}(t)=(\lambda_j/\mu_{ij})t$ for all $j\in\mJ$, $k\in\mK_j$, $i\in\mI_j$, $t\in\R_+$.
\end{proposition}
The proof of Proposition \ref{p_fluid} follows from standard methodology and so we skip it. For a similar proof, see the proof of Proposition 1 in \citet{ozk19}. We will use Proposition \ref{p_fluid} to prove a weak convergence result for the diffusion scaled processes.

For all $j\in\mJ$, $k\in\mK_j$, $i\in\mI$, and $t\in\R_+$, let
\begin{align*}
&\hat{X}_k^r(t) := \hat{A}_j^r(t) -  \hat{S}_k^r\circ\bar{T}_k^r(t) + r\left(\lambda_j^r-\mu_k \right)t,\\
&\hat{X}_i^r(t) := \sum_{j\in\mJ_i} \frac{1}{\mu_{ij}}\left(\hat{A}_j^r(t) -  \hat{S}_{ij}^r\circ\bar{T}_{ij}^r(t)\right) + r\left( \Bigg(\sum_{j\in\mJ_i} \frac{\lambda_j^r}{\mu_{ij}} \Bigg)-1\right)t.
\end{align*}
After some algebra, for all $j\in\mJ$, $k\in\mK_j$, $i\in\mI$, and $t\in\R_+$, we have
\begin{equation*}
\hat{Q}_k^r(t)=\hat{X}_k^r(t) + \mu_k \hat{I}_k^r(t),\qquad \hat{W}_i^r(t)=\hat{X}_i^r(t) + \hat{I}_i^r(t).
\end{equation*}
Under any admissible policy, by \eqref{eq_rm}, for all $j\in\mJ$, $k\in\mK_j$, and $i\in\mI$, we have
\begin{equation*}
\left(\hat{Q}_k^r, \mu_k\hat{I}_k^r\right) = \left(\Phi,\Psi\right)\left(\hat{X}_k^r\right),\qquad \left(\hat{W}_i^r, \hat{I}_i^r\right) = \left(\Phi,\Psi\right)\left(\hat{X}_i^r\right).
\end{equation*}

Let $a_1:= I+ \sum_{j\in\mJ} |\mK_j^H| $ and $0^{(a_1)}$ denote the origin in $\R^{a_1}$. Let us define the $a_1$-dimensional vector $\Theta:=\left(\theta_k,j\in\mJ,k\in\mK_j^H,\;\theta_i,i\in\mI\right)$ and the $(a_1\times a_1)$-dimensional positive definite matrix $\Sigma$ such that
\begin{align*}
& \Sigma_{kl}:= \lambda_j\left( \beta_j^2 + \sigma_k^2\I(k=l)\right)\qquad\text{for all $j\in\mJ$ and $k,l\in\mK_j^H$},\\
& \Sigma_{ki}:= \lambda_j\beta_j^2/\mu_{ij} \qquad\text{for all $j\in\mJ$ and $k\in\mK_j^H$ and $i\in\mI_j$},\\
& \Sigma_{in}:= \sum_{j\in\mJ_i\cap\mJ_n}\lambda_j\left( \beta_j^2 + \sigma_{ij}^2\I(i=n)\right) /(\mu_{ij}\mu_{nj}) \qquad\text{for all $i,n\in\mI$},
\end{align*}
and all of the remaining components of $\Sigma$ are equal to 0. Let $R$ be a $(a_1\times a_1)$-dimensional diagonal matrix such that $R_{kk}:=\mu_k$ and $R_{ii}:=1$ for all $j\in\mJ$, $k\in\mK_j^H$, and $i\in\mI$. Then, we have the following weak convergence result.

\begin{proposition}\label{p_diffusion}
Let $\pi=\{\pi^r,r\in\N_+\}$ be an arbitrary sequence of admissible policies. Then,
\begin{equation*}
\left( \hat{Q}_k^r, j\in\mJ,k\in\mK_j,\; \hat{W}_i^{\pi,r},i\in\mI \right)\Rightarrow \left( \ot{Q}_k, j\in\mJ,k\in\mK_j,\; \ot{W}_i,i\in\mI \right),
\end{equation*}
where $\ot{Q}_k=\bm{0}$ for all $j\in\mJ$ and $k\in\mK_j^L$ and $\left( \ot{Q}_k, j\in\mJ,k\in\mK_j^H,\;\ot{W}_i,i\in\mI \right)$ is a semimartingale reflected Brownian motion (SRBM) associated with the data $\left(\R_+^{a_1},\Theta, \Sigma, R,0^{(a_1)} \right)$. $\R_+^{a_1}$ is the state space of the SRBM; $\Theta$ and $\Sigma$ are the drift vector and the covariance matrix of the underlying Brownian motion of the SRBM, respectively; $R$ is the reflection matrix; and $0^{(a_1)}$ is the starting point of the SRBM.
\end{proposition}

The formal definition of an SRBM can be found in Definition 3.1 of \citet{wil98b}. Since the proof of Proposition \ref{p_diffusion} follows from standard methodology, we skip it. For a similar proof, see the proof of Proposition 2 in \citet{ozk19}.

Proposition \ref{p_diffusion} implies that the diffusion scaled workload processes in the shared servers (see \eqref{eq_workload}) converge to the same limit under any sequence of admissible policies. Therefore, the important question is how to split those workloads to the buffers in front of the shared servers in order to minimize the cost.

Next, we will derive an asymptotic lower bound on the performance of any admissible policy.

\section{Asymptotic Lower Bound}\label{s_alb}

We derive an asymptotic lower bound on the performance of any sequence of admissible policies with respect to the objective \eqref{eq_obj_2}. We construct an approximating DCP in Section \ref{s_dcp} and derive the asymptotic lower bound by the solution of the aforementioned DCP in Section \ref{s_alb2}.

\subsection{Approximating Diffusion Control Problem}\label{s_dcp}

In this section, we construct an approximating DCP whose solution will help us to derive an asymptotic lower bound with respect to the objective \eqref{eq_obj_2} in Section \ref{s_alb2}.

By \eqref{eq_sync_1}, for all $j\in\mJ$, $k\in\mK_j$, $i\in\mI_j$, $t\in\R_+$, and $r\in\N_+$, we have
\begin{equation}\label{eq_sync_3}
Q_k^{(1),r}(t) = N_j^r(t) - Q_k^r(t),\qquad  Q_{ij}^{(1),r}(t) = N_j^r(t) - Q_{ij}^r(t).
\end{equation}
By \eqref{eq_sync_2} and \eqref{eq_sync_3}, for all $j\in\mJ$ and $t\in\R_+$, we have
\begin{align*}
& \min_{k\in\mK_j} \left( N_j^r(t) - Q_k^r(t) \right) \wedge \min_{i\in\mI_j} \left( N_j^r(t) - Q_{ij}^r(t) \right)=0,\\
\implies & \left( N_j^r(t)- \max_{k\in\mK_j} Q_k^r(t) \right) \wedge  \left( N_j^r(t)- \max_{i\in\mI_j} Q_{ij}^r(t) \right)=0,\\
\implies & N_j^r(t) = \max_{k\in\mK_j} Q_k^r(t) \vee  \max_{i\in\mI_j} Q_{ij}^r(t).
\end{align*}
Then, parallel to the objective \eqref{eq_obj_2}, for any given $t\in\R_+$ and $r\in\N_+$, let us consider the diffusion scaled objective of minimizing
\begin{equation}\label{eq_obj_3}
\sum_{j\in\mJ} h_j \hat{N}_j^r(t) =\sum_{j\in\mJ} h_j \left(\max_{k\in\mK_j} \hat{Q}_k^r(t) \vee  \max_{i\in\mI_j} \hat{Q}_{ij}^r(t)\right).
\end{equation}
At this point, let us assume that
\begin{equation}\label{eq_hc}
\left( \hat{Q}_{ij}^r, j\in\mJ,i\in\mI_j \right)\Rightarrow \left( \ot{Q}_{ij}, j\in\mJ,i\in\mI_j  \right).
\end{equation}
Then, by \eqref{eq_workload}, \eqref{eq_obj_3}, \eqref{eq_hc}, and Proposition \ref{p_diffusion}, we construct the following DCP: For any $t\in\R_+$,
\begin{subequations}\label{eq_dcp}
\begin{align}
\min\;& \sum_{j\in\mJ} h_j \left(\max_{k\in\mK_j} \ot{Q}_k(t) \vee  \max_{i\in\mI_j} \ot{Q}_{ij}(t)\right)\label{eq_dcp_1}\\
\text{such that (s.t.) } & \sum_{j\in\mJ_i} \frac{\ot{Q}_{ij}(t)}{\mu_{ij}} = \ot{W}_i(t)\qquad\forall i\in\mI,\label{eq_dcp_2}\\
& \ot{Q}_{ij}(t)\geq 0\qquad\forall j\in\mJ,\;i\in\mI_j,\label{eq_dcp_3}
\end{align}
\end{subequations}
where the decision variables are $\big(\ot{Q}_{ij}(t),j\in\mJ,i\in\mI_j\big)$. The objective \eqref{eq_dcp_1} minimizes the total holding cost rate at time $t$. The constraints \eqref{eq_dcp_2} and \eqref{eq_dcp_3} state that we should split the workload of each shared server to the buffers in front of that shared server in order to minimize the total holding cost. For fixed $t\in\R_+$, the DCP \eqref{eq_dcp} has linear constraints and a convex objective, thus it is a convex problem. Furthermore, we can linearize the DCP \eqref{eq_dcp}. Let $a_2:= I+ \sum_{j\in\mJ} K_j $ and  $a_3:= J+ \sum_{j\in\mJ} |\mI_j| $. 

\begin{lemma}\label{l_lp}
Let $(h_j,\mu_{ij},j\in\mJ,i\in\mI_j)$ be constants such that $h_j\in\R_+$ and $\mu_{ij}>0$ for all $j\in\mJ$ and $i\in\mI_j$. For given $\bm{b}:=(q_k,j\in\mJ,k\in\mK_j, w_i,i\in\mI)\in\R_+^{a_2}$, consider the convex problem
\begin{subequations}\label{eq_cvx}
\begin{align}
\min\;& \sum_{j\in\mJ} h_j \left(\max_{k\in\mK_j} q_k\vee \max_{i\in\mI_j} q_{ij}\right)\label{eq_cvx_1}\\
\text{s.t. }\; & \sum_{j\in\mJ_i} \frac{q_{ij}}{\mu_{ij}} = w_i\qquad\forall i\in\mI,\label{eq_cvx_2}\\
& q_{ij}\geq 0\qquad\forall j\in\mJ,\;i\in\mI_j,\label{eq_cvx_3}
\end{align}
\end{subequations}
where the decision variables are $\big(q_{ij},j\in\mJ,i\in\mI_j\big)$. Next, consider the LP
\begin{subequations}\label{eq_lp}
\begin{align}
\min\;& \sum_{j\in\mJ} h_j y_j\label{eq_lp_1}\\
\text{s.t. }\; & y_j\geq q_k\qquad\forall j\in\mJ,\;k\in\mK_j,\label{eq_lp_2}\\
& y_j- q_{ij}\geq 0\qquad\forall j\in\mJ,\;i\in\mI_j,\label{eq_lp_3}\\
& \sum_{j\in\mJ_i} \frac{q_{ij}}{\mu_{ij}} = w_i\qquad\forall i\in\mI,\label{eq_lp_4}\\
& q_{ij}\geq 0\qquad\forall j\in\mJ,\;i\in\mI_j,\label{eq_lp_5}
\end{align}
\end{subequations}
where the decision variables are $\big( y_j,q_{ij},j\in\mJ,i\in\mI_j\big)$. Then, we have the following results:
\begin{enumerate}
\item Let $\big( \ot{y}_j,\ot{q}_{ij},j\in\mJ,i\in\mI_j\big)$ be an arbitrary optimal solution of the LP \eqref{eq_lp}. Then, $\big( \ot{q}_{ij},j\in\mJ,i\in\mI_j\big)$ is an optimal solution of the convex problem \eqref{eq_cvx}. Moreover, the optimal objective function value of the convex problem \eqref{eq_cvx} and the LP \eqref{eq_lp} are the same.

\item Let $z:\R_+^{a_2}\rightarrow \R_+$ be such that $z(\bm{b})$ denotes the optimal objective function value of the LP \eqref{eq_lp} for all $\bm{b}\in\R_+^{a_2}$. Then, for any given $\bm{b}^{(1)}:=(q_k^{(1)},j\in\mJ,k\in\mK_j, w_i^{(1)},i\in\mI)\in\R_+^{a_2}$ and $\bm{b}^{(2)}:=(q_k^{(2)},j\in\mJ,k\in\mK_j, w_i^{(2)},i\in\mI)\in\R_+^{a_2}$,
\begin{equation*}
\left| z\big(\bm{b}^{(1)}\big)-z\big(\bm{b}^{(2)}\big)\right| \leq C_1 \big| \bm{b}^{(1)}-\bm{b}^{(2)} \big|_{\infty},
\end{equation*}
where $C_1>0$ is a constant dependent on the objective coefficients and left-hand-side (LHS) parameters of the constraints of the LP \eqref{eq_lp}.

\item For given $\bm{b}=(q_k,j\in\mJ,k\in\mK_j, w_i,i\in\mI)\in\R_+^{a_2}$, consider the QP:
\begin{subequations}\label{eq_qp}
\begin{align}
\min\;& \sum_{j\in\mJ} y_j^2+ \sum_{j\in\mJ}\sum_{i\in\mI_j} q_{ij}^2 \label{eq_qp_1}\\
\text{s.t. }\; & y_j\geq q_k\qquad\forall j\in\mJ,\;k\in\mK_j,\label{eq_qp_2}\\
& y_j- q_{ij}\geq 0\qquad\forall j\in\mJ,\;i\in\mI_j,\label{eq_qp_3}\\
& \sum_{j\in\mJ_i} \frac{q_{ij}}{\mu_{ij}} = w_i\qquad\forall i\in\mI,\label{eq_qp_4}\\
&\sum_{j\in\mJ} h_j y_j \leq z(\bm{b}),\label{eq_qp_5}\\
& q_{ij}\geq 0\qquad\forall j\in\mJ,\;i\in\mI_j,\label{eq_qp_6}
\end{align}
\end{subequations}
where the decision variables are $\big(y_j,q_{ij},j\in\mJ,i\in\mI_j\big)$. For each $\bm{b}\in\R_+^{a_2}$, there exists a unique optimal solution of the QP \eqref{eq_qp}. Let $\big( y_j^{(1)},q_{ij}^{(1)},j\in\mJ,i\in\mI_j\big)$ and $\big( y_j^{(2)},q_{ij}^{(2)},j\in\mJ,i\in\mI_j\big)$ be the unique optimal solutions of the QP \eqref{eq_qp} under $\bm{b}^{(1)}$ and $\bm{b}^{(2)}$, respectively, where $\bm{b}^{(1)}\in\R_+^{a_2}$ and $\bm{b}^{(2)}\in\R_+^{a_2}$. Then,
\begin{equation*}
\max_{j\in\mJ,i\in\mI_j} \left| q_{ij}^{(1)}-q_{ij}^{(2)}\right| \leq C_2 \big| \bm{b}^{(1)}-\bm{b}^{(2)} \big|_{\infty},
\end{equation*}
where $C_2>0$ is a constant dependent on the LHS parameters of the constraints of the QP \eqref{eq_qp}.

\end{enumerate}
\end{lemma}
The proof of Lemma \ref{l_lp} is presented in Appendix \ref{l_lp_proof}. The first part of Lemma \ref{l_lp} states that we can solve the convex problem \eqref{eq_cvx} efficiently by solving the LP \eqref{eq_lp}. The second part of Lemma \ref{l_lp} states that the optimal objective function value of the LP \eqref{eq_lp} is Lipschitz continuous in the RHS parameter $\bm{b}\in\R_+^{a_2}$. Because we will solve LP \eqref{eq_lp} regularly over time (at discrete time epochs) and LP \eqref{eq_lp} may have multiple optimal solutions at some time epochs, we need to choose an optimal solution among the set of optimal solutions at those time epochs such that the optimal solutions that we will use over time will not fluctuate a lot. The third part of Lemma \ref{l_lp} presents a method to achieve the aforementioned goal. For given $\bm{b}\in\R_+^{a_2}$, QP \eqref{eq_qp} finds the optimal solution of the LP \eqref{eq_lp} with the smallest Euclidean norm. Because QP \eqref{eq_qp} is convex, it is solvable in polynomial time (see \citet{vav08}). The third part of Lemma \ref{l_lp} states that the optimal solution of the QP \eqref{eq_qp} is unique and Lipschitz continuous in the RHS parameter $\bm{b}\in\R_+^{a_2}$. A direct consequence of the third part of Lemma \ref{l_lp} is the following Lipschitz continuity result.

\begin{lemma}\label{l_reg}
For any given nonnegative parameter process $(\bm{b}(t),t\in\R_+)\in\D^{a_2}$, let $(\ot{y}_{j}(t),\ot{q}_{ij}(t),j\in\mJ,i\in\mI_j,t\in\R_+)\in\D^{a_3}$ denote the optimal solution process associated with the LP \eqref{eq_lp} selected by the QP \eqref{eq_qp}. For all $s,t\in\R_+$, we have
\begin{equation*}
\max_{j\in\mJ,i\in\mI_j} \left| \ot{q}_{ij}(s)-\ot{q}_{ij}(t)\right| \leq C_2 \big| \bm{b}(s)-\bm{b}(t) \big|_{\infty}.
\end{equation*}
\end{lemma}

\subsection{Asymptotic Lower Bound with respect to the Objective \eqref{eq_obj_2}}\label{s_alb2}

We prove that the optimal objective function value of the DCP \eqref{eq_dcp} provides an asymptotic lower bound on the performance of any admissible policy with respect to the objective \eqref{eq_obj_2}.

\begin{theorem}\label{t_lb}
Let $\pi=\{\pi^r,r\in\N_+\}$ be an arbitrary sequence of admissible policies. Then, for all $\ep>0$ and $t \in\R_+$, we have
\begin{equation*}
\liminf_{r\rightarrow\infty} \pr\left( \sum_{j\in\mJ} h_j \hat{N}_j^{\pi,r}(t) > \ep\right) \geq \pr\left( z\Big( \ot{Q}_k(t),j\in\mJ,k\in\mK_j,\ot{W}_i(t),i\in\mI\Big) > \ep\right),
\end{equation*}
where $\Big( \ot{Q}_k,j\in\mJ,k\in\mK_j,\ot{W}_i,i\in\mI\Big)$ is defined in Proposition \ref{p_diffusion}.
\end{theorem}
The proof of Theorem \ref{t_lb} is presented in Appendix \ref{t_lb_proof}. 

We call a sequence of admissible policies \textit{asymptotically optimal} with respect to the objective \eqref{eq_obj_2}, if it achieves the asymptotic lower bound in Theorem \ref{t_lb}. Next, we will formally introduce the proposed policy. Then, we will prove that the proposed policy is asymptotically optimal.

\section{Proposed Policy}\label{s_policy}

By Proposition \ref{p_diffusion}, Lemma \ref{l_lp}, and Theorem \ref{t_lb}, if an admissible policy keeps the diffusion scaled number of jobs in the buffers in front of the shared servers close to an optimal LP \eqref{eq_lp} solution under the LP parameters $(\hat{Q}_k^r(t),j\in\mJ,k\in\mK_j, \hat{W}_i^r(t),i\in\mI)$ at all times for sufficiently large $r$, then that policy is a good candidate for an asymptotically optimal policy. Therefore, the policy that we will propose should track the optimal LP \eqref{eq_lp} solution at all times. Specifically, at each shared server, we will compare the number of jobs in front of that shared server with the optimal LP \eqref{eq_lp} solution selected by the QP \eqref{eq_qp}, and then determine a scheduling rule in the shared server which makes the number of jobs in front of that shared server close to that optimal LP \eqref{eq_lp} solution. Then, we will resolve the LP \eqref{eq_lp} and then the QP \eqref{eq_qp} and repeat the same procedure. We call the time between successively solving the LP \eqref{eq_lp} for a shared server as the \textit{review period} for that shared server. At each review period, the shared server takes action in order to make the numbers of the job types that it processes close to the optimal LP \eqref{eq_lp} solution.

First, we will introduce some additional notation below. Then, we will explain the intuition behind our proposed policy in Section \ref{s_intuition}. Next, we will formally introduce the proposed policy in Section \ref{s_p_def}. Finally, we will prove the asymptotic optimality of the proposed policy in Section \ref{s_ao}.

Let us fix an arbitrary $r\in\N_+$ and a sample path. Let $(y_j^{*,r},q_{ij}^{*,r},j\in\mJ,i\in\mI_j)\in\D^{a_3}$ denote the optimal solution process of the LP \eqref{eq_lp} under the parameters $(Q_k^r,j\in\mJ,k\in\mK_j, W_i^r,i\in\mI)\in\D^{a_2}$ selected by the QP \eqref{eq_qp}. By \eqref{eq_workload} and \eqref{eq_lp_4}, we have
\begin{equation}\label{eq_we}
\sum_{j\in\mJ_i} \frac{Q_{ij}^r(t)}{\mu_{ij}}= \sum_{j\in\mJ_i} \frac{q_{ij}^{*,r}(t)}{\mu_{ij}}=W_i^r(t),\qquad\forall i\in\mI,\;t\in\R_+.
\end{equation}
For all $i\in\mI$ and $t\in\R_+$, let 
\begin{align*}
\mJ_i^{>,r}(t)&:=\{j\in\mJ_i: Q_{ij}^r(t)>\lceil q_{ij}^{*,r}(t) \rceil \},\\
\mJ_i^{\leq,r}(t)&:=\{j\in\mJ_i: Q_{ij}^r(t)\leq \lceil q_{ij}^{*,r}(t) \rceil \}.
\end{align*}
Then, $\{\mJ_i^{>,r}(t),\mJ_i^{\leq,r}(t)\}$ is a disjoint partition of $\mJ_i$ for all $i\in\mI$ and $t\in\R_+$.

\subsection{Intuition Behind the Proposed Policy}\label{s_intuition}

In this section, by non-rigorous arguments, we derive some intuition for the control policy that we will propose. Let us consider an arbitrary shared server $i\in\mI$ at an arbitrary time $t\in\R_+$. For simplicity, let us assume that $q_{ij}^{*,r}(t)$ is an integer for all $j\in\mJ_i$. Suppose that there exists a $j\in\mJ_i$ such that $Q_{ij}^r(t)\neq q_{ij}^{*,r}(t)$. Then, $\mJ_i^{>,r}(t)\neq \emptyset$ and $\mJ_i^{\leq,r}(t)\neq \emptyset$ by \eqref{eq_we}. We want the shared server $i$ to decrease the number of jobs in the buffers associated with $j\in\mJ_i^{>,r}(t)$ from $Q_{ij}^r(t)$ to $q_{ij}^{*,r}(t)$, while keeping the number of jobs in the buffers associated with $j\in\mJ_i^{\leq,r}(t)$ less than or equal to $q_{ij}^{*,r}(t)$. Let $L_i^r(t)\in\R_+$ denote the expected length of the review period for given $(Q_{ij}^r(t),j\in\mJ)$. Then, $L_i^r(t)$ should satisfy the equalities
\begin{align}
L_i^r(t) &= \sum_{j\in\mJ_i^{>,r}(t)} \frac{Q_{ij}^r(t)-q_{ij}^{*,r}(t)}{\mu_{ij}} + L_i^r(t)\sum_{j\in\mJ_i^{>,r}(t)} \frac{\lambda_j^r}{\mu_{ij}} + \sum_{j\in\mJ_i^{\leq,r}(t)} \frac{\left(\lambda_j^r L_i^r(t) - q_{ij}^{*,r}(t)+Q_{ij}^r(t)\right)^+}{\mu_{ij}}\label{eq_review}\\
&= \sum_{j\in\mJ_i}\frac{\left(\lambda_j^r L_i^r(t) -q_{ij}^{*,r}(t)+Q_{ij}^r(t)\right)^+}{\mu_{ij}}.\label{eq_l_review_3}
\end{align}
Notice that \eqref{eq_l_review_3} is a compact version of the RHS of \eqref{eq_review}. The first term in the RHS of \eqref{eq_review} denotes the average time that the shared server $i$ should spend to deplete the excess jobs in the set $j\in\mJ_i^{>,r}(t)$. In the mean time, there will be external type $j$ job arrivals for all $j\in\mJ_i$. Hence, the second term in the RHS of \eqref{eq_review} denotes the average time that the shared server $i$ should spend to process the excess jobs due to the external job arrivals associated with the jobs in the set $\mJ_i^{>,r}(t)$. Finally, the third term in the RHS of \eqref{eq_review} denotes the average time that the shared server $i$ should spend to process the jobs in the set $\mJ_i^{\leq,r}(t)$ if the average number of external job arrivals associated with the job type $j\in\mJ_i^{\leq,r}(t)$ is greater than $q_{ij}^{*,r}(t) - Q_{ij}^r(t)$. Then, we have the following result.

\begin{lemma}\label{l_review}
If $\lambda_j^r=\lambda_j$ for all $j\in\mJ$, that is, if the arrival rates are equal to the limiting ones, then, for all $i\in\mI$, $L_i^r(t)\in\R_+$ is a solution of the equality \eqref{eq_review} if and only if
\begin{equation}\label{eq_l_review}
L_i^r(t) \geq \max_{j\in\mJ_i^{\leq,r}(t)} \frac{q_{ij}^{*,r}(t)-Q_{ij}^r(t)}{\lambda_j}.
\end{equation}
\end{lemma}
The proof of Lemma \ref{l_review} is presented in Appendix \ref{l_review_proof}. Lemma \ref{l_review} provides a lower bound on the expected length of the review period under the limiting arrival rates. However, we do not want the length of the review period to be very long because otherwise at the end of the review period, the system state can be far away from the optimal LP \eqref{eq_lp} solution. Therefore, intuitively, it is better to have the expected length of the review period as short as possible. Hence, we choose 
\begin{equation}\label{eq_l_review_2}
L_i^r(t) = \max_{j\in\mJ_i^{\leq,r}(t)} \frac{q_{ij}^{*,r}(t)-Q_{ij}^r(t)}{\lambda_j}. 
\end{equation}

By \eqref{eq_l_review_3} and \eqref{eq_l_review_2}, under the limiting arrival rates, the shared server $i$ does not allocate any time during the review period for the job types in the set
\begin{equation*}
\mJ_i^{<,r}(t):= \argmax_{j\in\mJ_i^{\leq,r}(t)} \frac{q_{ij}^{*,r}(t)-Q_{ij}^r(t)}{\lambda_j}.
\end{equation*}
This is because the length of the review period is short enough such that for all $j\in\mJ_i^{<,r}(t)$, the number of external job arrivals to the buffer $ij$ will not make the number of jobs waiting in that buffer greater than $q_{ij}^{*,r}(t)$ at the end of the review period.

Suppose that $(Q_{ij}^r,j\in\mJ_i)$ is close to $(q_{ij}^*,j\in\mJ_i)$ at the beginning of the review period. Then, the length of the review period will be short by \eqref{eq_l_review_2}. Hence, we expect the process $(Q_k^r,j\in\mJ,k\in\mK_j, W_i^r,i\in\mI)$ to not to change significantly during the review period. By Lemma \ref{l_reg}, the optimal LP \eqref{eq_lp} solution will not change significantly during the review period. Hence, we expect the number of jobs in the buffers in front of the shared server $i$ to be close to the optimal LP \eqref{eq_lp} solution at the end of the review period too. Consequently, we expect $(Q_{ij}^r(t),j\in\mJ_i)$ to be close to $(q_{ij}^{*,r}(t),j\in\mJ_i)$ for all $t\in\R_+$. If we repeat this procedure at each shared server, then we expect to achieve the asymptotic lower bound in Theorem \ref{t_lb}.

Based on this intuition, we formally propose a control policy in the following section.

\subsection{Formal Definition of the Proposed Policy}\label{s_p_def}

We propose a continuous-review, state dependent, and non-preemptive control policy.

\begin{definition}\label{d_policy}
For all $i\in\mI$, the proposed policy for the shared server $i$ is the following:
\begin{description}[align=left]
\item [Step 0] (\textit{Initialization}) Go to Step 1.

\item [Step 1] Let $t\in\R_+$ denote the current time. Solve the LP \eqref{eq_lp} and then the QP \eqref{eq_qp}. If $\mJ_i^{>,r}(t)=\emptyset$, go to Step 2. Otherwise, go to Step 3.

\item [Step 2] Let $t\in\R_+$ denote the current time. If there are not any jobs waiting in front of the shared server $i$ at time $t$, then the server processes the first job that externally arrives after time $t$. Otherwise, the shared server $i$ processes an arbitrary job among the jobs waiting at the head of the buffers $\{ij: j\in\mJ_i,\; Q_{ij}^r(t)>0\}$. At the first service completion epoch in the shared server $i$ after time $t$, go to Step 1.

\item [Step 3] Let $t\in\R_+$ denote the current time. Because $\mJ_i^{>,r}(t)\neq\emptyset$, there exists a $j\in\mJ_i$ such that $Q_{ij}^r(t)>\lceil q_{ij}^{*,r}(t) \rceil\geq q_{ij}^{*,r}(t)$. This implies that there exists an $l\in\mJ_i$ such that $Q_{il}^r(t)< q_{il}^{*,r}(t)\leq \ru{q_{il}^{*,r}(t)}$ by \eqref{eq_we}. Hence, $\mJ_i^{\leq,r}(t)\neq \emptyset$ and $\mJ_i^{<,r}(t)\neq \emptyset$ by definition. Let us choose an arbitrary $m\in\mJ_i^{<,r}(t)$. The shared server $i$ first processes the excess $Q_{ij}^r(t)- \lceil q_{ij}^{*,r}(t) \rceil$ jobs in the buffers associated with the job types in $\mJ_i^{>,r}(t)$ in an admissible and non-preemptive way. Let $t_1\geq t$ denote the first time when those excess jobs are processed. During the interval $[t,t_1]$, if there are external job arrivals such that $Q_{ij}^r(t_1)> \lceil q_{ij}^{*,r}(t) \rceil$ for some $j\in\mJ_i\backslash\{m\}$, then the shared server $i$ should process those excess $Q_{ij}^r(t_1)- \lceil q_{ij}^{*,r}(t) \rceil$ jobs in an admissible and non-preemptive way. Let $t_2\geq t_1$ denote the first time when those excess jobs are processed. During the interval $[t_1,t_2]$, if there are external job arrivals such that $Q_{ij}^r(t_2)> \lceil q_{ij}^{*,r}(t) \rceil$ for some $j\in\mJ_i\backslash\{m\}$, then the shared server $i$ should process those excess $Q_{ij}^r(t_2)- \lceil q_{ij}^{*,r}(t) \rceil$ jobs in an admissible and non-preemptive way. The shared server $i$ continues processing the jobs in the same way until
\begin{align*}
&\inf\left\{s\geq t: Q_{ij}^r(s)=\lceil q_{ij}^{*,r}(t) \rceil,\forall j\in\mJ_i^{>,r}(t),\quad Q_{ij}^r(s)\leq \lceil q_{ij}^{*,r}(t) \rceil,\forall j\in\mJ_i^{\leq,r}(t)\backslash\{m\}\right\}\\
&\hspace{3cm} = \inf\Bigg\{s\geq t: \sum_{j\in\mJ_i\backslash\{m\}} \sum_{x=1}^{\left(Q_{ij}^r(t)- \lceil q_{ij}^{*,r}(t) \rceil+A_j^r(s)-A_j^r(t)\right)^+} v_{ij(S_{ij}(T_{ij}^r(t))+x)} \leq s-t \Bigg\}.
\end{align*}
At time $s$, go to Step 1.

\end{description}
\end{definition}

The proposed policy is the simultaneous implementation of the control policy defined in Definition \ref{d_policy} in all of the shared servers. Observe that Steps 0 and 1 are done instantaneously and both Step 2 and Step 3 are review periods for the shared server $i$. By definition, Step 2 lasts at most as much as the sum of a residual inter-arrival time and a service time. Hence, Step 2 does not last very long (specifically, we will prove that the length of Step 2 is $o_p(r)$ in Lemma \ref{l_tail_2}, where $o_p(\cdot)$ denotes the little-$o$ in probability).

In Step 3, the shared server $i$ works on at most $|\mJ_i|-1$ number of type of jobs. Hence, it acts like a light traffic queue by Assumption \ref{a_ht} Part 2. Therefore, given that the system state is not very far away from the optimal LP \eqref{eq_lp} solution at the beginning of Step 3, the shared server $i$ quickly completes Step 3. Because Step 3 does not last very long (see Lemma \ref{l_tail_3}), the number of type $m$ jobs ($m\in\mJ_i^{<,r}(t)$), that is, the number of the job type that the shared server $i$ does not process in Step 3, will not grow significantly. Consequently, at the end of Step 2 or 3, the number of jobs in front of the shared server $i$ will be close to the optimal LP \eqref{eq_lp} solution.

\subsection{Asymptotic Optimality of the Proposed Policy}\label{s_ao}

In this section, we prove that the proposed policy is asymptotically optimal with respect to the objective \eqref{eq_obj_2}. Then, we show that this result implies asymptotic optimality with respect to the objective \eqref{eq_obj_1}.

\begin{theorem}\label{t_ao_1}
Consider the proposed policy defined in Definition \ref{d_policy}. Then, for all $\ep>0$ and $t\in\R_+$, we have
\begin{equation*}
\lim_{r\rightarrow\infty} \pr\left( \sum_{j\in\mJ} h_j \hat{N}_j^{r}(t) > \ep\right) = \pr\left( z\Big( \ot{Q}_k(t),j\in\mJ,k\in\mK_j,\ot{W}_i(t),i\in\mI\Big) > \ep\right),
\end{equation*}
where $\Big( \ot{Q}_k,j\in\mJ,k\in\mK_j,\ot{W}_i,i\in\mI\Big)$ is defined in Proposition \ref{p_diffusion}.
\end{theorem}
The proof of Theorem \ref{t_ao_1} is presented in Appendix \ref{t_ao_1_proof}. Theorem \ref{t_ao_1} states that the proposed policy achieves the asymptotic lower bound in Theorem \ref{t_lb}, thus it is asymptotically optimal with respect to the objective \eqref{eq_obj_2}. This result also implies asymptotic optimality with respect to the objective \eqref{eq_obj_1} as formally stated below.

\begin{theorem}\label{t_ao_2}
Let $\pi = \{\pi^r,r\in\N_+\}$ be an arbitrary sequence of admissible policies and $*$ denote the proposed policy. Then,
\begin{align*}
\lim_{r\rightarrow\infty}\sum_{j\in\mJ} h_j  \E\left[ \int_0^\infty \e^{-\delta t} \hat{N}_j^{*,r}(t) \dr t\right] &= \sum_{j\in\mJ} h_j \E\left[  \int_0^\infty \e^{-\delta t} z\Big( \ot{Q}_k(t),j\in\mJ,k\in\mK_j,\ot{W}_i(t),i\in\mI\Big)\dr t\right] \\
&\leq \liminf_{r\rightarrow\infty} \sum_{j\in\mJ} h_j  \E\left[ \int_0^\infty \e^{-\delta t} \hat{N}_j^{\pi,r}(t) \dr t\right],
\end{align*}
where $\Big( \ot{Q}_k,j\in\mJ,k\in\mK_j,\ot{W}_i,i\in\mI\Big)$ is defined in Proposition \ref{p_diffusion}.
\end{theorem}
The proof of Theorem \ref{t_ao_2} follows from Theorems \ref{t_lb} and \ref{t_ao_1} and a uniform integrability result and is very similar to the proof of Theorem 3 in \citet{ozk19} and the proof of Theorem 5.3 in \citet{bel01}. Hence, we skip it.

\section{Extensions}\label{s_ext}

We extend the empty initial system assumption in Section \ref{s_ext_i}, instantaneous fork and join operations assumption in Section \ref{s_ext_ni}, and the network structure in Section \ref{s_ext_ns}.

\subsection{Non-Empty Initial System}\label{s_ext_i}

We extend the empty initial system assumption with the following one:
\begin{assumption}\label{a_initial} 
For all $r\in\N_+$, $\bm{Q}^r(0):= \big(Q_k^r(0),\; Q_{ij}^{r}(0),\; Q_k^{(1),r}(0), \; Q_{ij}^{(1),r}(0),\; j\in\mJ, k\in\mK_j,i\in\mI_j\big)$ is a random vector independent of all other stochastic primitives and takes values in $\N^{a_4}$, where $a_4:= 2\sum_{j\in\mJ} \left(K_j + |\mI_j|\right)$. Furthermore, 
\begin{enumerate}
\item $r^{-2}\bm{Q}^r(0)\xrightarrow{a.s.} 0^{(a_4)}$ and $r^{-1}\bm{Q}^r(0)\Rightarrow \ot{\bm{Q}}(0)$ such that $\ot{Q}_k(0)= 0$ for all $j\in\mJ$ and $k\in\mK_j^L$.

\item There exists an $n_1\in\N_+$ such that
\begin{equation*}
\sup_{r\geq n_1} \E\left[ \left(r^{-1}\bm{Q}^r(0)\right)^2\right] \in\R^{a_4}_+.
\end{equation*}

\item For all $\ep>0$,
\begin{equation*}
\pr\left(\max_{j\in\mJ,i\in\mI_j} \left| Q_{ij}^{r}(0) - q_{ij}^{*,r}(0)\right| >\ep r \right) \rightarrow 0.
\end{equation*}

\item For all $\ep>0$, there exist $n_2(\ep)\in\N_+$ such that if $r\geq n_2(\ep)$,
\begin{equation*}
\pr\left(Q_k^r(0) >\frac{((\mu_k-\lambda_j^r)\wedge\ep) r}{7} \right) \leq C_3 r^5\e^{-C_4 r}\qquad\forall j\in\mJ,\;k\in\{k\in\mK_j: \lambda_j<\mu_k\},
\end{equation*}
where $C_3$ and $C_4$ are strictly positive constants independent of $r$.

\end{enumerate}
\end{assumption}
We need Assumption \ref{a_initial} Part 1 to prove Propositions \ref{p_fluid} and \ref{p_diffusion}. Assumption \ref{a_initial} Part 2 is a uniform integrability condition which is used to prove Theorem \ref{t_ao_2}. We need Assumption \ref{a_initial} Part 3 to prove Lemma \ref{l_good_set}. Finally, we need Assumption \ref{a_initial} Part 4 to prove Lemma \ref{l_tail_4}.

\subsection{Non-Instantaneous Fork and Join Operations}\label{s_ext_ni}

So far, we assume that the fork and join operations are done instantaneously. However, we can extend this assumption in the following way. Suppose that for all $j\in\mJ$, there exist a fork server and a join server which make the fork and join operations for the type $j$ jobs, respectively, and there exists an infinite capacity buffer in front of the fork server (see for example the networks in Figures \ref{fj_network_0} and \ref{fj_network_1}). Furthermore, the fork server can be in either heavy or light traffic but the join server must be in light traffic. Then, all of our results hold under this extension (see \citet{ozk19} for an explicit and rigorous extension). 

It is crucial for the join servers to be in light traffic because otherwise there will be workload in front of the join servers because of not only the synchronization constraints but also the tight processing capacity. Hence, we will have workload constraints associated with the join servers in the DCP \eqref{eq_dcp}. However, those workload processes depend on the scheduling control in the shared servers nonlinearly. Consequently, the resulting DCP will be very complicated and it is not clear how to solve that DCP and interpret a control policy from it. An interesting and challenging future research topic is to derive an asymptotically optimal control policy when some of the join servers are in heavy traffic.

\subsection{Extensions of the Network Structure}\label{s_ext_ns}

Consider an arbitrary job type $j\in\mJ$ and a dedicated server $k\in\mK_j$. We can replace the dedicated server $k$ and the buffer in front of it with an arbitrary open queueing network with private servers and no control. Let $Q_k^r(t)$ denote the total number of jobs at time $t\in\R_+$ in that queueing network. As long as Proposition \ref{p_diffusion} can be extended with the weak convergence of the process $\hat{Q}_k^r$ and Lemma \ref{l_tail_4} can be extended by including the process $Q_k^r$, all of the results in the paper continue to hold under this extension.

Next, let us consider an arbitrary job type $j\in\mJ$ and a shared server $i\in\mI_j$. We can insert an arbitrary open queueing network with private servers and no control between the fork operation of type $j$ jobs and the shared server $i$. Let $Q_{ij}^{(2),r}(t)$ denote the total number of jobs at time $t\in\R_+$ in that queueing network. As long as Proposition \ref{p_diffusion} can be extended with the weak convergence of the process $\hat{Q}_{ij}^{(2),r}$, Lemma \ref{l_tail_1} (specifically \eqref{eq_t1_1}) can be extended with the departure process from the aforementioned queueing network, and Lemma \ref{l_tail_4} can be extended by including the process $Q_{ij}^{(2),r}$, all of the results in the paper continue to hold. The only difference is that the constraint \eqref{eq_lp_3} of the LP \eqref{eq_lp} and the constraint \eqref{eq_qp_3} of the QP \eqref{eq_qp} should be modified as
\begin{equation*}
y_j- q_{ij}\geq q_{ij}^{(2)}\qquad\forall j\in\mJ,\;i\in\mI_j,
\end{equation*}
where $q_{ij}^{(2)}$ is a parameter associated with $Q_{ij}^{(2),r}$.

The complicated case is when there are heavy-traffic queues after the shared servers. By a similar argument presented in Section \ref{s_ext_ni}, it is not clear either what the proposed policy should be or how to prove an asymptotic optimality result in that case. An excellent topic for future research is to develop control policies for the broader class of fork-join networks with multiple job types described in \citet{ngu94}. More specifically, that paper assumes FCFS scheduling, but we believe other control policies can lead to better performance.

%
%
%

\begin{APPENDICES}


\section{Lemma Proofs}\label{s_l_proofs}

Sections \ref{l_lp_proof} and \ref{l_review_proof} present the proofs of Lemmas \ref{l_lp} and \ref{l_review}, respectively.

\subsection{Proof of Lemma \ref{l_lp}}\label{l_lp_proof}

Notice that there exists an optimal solution of the LP \eqref{eq_lp} for all $\bm{b}\in\R_+^{a_2}$. Let $\big( q_{ij},j\in\mJ,i\in\mI_j\big)$ be an arbitrary feasible point of the convex problem \eqref{eq_cvx} and let us define $y_j:=\max_{k\in\mK_j} q_k \vee  \max_{i\in\mI_j} q_{ij}$ for all $j\in\mJ$. Then, $\big( y_j, q_{ij},j\in\mJ,i\in\mI_j\big)$ is a feasible point of the LP \eqref{eq_lp} with the same objective function value. Therefore, for all $\big( q_{ij},j\in\mJ,i\in\mI_j\big)$ which is a feasible point of the convex problem \eqref{eq_cvx}, we have
\begin{equation}\label{eq_l_lp_1}
\sum_{j\in\mJ} h_j \left(\max_{k\in\mK_j} q_k\vee \max_{i\in\mI_j} q_{ij}\right) \geq z(\bm{b}).
\end{equation}
In other words, the optimal objective function value of the LP \eqref{eq_lp} is a lower bound on the objective function value of any feasible point of the convex problem \eqref{eq_cvx}. 

Let $\big( \ot{y}_j,\ot{q}_{ij},j\in\mJ,i\in\mI_j\big)$ be an arbitrary optimal solution of the LP \eqref{eq_lp}. By \eqref{eq_lp_2} and \eqref{eq_lp_3}, $\ot{y}_j\geq \max_{k\in\mK_j} q_k\vee \max_{i\in\mI_j} \ot{q}_{ij}$ for all $j\in\mJ$ and so we can choose $\ot{y}_j= \max_{k\in\mK_j} q_k\vee \max_{i\in\mI_j} \ot{q}_{ij}$ for all $j\in\mJ$ without loss of generality by \eqref{eq_lp_1}. Notice that, $\big(\ot{q}_{ij},j\in\mJ,i\in\mI_j\big)$ is a feasible point of the convex problem \eqref{eq_cvx} with the objective function value $\sum_{j\in\mJ} h_j \left(\max_{k\in\mK_j} q_k\vee \max_{i\in\mI_j} \ot{q}_{ij}\right)$. By \eqref{eq_lp_2} and \eqref{eq_lp_3}, we have 
\begin{equation}\label{eq_l_lp_2}
z(\bm{b}) = \sum_{j\in\mJ} h_j \ot{y}_j = \sum_{j\in\mJ} h_j \left(\max_{k\in\mK_j} q_k\vee \max_{i\in\mI_j} \ot{q}_{ij}\right).
\end{equation}
Therefore, $\big(\ot{q}_{ij},j\in\mJ,i\in\mI_j\big)$ is an optimal solution of the convex problem \eqref{eq_cvx} with the objective function value $z(\bm{b}) $ by \eqref{eq_l_lp_1} and \eqref{eq_l_lp_2}. 

The second part of Lemma \ref{l_lp} follows directly from Equation (10.22) of \citet{sch98}. Finally, the third part of Lemma \ref{l_lp} follows directly from Proposition 4.1.d of \citet{han12}.

\subsection{Proof of Lemma \ref{l_review}}\label{l_review_proof}

By \eqref{eq_we}, we have
\begin{equation}\label{eq_l_review_p_1}
\sum_{j\in\mJ_i^{>,r}(t)} \frac{Q_{ij}^r(t)-q_{ij}^{*,r}(t)}{\mu_{ij}}= \sum_{j\in\mJ_i^{\leq,r}(t)} \frac{q_{ij}^{*,r}(t)-Q_{ij}^r(t)}{\mu_{ij}},\qquad\forall i\in\mI,\;t\in\R_+.
\end{equation}
By Assumption \ref{a_ht} Part 2 and \eqref{eq_l_review_p_1}, the RHS of \eqref{eq_review} is equal to
\begin{align*}
&\sum_{j\in\mJ_i^{\leq,r}(t)} \frac{q_{ij}^{*,r}(t)-Q_{ij}^r(t)}{\mu_{ij}} +L_i^r(t)\left(1-\sum_{j\in\mJ_i^{\leq,r}(t)} \frac{\lambda_j}{\mu_{ij}}\right)  + \sum_{j\in\mJ_i^{\leq,r}(t)} \frac{\left(\lambda_j L_i^r(t) - q_{ij}^{*,r}(t)+Q_{ij}^r(t)\right)^+}{\mu_{ij}}\\
&\hspace{3cm}= L_i^r(t)+ \sum_{j\in\mJ_i^{\leq,r}(t)} \frac{q_{ij}^{*,r}(t)-Q_{ij}^r(t)-\lambda_jL_i^r(t)+ \left(\lambda_j L_i^r(t) - q_{ij}^{*,r}(t)+Q_{ij}^r(t)\right)^+}{\mu_{ij}}\\
&\hspace{3cm}= L_i^r(t)+ \sum_{j\in\mJ_i^{\leq,r}(t)} \frac{\left(q_{ij}^{*,r}(t)-Q_{ij}^r(t)-\lambda_jL_i^r(t)\right)^+}{\mu_{ij}}.
\end{align*}
Therefore, in order for \eqref{eq_review} to hold, it must be the case that $q_{ij}^{*,r}(t)-Q_{ij}^r(t)-\lambda_jL_i^r(t)\leq 0$ for all $j\in\mJ_i^{\leq,r}(t)$ which holds if and only if \eqref{eq_l_review} holds.

\section{Theorem Proofs}\label{s_t_proofs}

Sections \ref{t_lb_proof} and \ref{t_ao_1_proof} present the proofs of Theorems \ref{t_lb} and \ref{t_ao_1}, respectively.

\subsection{Proof of Theorem \ref{t_lb}}\label{t_lb_proof}

Let us fix an arbitrary sequence of admissible policies $\pi=\{\pi^r,r\in\N_+\}$ and arbitrary $\ep>0$ and $t \in\R_+$. By \eqref{eq_queue_3}, \eqref{eq_diffusion}, and \eqref{eq_workload},
\begin{equation*}
\sum_{j\in\mJ_i} \frac{\hat{Q}_{ij}^{\pi,r}(t)}{\mu_{ij}} = \hat{W}_i^{\pi,r}(t)\quad\forall i\in\mI,\hspace{2cm}\hat{Q}_{ij}^{\pi,r}(t)\geq 0\quad\forall j\in\mJ,i\in\mI_j.
\end{equation*}
Therefore, $\big(\hat{Q}_{ij}^{\pi,r}(t), j\in\mJ, i\in\mI_j\big)$ is a feasible point of the convex problem \eqref{eq_cvx} under the parameters $\big( \hat{Q}_k^{r}(t),j\in\mJ,k\in\mK_j,\hat{W}_i^{\pi,r}(t),i\in\mI\big)$. By Lemma \ref{l_lp} Part 1, we have
\begin{equation}\label{eq_t_lb_1}
\sum_{j\in\mJ} h_j \left(\max_{k\in\mK_j} \hat{Q}_k^{r}(t) \vee  \max_{i\in\mI_j} \hat{Q}_{ij}^{\pi,r}(t)\right) \geq  z\Big( \hat{Q}_k^{r}(t),j\in\mJ,k\in\mK_j,\hat{W}_i^{\pi,r}(t),i\in\mI\Big),
\end{equation}
which holds for all sample paths. Then,
\begin{align}
\liminf_{r\rightarrow\infty}\pr\left( \sum_{j\in\mJ} h_j \hat{N}_j^{\pi,r} (t) > \ep\right) &= \liminf_{r\rightarrow\infty}\pr\left( \sum_{j\in\mJ} h_j \left(\max_{k\in\mK_j} \hat{Q}_k^{r}(t) \vee  \max_{i\in\mI_j} \hat{Q}_{ij}^{\pi,r}(t)\right) > \ep\right) \label{eq_t_lb_2}\\
&\geq \liminf_{r\rightarrow\infty}\pr\left( z\Big( \hat{Q}_k^{r}(t),j\in\mJ,k\in\mK_j,\hat{W}_i^{\pi,r}(t),i\in\mI\Big)> \ep\right)\label{eq_t_lb_3}\\
&=\pr\left( z\Big( \ot{Q}_k(t),j\in\mJ,k\in\mK_j,\ot{W}_i(t),i\in\mI\Big)> \ep\right),\label{eq_t_lb_4}
\end{align}
where \eqref{eq_t_lb_2} is by \eqref{eq_obj_3}, \eqref{eq_t_lb_3} is by \eqref{eq_t_lb_1}, and \eqref{eq_t_lb_4} is by Proposition \ref{p_diffusion}, Lemma \ref{l_lp}, and Theorems 3.4.3 and 11.6.6 of \citet{whi02}. Specifically, weak convergence result in Proposition \ref{p_diffusion} implies weak convergence of the associated finite dimensional distributions by Theorem 11.6.6 of \citet{whi02}. Because $z$ is Lipschitz continuous (see Lemma \ref{l_lp} Part 2), the convergence result in \eqref{eq_t_lb_4} follows from continuous mapping theorem (see Theorem 3.4.3 of \citet{whi02}).

\subsection{Proof of Theorem \ref{t_ao_1}}\label{t_ao_1_proof}

Let $Z$ be a mapping from $\D^{a_2}$ such that $Z(f)(t):=z(f(t))$ for all $f\in\D^{a_2}$ and $t\in\R_+$. Then, $Z$ is the process version of $z$. Since $z$ is Lipschitz continuous (see Lemma \ref{l_lp} Part 2), $Z$ maps the functions from $\D^{a_2}$ to $\D$, that is, $Z:\D^{a_2}\rightarrow\D$. Let $d(\cdot)$ denote the Skorokhod distance (see Equation (12.13) of \citet{bil99}). For arbitrary $X,Y\in\D^{a_2}$, because $z$ is Lipschitz continuous (see Lemma \ref{l_lp} Part 2), one can see that $d(Z(X),Z(Y))\leq (C_1\vee 1)d(X,Y)$. Therefore, $Z$ is also Lipschitz continuous. By Proposition \ref{p_diffusion} and continuous mapping theorem (see Theorem 3.4.3 of \citet{whi02}), we have
\begin{equation*}
Z\left( \hat{Q}_k^r, j\in\mJ,k\in\mK_j,\; \hat{W}_i^r,i\in\mI \right)\Rightarrow Z\left( \ot{Q}_k, j\in\mJ,k\in\mK_j,\; \ot{W}_i,i\in\mI \right).
\end{equation*}
Moreover, we have the following proposition whose proof is presented in Section \ref{p_main_proof}.

\begin{proposition}\label{p_main}
Let us fix arbitrary $\ep,T>0$. Under the proposed policy (see Definition \ref{d_policy}),
\begin{equation*}
\lim_{r\rightarrow\infty}\pr\left( \bigg\Vert \sum_{j\in\mJ} h_j \hat{N}_j^{r} - Z\left( \hat{Q}_k^r, j\in\mJ,k\in\mK_j,\; \hat{W}_i^r,i\in\mI \right)\bigg\Vert_{T}>\ep\right)=0.
\end{equation*}
\end{proposition}
By Proposition \ref{p_main} and convergence-together theorem (see Theorem 11.4.7 of \citet{whi02}), we have the following weak convergence result associated with the proposed policy:
\begin{equation}\label{eq_wcp}
\sum_{j\in\mJ} h_j \hat{N}_j^{r} \Rightarrow Z\left( \ot{Q}_k, j\in\mJ,k\in\mK_j,\; \ot{W}_i,i\in\mI \right).
\end{equation}
Finally, Theorem \ref{t_ao_1} follows from \eqref{eq_wcp}.

\section{Proof of Proposition \ref{p_main}}\label{p_main_proof}

Let us fix an arbitrary $\ep,T>0$. Let $\hat{q}_{ij}^{*,r}(t):=r^{-1}q_{ij}^{*,r}(r^2t)$ denote the diffusion scaled version of the optimal solution process for all $j\in\mJ$, $i\in\mI_j$, $t\in\R_+$, and $r\in\N_+$. By \eqref{eq_obj_3} and Lemma \ref{l_lp} Part 1, the probability in Proposition \ref{p_main} is equal to
\begin{align}
&\pr\left( \Bigg\Vert \sum_{j\in\mJ} h_j \left(\max_{k\in\mK_j} \hat{Q}_k^r \vee  \max_{i\in\mI_j} \hat{Q}_{ij}^r - \max_{k\in\mK_j} \hat{Q}_k^r \vee  \max_{i\in\mI_j} \hat{q}_{ij}^{*,r}\right)\Bigg\Vert_{T}>\ep\right)\nonumber\\
&\hspace{3cm}\leq \pr\left( \sum_{j\in\mJ}h_j \Bigg\Vert  \left(\max_{k\in\mK_j} \hat{Q}_k^r \vee  \max_{i\in\mI_j} \hat{Q}_{ij}^r - \max_{k\in\mK_j} \hat{Q}_k^r \vee  \max_{i\in\mI_j} \hat{q}_{ij}^{*,r}\right)\Bigg\Vert_{T}>\ep\right)\nonumber\\
&\hspace{3cm}\leq \pr\left( \sum_{j\in\mJ}h_j \left\Vert  \max_{i\in\mI_j} \hat{Q}_{ij}^r -  \max_{i\in\mI_j} \hat{q}_{ij}^{*,r}\right\Vert_{T}>\ep\right)\nonumber\\
&\hspace{3cm}\leq \pr\left( \sum_{j\in\mJ,i\in\mI_j}h_j \left\Vert   \hat{Q}_{ij}^r - \hat{q}_{ij}^{*,r}\right\Vert_{T}>\ep\right)\nonumber\\
&\hspace{3cm}\leq \pr\left( \sum_{j\in\mJ,i\in\mI_j}\left\Vert   \hat{Q}_{ij}^r - \hat{q}_{ij}^{*,r}\right\Vert_{T}>\frac{\ep}{\max_{j\in\mJ}h_j}\right)\nonumber\\
&\hspace{3cm}= \pr\left( \sum_{i\in\mI}\sum_{j\in\mJ_i}\left\Vert   \hat{Q}_{ij}^r - \hat{q}_{ij}^{*,r}\right\Vert_{T}>\frac{\ep}{\max_{j\in\mJ}h_j}\right)\label{eq_p_main_1}\\
&\hspace{3cm}\leq \sum_{i\in\mI} \pr\left( \sum_{j\in\mJ_i}\left\Vert   \hat{Q}_{ij}^r - \hat{q}_{ij}^{*,r}\right\Vert_{T}>\ep_1\right),\label{eq_p_main_2}
\end{align}
where \eqref{eq_p_main_1} is by the fact that $ \bigcup_{j\in\mJ}\mI_j=\bigcup_{i\in\mI}\mJ_i $ and \eqref{eq_p_main_2} is because $\ep_1:=\ep/(I\max_{j\in\mJ}h_j)$. Therefore, it is enough to prove that \eqref{eq_p_main_2} converges to 0, which implies that the proposed policy should keep the number of jobs in front of the shared server $i$ close to the optimal solution process at all times for all $i\in\mI$. 

For notational convenience, let us define
\begin{equation*}
\bar{\lambda} := \max_{j\in\mJ} \lambda_j,\quad\bar{\mu} := \max_{j\in\mJ,i\in\mI_j} \mu_{ij},\quad \underline{\lambda} := \min_{j\in\mJ} \lambda_j,\quad \underline{\mu} := \min_{j\in\mJ,i\in\mI_j} \mu_{ij}.
\end{equation*}
Let $\tau_{in}^r:\Omega\rightarrow\R_+\cup\{\infty\}$ denote the start time of the $n$th review period (Step 2 or 3) in the shared server $i$ under the proposed policy for all $i\in\mI$ and $n,r\in\N_+$. For completeness, if $\tau_{in}^r(\omega)=\infty$ for some $i\in\mI$, $n,r\in\N_+$, and $\omega\in\Omega$, then $\tau_{im}^r(\omega):=\infty$ for all $m>n$. Then, $\tau_{i1}^r(\omega)=0$ and $\tau_{i(n+1)}^r(\omega) \geq\tau_{in}^r(\omega)$ for all $i\in\mI$, $n,r\in\N_+$, and $\omega\in\Omega$. Let $M^r:=1+\ru{(1+\bar{\mu})Jr^2T}$ for all $r\in\N_+$. Because $\tau_{in}^r$ is a service completion epoch in the shared server $i$ for all $n\geq 2$ and $i\in\mI$, we have
\begin{align}
\sum_{i\in\mI}\pr\left( \tau_{iM^r}^r \leq r^2T\right) &\leq \sum_{i\in\mI}\pr\left( \sum_{j\in\mJ_i} S_{ij}(r^2T) \geq M^r -1\right)\nonumber \\
& \leq \sum_{i\in\mI}\sum_{j\in\mJ_i}\pr\left( S_{ij}(r^2T) \geq (1+\bar{\mu})r^2T\right)\rightarrow 0,\label{eq_p_main_3}
\end{align}
where \eqref{eq_p_main_3} is by functional strong law of large numbers (FSLLN) for renewal processes (see Theorem 5.10 of \citet{che01}). The convergence result in \eqref{eq_p_main_3} implies that there are at most $O(r^2)$ review periods in the interval $[0,r^2T]$ in each shared server with a high probability when $r$ is sufficiently large, where $O(\cdot)$ denotes the big-$O$ notation.

With the convention that $\infty - \infty := \infty$, let us define the following sets for all $\ep_2>0$, $i\in\mI$, and $n,r\in\N_+$:
\begin{subequations}\label{eq_good_set}
\begin{align}
&\mA_{in}^{(1),r} (\ep_2):= \left\{ \tau_{in}^r > r^2T \right\}\cup \left\{ \tau_{i(n+1)}^r - \tau_{in}^r \leq \ep_2 r \right\},\label{eq_good_set_1}\\
&\mA_{in}^{(2),r} (\ep_2):= \left\{ \tau_{in}^r > r^2T \right\}\cup\bigcap_{j\in\mJ_i} \bigg\{ \sup_{\tau_{in}^r\leq t\leq \tau_{i(n+1)}^r}\left|Q_{ij}^r(t) -Q_{ij}^r(\tau_{in}^r)  \right| \leq C_5\ep_2 r \bigg\},\label{eq_good_set_2}\\
&\mA_{in}^{(3),r} (\ep_2):=  \left\{ \tau_{in}^r > r^2T \right\}\cup\bigcap_{j\in\mJ_i} \bigg\{ \sup_{\tau_{in}^r\leq t\leq \tau_{i(n+1)}^r}\left|q_{ij}^{*,r}(t) -q_{ij}^{*,r}(\tau_{in}^r)  \right| \leq C_6\ep_2 r \bigg\},\label{eq_good_set_3}\\
&\mA_{in}^{(4),r} (\ep_2):= \left\{ \tau_{in}^r > r^2T \right\}\cup \bigcap_{j\in\mJ_i} \bigg\{\left|Q_{ij}^r(\tau_{i(n+1)}^r) - q_{ij}^{*,r}(\tau_{i(n+1)}^r)  \right| \leq C_7\ep_2 r \bigg\},\label{eq_good_set_4}\\
&\mA_{in}^r (\ep_2):=\bigcap_{l=1}^4 \mA_{in}^{(l),r}(\ep_2),\label{eq_good_set_5}
\end{align}
\end{subequations}
where $C_5$, $C_6$, and $C_7$ are arbitrary strictly positive constants such that
\begin{equation}\label{eq_constants}
C_5 > \max\left\{4\bar{\lambda}, 2\bar{\mu}\right\},\qquad \left(3+ 2(J-1)\frac{\bar{\mu}}{\underline{\mu}}\right)C_6 < C_7<0.5\underline{\lambda}.
\end{equation}
We let $\mA_{i0}^r(\ep_2):=\Omega$ for all $i\in\mI$, $\ep_2>0$, and $r\in\N_+$ for completeness. 

The event in \eqref{eq_good_set_1} implies that the length of a review period is short in the shared server $i$ in $[0,r^2T]$. The event in \eqref{eq_good_set_2} implies that the queue length processes associated with the buffers in front of the shared server $i$ do not change a lot during a review period in $[0,r^2T]$. The event in \eqref{eq_good_set_3} implies that the optimal LP \eqref{eq_lp} solution does not change a lot during a review period in $[0,r^2T]$. The event in \eqref{eq_good_set_4} implies that the queue lengths in the buffers in front of the shared server $i$ do not deviate a lot from the optimal LP \eqref{eq_lp} solution at the end of a review period in $[0,r^2T]$. Finally, the following result states that the aforementioned events are realized jointly in the review periods $\{1,2,\ldots,M^r\}$ with high probability when $r$ is large.

\begin{lemma}\label{l_good_set}
For all $\ep_2>0$ and $i\in\mI$, we have
\begin{equation*}
\pr\left(\bigcap_{n=1}^{M^r} \mA_{in}^r(\ep_2)\right) \rightarrow 1.
\end{equation*}

\end{lemma}
The proof of Lemma \ref{l_good_set} is presented in Appendix \ref{l_good_set_proof}.

Let $\ep_2>0$ be such that $J(C_5+C_6+C_7)\ep_2\leq \ep_1$. Then, the probability in \eqref{eq_p_main_2} is less than or equal to
\begin{align}
&\sum_{i\in\mI}\pr\left( \sum_{j\in\mJ_i}\left\Vert   \hat{Q}_{ij}^r - \hat{q}_{ij}^{*,r}\right\Vert_{T}>\ep_1,\; \tau_{iM^r}^r > r^2T,\;  \bigcap_{n=1}^{M^r} \mA_{in}^r(\ep_2) \right)\label{eq_p_main_4}\\
&\hspace{6cm} +  \sum_{i\in\mI}\pr\left( \tau_{iM^r}^r\leq r^2T\right) +  \sum_{i\in\mI}\pr\left( \left(\bigcap_{n=1}^{M^r} \mA_{in}^r(\ep_2)\right)^c\; \right),\label{eq_p_main_5}
\end{align}
where the superscript $c$ denotes complement of the associated set. The sums in \eqref{eq_p_main_5} converge to 0 by \eqref{eq_p_main_3} and Lemma \ref{l_good_set}, respectively. Hence, it is enough to prove that the probability in \eqref{eq_p_main_4} converges to 0. The probability in \eqref{eq_p_main_4} is equal to
\begin{align}
& \sum_{i\in\mI}\pr\Bigg(\sum_{j\in\mJ_i} \sup_{0\leq t\leq r^2T}\left| Q_{ij}^r(t) - q_{ij}^{*,r}(t)\right|>\ep_1r,\; \tau_{iM^r}^r > r^2T,\;  \bigcap_{n=1}^{M^r} \mA_{in}^r(\ep_2)\Bigg)\nonumber\\
&\hspace{1cm}\leq \sum_{i\in\mI}\pr\Bigg( \bigcup_{n=1}^{M^r} \Bigg\{\sum_{j\in\mJ_i} \sup_{\tau_{in}^r\leq t \leq \tau_{i(n+1)}^r}\left| Q_{ij}^r(t) - q_{ij}^{*,r}(t)\right|>\ep_1r,\;\tau_{in}^r\leq r^2T \Bigg\},\; \bigcap_{n=1}^{M^r} \mA_{in}^r(\ep_2) \Bigg)\nonumber\\
&\hspace{1cm}\leq \sum_{i\in\mI}\pr\Bigg( \bigcup_{n=1}^{M^r} \Bigg\{\sum_{j\in\mJ_i} \bigg( \sup_{\tau_{in}^r\leq t \leq \tau_{i(n+1)}^r}\left| Q_{ij}^r(t) - Q_{ij}^r(\tau_{in}^r) \right| + \left| Q_{ij}^r(\tau_{in}^r) - q_{ij}^{*,r}(\tau_{in}^r) \right|\nonumber\\ 
&\hspace{3cm} + \sup_{\tau_{in}^r\leq t \leq \tau_{i(n+1)}^r}\left| q_{ij}^{*,r}(t) - q_{ij}^{*,r}(\tau_{in}^r) \right|\bigg)>\ep_1r,\;\tau_{in}^r\leq r^2T \Bigg\},\; \bigcap_{n=1}^{M^r} \mA_{in}^r(\ep_2) \Bigg).\label{eq_p_main_6}
\end{align}
In the set $\{\tau_{in}^r\leq r^2T\}\cap \mA_{i(n-1)}^r(\ep_2) \cap \mA_{in}^r(\ep_2)$, 
\begin{align}
&\sum_{j\in\mJ_i} \bigg( \sup_{\tau_{in}^r\leq t \leq \tau_{i(n+1)}^r}\left| Q_{ij}^r(t) - Q_{ij}^r(\tau_{in}^r) \right| + \left| Q_{ij}^r(\tau_{in}^r) - q_{ij}^{*,r}(\tau_{in}^r) \right| + \sup_{\tau_{in}^r\leq t \leq \tau_{i(n+1)}^r}\left| q_{ij}^{*,r}(t) - q_{ij}^{*,r}(\tau_{in}^r) \right| \bigg)\nonumber\\
&\hspace{10cm} \leq J(C_5+C_6+C_7)\ep_2 r \leq \ep_1 r\label{eq_p_main_7}
\end{align}
for all $n\in\{1,2,\ldots,M^r\}$, $i\in\mI$, and $r\in\N_+$ by \eqref{eq_good_set}. Hence, the event inside the probability in \eqref{eq_p_main_6} is equal to $\emptyset$ by \eqref{eq_p_main_7} for all $i\in\mI$. Therefore, the sum in \eqref{eq_p_main_6} is equal to 0.

\section{Proof of Lemma \ref{l_good_set}}\label{l_good_set_proof}

Let us define
\begin{equation}\label{eq_gsp_0}
B^r:=\left\{ \max_{j\in\mJ} A_j^r(r^2T) \vee \max_{i\in\mI,j\in\mJ_i} S_{ij}(r^2T) \vee \max_{j\in\mJ,k\in\mK_j} S_k(r^2T) \leq \rd{C_8r^2T} - 1 \right\},
\end{equation}
where 
\begin{equation*}
C_8 := 1+ \bar{\lambda} \vee \bar{\mu} \vee \max_{j\in\mJ,k\in\mK_j} \mu_k.
\end{equation*}
By FSLLN, we have $\pr\left(B^r\right)\rightarrow 1$.

We present the following lemmas which will be useful later. The first one provides an exponential tail bound for renewal processes.

\begin{lemma}\label{l_tail_1}
Let us fix arbitrary $a>0$ and $b>0$. There exists an $r_1\in\N_+$ such that if $r\geq r_1$, then for all $j\in\mJ$, $i\in\mI_j$, $k\in\mK_j$, $l\in\mI$, and $n\in\{1,2,\ldots,M^r\}$, we have
\begin{subequations}\label{eq_t1_0}
\begin{align}
&\pr\left( \sup_{0\leq t\leq ar} \left| A_{j}^{r}(\tau_{ln}^r + t) -  A_{j}^{r}(\tau_{ln}^r) -\lambda_j^r t \right| >br,\;\tau_{ln}^r \leq r^2T,\;B^r\right) \leq C_9r^{2}\e^{-C_{10}r},\label{eq_t1_1}\\
&\pr\left( \sup_{0\leq t\leq ar} \left| S_{ij}(T_{ij}^r(\tau_{ln}^r) + t) -  S_{ij}(T_{ij}^r(\tau_{ln}^r)) -\mu_{ij}t \right| >br,\;\tau_{ln}^r \leq r^2T,\;B^r\right) \leq C_9r^{2}\e^{-C_{10}r},\label{eq_t1_2}\\
&\pr\left( \sup_{0\leq t\leq ar} \left| S_k(T_k^r(\tau_{ln}^r) + t) -  S_k(T_k^r(\tau_{ln}^r)) -\mu_k t \right| >br,\;\tau_{ln}^r \leq r^2T,\;B^r\right) \leq C_9r^{2}\e^{-C_{10}r},\label{eq_t1_3}
\end{align}
\end{subequations}
where $C_9$ and $C_{10}$ are strictly positive constants independent of $i$, $j$, $k$, $l$, $n$, and $r$.
\end{lemma}
The proof of Lemma \ref{l_tail_1} is presented in E-companion \ref{l_tail_1_proof}. 

The second lemma states that the length of Step 2 in Definition \ref{d_policy} is short with high probability when $r$ is large. By Assumption \ref{a_ht} Part 1, there exists an $r_2\in\N_+$ such that if $r\geq r_2$, then $0.5\lambda_j <\lambda_j^r < 2\lambda_j$ for all $j\in\mJ$.

\begin{lemma}\label{l_tail_2}
For all $i\in\mI$, $n\in\{1,2,\ldots,M^r\}$, $b>0$, and $r\geq r_2$, if the $n$th review period in the shared server $i$ is Step 2 in Definition \ref{d_policy}, then
\begin{equation*}
\pr\left(\tau_{i(n+1)}^r - \tau_{in}^r>b r,\; \tau_{in}^r\leq r^2T,\;B^r\right)\leq C_{11}r^2\e^{-C_{12}r},
\end{equation*}
where $C_{11}$ and $C_{12}$ are strictly positive constants independent of $i$, $n$, and $r$.
\end{lemma}
The proof of Lemma \ref{l_tail_2} is presented in E-companion \ref{l_tail_2_proof}. 

The third lemma states that the length of Step 3 in Definition \ref{d_policy} is short and the buffer content of the job type that is not processed in Step 3 does not grow a lot in the review period with high probability when $r$ is large. 

\begin{lemma}\label{l_tail_3}
Fix arbitrary $i\in\mI$, $n\in\{1,2,\ldots,M^r\}$, $a>0$, and $b>0$. Suppose that the $n$th review period in the shared server $i$ is Step 3 in Definition \ref{d_policy}. Without loss of generality, let $m\in\mJ_i^{<,r}(\tau_{in}^r)$ denote the job type that the shared server $i$ does not process in the $n$th review period. Then, there exists an $r_3\in\N_+$ such that $r_3$ is independent of $n$ and if $r\geq r_3$, we have
\begin{align}
&\pr\Bigg(\left\{\tau_{i(n+1)}^r - \tau_{in}^r>\frac{\ru{q_{im}^{*,r}(\tau_{in}^r)}-Q_{im}^r(\tau_{in}^r)}{\lambda_m}+ br \right\}\cup \left\{Q_{im}^r(\tau_{i(n+1)}^r) -\ru{q_{im}^{*,r}(\tau_{in}^r)} > 2b\lambda_m r\right\},\nonumber\\
&\hspace{4.5cm} \max_{j\in\mJ_i}\left|Q_{ij}^r(\tau_{in}^r) -\ru{q_{ij}^{*,r}(\tau_{in}^r)}\right| \leq ar,\;\tau_{in}^r\leq r^2T,\;B^r\Bigg) \leq C_{13}r^{4}\e^{-C_{14}r},\label{eq_l_tail_3_0}
\end{align}
where $C_{13}$ and $C_{14}$ are strictly positive constants independent of $n$ and $r$.
\end{lemma}
The proof of Lemma \ref{l_tail_3} is presented in E-companion \ref{l_tail_3_proof}. 

The fourth lemma states that the workload amounts in the shared servers and the number of jobs waiting in front of the dedicated servers do not fluctuate a lot within a time interval with length $O(r)$ with high probability when $r$ is large.

\begin{lemma}\label{l_tail_4}
Fix arbitrary $a>0$ and $b>0$. There exists an $r_4\in\N_+$ such that if $r\geq r_4$, then for all $l\in\mI$ and $n\in\{1,2,\ldots,M^r\}$, we have
\begin{align}
&\pr\Bigg(\sup_{0\leq t\leq ar } \left(\max_{i\in\mI} \left|W_i^{r}(\tau_{ln}^r+t)- W_i^{r}(\tau_{ln}^r) \right| \vee \max_{j\in\mJ,k\in\mK_j} \left|Q_k^{r}(\tau_{ln}^r+ t)- Q_k^{r}(\tau_{ln}^r) \right|\right)>b r, \nonumber\\
&\hspace{10cm} \tau_{ln}^r \leq r^2T,\;B^r\Bigg)\leq C_{15}r^5\e^{-C_{16}r},\label{eq_t4_0}
\end{align}
where $C_{15}$ and $C_{16}$ are strictly positive constants independent of $l$, $n$, and $r$.
\end{lemma}
The proof of Lemma \ref{l_tail_4} is presented in E-companion \ref{l_tail_4_proof}.

Because $\pr\left(B^r\right)\rightarrow 1$, proving Lemma \ref{l_good_set} is equivalent to proving 
\begin{equation*}
\pr\left( \bigcup_{n=1}^{M^r}\left(\mA_{in}^r(\ep_2)\right)^c,\;B^r\right)\rightarrow 0,\qquad\forall \ep_2>0,\; i\in\mI.
\end{equation*}

Let us fix arbitrary $\ep_2>0$ and $i\in\mI$. Let $\{A_n,n\in\N_+\}$ be an arbitrary sequence of sets. One can see that
\begin{equation}\label{eq_gsp_1}
\bigcup_{n=1}^M A_n=A_1 \cup \bigcup_{n=2}^M \left(A_n\cap A_{n-1}^c \cap A_{n-2}^c \cap \ldots \cap A_1^c\right),\qquad\forall M\in\N_+.
\end{equation}
Therefore, we have
\begin{equation*}
\pr\left( \bigcup_{n=1}^{M^r}\left(\mA_{in}^r\right)^c,B^r\right)\leq \sum_{n=1}^{M^r} \pr\left(\left(\mA_{in}^r\right)^c \cap \mA_{i(n-1)}^r \cap B^r\right).
\end{equation*}

Let us fix an arbitrary $n\in\{1,2,\ldots,M^r\}$. By \eqref{eq_good_set_5} and \eqref{eq_gsp_1}, we have
\begin{subequations}\label{eq_gsp_2}
\begin{align}
\pr\left(\left(\mA_{in}^r\right)^c \cap \mA_{i(n-1)}^r \cap B^r\right) &\leq \pr\left(\left(\mA_{in}^{(1),r}\right)^c \cap \mA_{i(n-1)}^r \cap B^r\right)\label{eq_gsp_21}\\
&\hspace{1cm} +\pr\left(\left(\mA_{in}^{(2),r}\right)^c \cap \mA_{in}^{(1),r} \cap B^r\right)\label{eq_gsp_22}\\
&\hspace{1cm} +\pr\left(\left(\mA_{in}^{(3),r}\right)^c \cap \mA_{in}^{(1),r} \cap B^r\right)\label{eq_gsp_23}\\
&\hspace{1cm} +\pr\left(\left(\mA_{in}^{(4),r}\right)^c \cap \mA_{in}^{(1),r}  \cap \mA_{in}^{(3),r} \cap \mA_{i(n-1)}^r \cap B^r\right).\label{eq_gsp_24}
\end{align}
\end{subequations}
We will consider the probabilities in the RHS of \eqref{eq_gsp_2} one by one.

\textbf{\textit{The probability in the RHS of \eqref{eq_gsp_21}:}} By \eqref{eq_good_set_1}, it is equal to
\begin{equation}\label{eq_a_1}
\pr\left( \tau_{i(n+1)}^r - \tau_{in}^r > \ep_2 r,\; \tau_{in}^r\leq r^2T,\;\mA_{i(n-1)}^r,\; B^r \right).
\end{equation}

Suppose that the $n$th review period in the shared server $i$ is Step 2 in Definition \ref{d_policy}. By Lemma \ref{l_tail_2}, if $r\geq r_2$, \eqref{eq_a_1} is less than or equal to 
\begin{equation}\label{eq_a_2}
C_{17}r^2\e^{-C_{18}r},
\end{equation}
where $C_{17}$ and $C_{18}$ are strictly positive constants independent of $i$, $n$, and $r$.

Suppose that the $n$th review period in the shared server $i$ is Step 3 in Definition \ref{d_policy}. Without loss of generality, let $m\in\mJ_i^{<,r}(\tau_{in}^r)$ denote the job type that the shared server $i$ does not process in the $n$th review period. There exists an $r_5\in\N_+$ such that if $r\geq r_5$,
\begin{equation*}
1+C_7\ep_2 r \leq 2C_7\ep_2 r.
\end{equation*}
Hence, by \eqref{eq_good_set_4}, in the set $\mA_{i(n-1)}^r \cap\{ \tau_{in}^r\leq r^2T\}$, if $r\geq r_5$,
\begin{equation*}
\max_{j\in\mJ_i}\left|Q_{ij}^r(\tau_{in}^r) -\ru{q_{ij}^{*,r}(\tau_{in}^r)}\right| \leq 1+C_7\ep_2 r \leq2C_7\ep_2 r.
\end{equation*}
Let $b:=(1-2C_7/ \underline{\lambda})\ep_2$. Then, $b>0$ by \eqref{eq_constants}. The probability in \eqref{eq_a_1} is less than or equal to
\begin{align}
& \pr\left( \tau_{i(n+1)}^r - \tau_{in}^r > \left(\frac{2C_7\ep_2}{\lambda_m}+ b \right) r,\; \tau_{in}^r\leq r^2T,\;\mA_{i(n-1)}^r,\; B^r \right)\nonumber\\
&\hspace{1cm}\leq \pr\Bigg( \tau_{i(n+1)}^r - \tau_{in}^r > \left(\frac{2C_7\ep_2}{\lambda_m}+b \right) r,\nonumber\\
&\hspace{5cm} \max_{j\in\mJ_i}\left|Q_{ij}^r(\tau_{in}^r) -\ru{q_{ij}^{*,r}(\tau_{in}^r)}\right| \leq 2C_7\ep_2r,\;\tau_{in}^r\leq r^2T,\;B^r \Bigg),\label{eq_a_2_1}
\end{align}
where \eqref{eq_a_2_1} holds if $r\geq r_5$. Let us invoke Lemma \ref{l_tail_3} by letting $a:=2C_7\ep_2$ so that we can derive that there exists an $r_6\geq r_5$ such that $r_6$ is independent of $n$ and if $r\geq r_6$, \eqref{eq_a_2_1} is less than or equal to 
\begin{equation}\label{eq_a_3}
C_{19}r^{4}\e^{-C_{20}r},
\end{equation}
where $C_{19}$ and $C_{20}$ are strictly positive constants independent of $n$ and $r$.

Therefore, by \eqref{eq_a_2} and \eqref{eq_a_3}, if $r\geq r_2 \vee r_6$, \eqref{eq_a_1} is less than or equal to
\begin{equation}\label{eq_a_4}
\left(C_{17}\vee C_{19}\right) r^{4}\e^{-\left(C_{18}\wedge C_{20}\right)r}.
\end{equation}

\textbf{\textit{The probability in \eqref{eq_gsp_22}:}} By \eqref{eq_good_set_2}, it is less than or equal to
\begin{align}
&\sum_{j\in\mJ_i}\pr\left( \sup_{\tau_{in}^r\leq t\leq \tau_{i(n+1)}^r}\left|Q_{ij}^r(t) -Q_{ij}^r(\tau_{in}^r)  \right| > C_5\ep_2 r,\;\tau_{in}^r \leq r^2T,\; \mA_{in}^{(1),r},\;B^r\right)\nonumber\\
&\hspace{0cm}\leq \sum_{j\in\mJ_i}\pr\left( \sup_{0\leq t\leq \ep_2 r}\left|Q_{ij}^r(\tau_{in}^r +t) -Q_{ij}^r(\tau_{in}^r)  \right| > C_5\ep_2 r,\;\tau_{in}^r \leq r^2T,\;B^r \right)\nonumber\\
&\hspace{0cm}= \sum_{j\in\mJ_i}\pr\bigg( \sup_{0\leq t\leq \ep_2 r}\left|A_{j}^r(\tau_{in}^r+t)-S_{ij}(T_{ij}^r(\tau_{in}^r+t)) - A_{j}^r(\tau_{in}^r)+S_{ij}(T_{ij}^r(\tau_{in}^r))\right| > C_5\ep_2 r, \nonumber\\
&\hspace{13cm} \tau_{in}^r \leq r^2T,\;B^r\bigg)\nonumber\\
&\hspace{0cm}\leq \sum_{j\in\mJ_i}\pr\left( S_{ij}(T_{ij}^r(\tau_{in}^r)+\ep_2 r)-S_{ij}(T_{ij}^r(\tau_{in}^r)) > 0.5C_5\ep_2 r,\;\tau_{in}^r \leq r^2T,\;B^r \right)\nonumber\\
&\hspace{2cm}+\sum_{j\in\mJ_i}\pr\left( A_{j}^r(\tau_{in}^r+\ep_2 r) - A_{j}^r(\tau_{in}^r) > 0.5C_5\ep_2 r,\;\tau_{in}^r \leq r^2T,\;B^r \right),\label{eq_b_0}\\
&\hspace{0cm}\leq\sum_{j\in\mJ_i}\pr\left( S_{ij}(T_{ij}^r(\tau_{in}^r)+\ep_2 r)-S_{ij}(T_{ij}^r(\tau_{in}^r)) -\mu_{ij}\ep_2 r > \left(0.5C_5-\bar{\mu}\right)\ep_2 r,\;\tau_{in}^r \leq r^2T,\;B^r \right)\nonumber\\
&\hspace{2cm}+\sum_{j\in\mJ_i}\pr\left( A_{j}^r(\tau_{in}^r+\ep_2 r) - A_{j}^r(\tau_{in}^r) -\lambda_j^r \ep_2 r > \left(0.5C_5-2\bar{\lambda}\right)\ep_2 r,\;\tau_{in}^r \leq r^2T,\;B^r \right),\label{eq_b_1}
\end{align}
where \eqref{eq_b_0} is by triangular inequality and the fact that $T_{ij}^r(\tau_{in}^r+\ep_2 r) \leq T_{ij}^r(\tau_{in}^r)+\ep_2 r$ and \eqref{eq_b_1} holds if $r\geq r_2$. By \eqref{eq_constants} and Lemma \ref{l_tail_1}, there exists an $r_7\in\N_+$ such that $r_7$ is independent of $i$ and $n$ and if $r\geq r_7$, the sum in \eqref{eq_b_1} is less than or equal to 
\begin{equation}\label{eq_b_5}
C_{21}r^2\e^{-C_{22}r}
\end{equation}
where $C_{21}$ and $C_{22}$ are strictly positive constants independent of $i$, $n$, and $r$. Finally, if $r\geq r_2\vee r_7$, the probability in \eqref{eq_gsp_22} is less than or equal to \eqref{eq_b_5}.

\textbf{\textit{The probability in \eqref{eq_gsp_23}:}} By \eqref{eq_good_set_3}, it is less than or equal to
\begin{align}
&\sum_{j\in\mJ_i}\pr\left( \sup_{\tau_{in}^r\leq t\leq \tau_{i(n+1)}^r}\left|q_{ij}^{*,r}(t)- q_{ij}^{*,r}(\tau_{in}^r) \right| > C_6\ep_2 r,\;\tau_{in}^r \leq r^2T,\; \mA_{in}^{(1),r},\;B^r\right)\label{eq_c_0}\\
&\hspace{0.5cm}\leq \sum_{j\in\mJ_i}\pr\Bigg(C_2\sup_{\tau_{in}^r\leq t\leq \tau_{i(n+1)}^r} \left(\max_{l\in\mI} \left|W_l^{r}(t)- W_l^{r}(\tau_{in}^r) \right| \vee \max_{x\in\mJ,k\in\mK_x} \left|Q_k^{r}(t)- Q_k^{r}(\tau_{in}^r) \right|\right)>C_6\ep_2 r,\nonumber\\
&\hspace{11cm} \tau_{in}^r \leq r^2T,\; \mA_{in}^{(1),r},\;B^r\Bigg),\label{eq_c_1}\\
&\hspace{0.5cm}\leq \sum_{j\in\mJ_i}\pr\Bigg(\sup_{0\leq t\leq \ep_2 r} \left(\max_{l\in\mI} \left|W_l^{r}(\tau_{in}^r+t)- W_l^{r}(\tau_{in}^r) \right| \vee \max_{x\in\mJ,k\in\mK_x} \left|Q_k^{r}(\tau_{in}^r+ t)- Q_k^{r}(\tau_{in}^r) \right|\right)\nonumber\\
&\hspace{10cm}>\frac{C_6}{C_2}\ep_2 r,\; \tau_{in}^r \leq r^2T,\;B^r\Bigg),\label{eq_c_2}
\end{align}
where \eqref{eq_c_1} is by Lemma \ref{l_reg}. By Lemma \ref{l_tail_4}, there exists an $r_8\in\N_+$ such that $r_8$ is independent of $i$ and $n$ and if $r\geq r_8$, the sum in \eqref{eq_c_2} is less than or equal to
\begin{equation}\label{eq_c_3}
C_{23}r^5 \e^{-C_{24}r},
\end{equation}
where $C_{23}$ and $C_{24}$ are strictly positive constants independent of $i$, $n$, and $r$.

\textbf{\textit{The probability in \eqref{eq_gsp_24}:}} By \eqref{eq_good_set_4}, it is less than or equal to
\begin{align}
&\sum_{j\in\mJ_i}\pr\left(\left|Q_{ij}^r(\tau_{i(n+1)}^r) - q_{ij}^{*,r}(\tau_{i(n+1)}^r)  \right| > C_7\ep_2 r,\;\tau_{in}^r \leq r^2T,\; \mA_{in}^{(1),r},\; \mA_{in}^{(3),r},\; \mA_{i(n-1)}^r,\;B^r\right)\nonumber\\
&\hspace{0cm} \leq\sum_{j\in\mJ_i}\pr\left(\left|Q_{ij}^r(\tau_{i(n+1)}^r) - q_{ij}^{*,r}(\tau_{in}^r)  \right| > (C_7-C_6)\ep_2 r,\;\tau_{in}^r \leq r^2T,\;\mA_{in}^{(3),r}, \mA_{i(n-1)}^r,B^r\right)\label{eq_d_1}\\
&\hspace{3cm} + \sum_{j\in\mJ_i}\pr\left(\left|q_{ij}^{*,r}(\tau_{i(n+1)}^r)- q_{ij}^{*,r}(\tau_{in}^r)  \right| > C_6\ep_2 r,\;\tau_{in}^r \leq r^2T,\; \mA_{in}^{(1),r},\;B^r\right),\label{eq_d_2}
\end{align}
where we use the fact that $C_7>C_6$ (see \eqref{eq_constants}). The sum in \eqref{eq_d_2} is less than or equal to the sum in \eqref{eq_c_0}. Therefore, if $r\geq r_8$, the sum in \eqref{eq_d_2} is less than or equal to the term in \eqref{eq_c_3}.

Next, let us consider the sum in \eqref{eq_d_1}. First, suppose that the $n$th review period is Step 2 of Definition \ref{d_policy}, which implies that $Q_{ij}^r(\tau_{in}^r) \leq\ru{ q_{ij}^{*,r}(\tau_{in}^r)}$ for all $j\in\mJ_i$. By \eqref{eq_we}, we can derive that 
\begin{equation}\label{eq_d_3}
\max_{j\in\mJ_i}\left|Q_{ij}^r(\tau_{in}^r) - q_{ij}^{*,r}(\tau_{in}^r) \right| \leq 1 \vee\bigg( \bar{\mu}\sum_{i\in\mI, j\in\mJ_i} \frac{1}{\mu_{ij}}\bigg)=:C_{25}.
\end{equation}
Recall that Step 2 in shared server $i$ ends with the first service completion in that server. Therefore, 
\begin{align}
&Q_{ij}^r(\tau_{in}^r) -1 \leq Q_{ij}^r(\tau_{i(n+1)}^r) \leq Q_{ij}^r(\tau_{in}^r)  + A_j^r(\tau_{i(n+1)}^r) -   A_j^r(\tau_{in}^r),\quad\forall j\in\mJ_i,\nonumber\\
\implies& \left|Q_{ij}^r(\tau_{i(n+1)}^r) - Q_{ij}^r(\tau_{in}^r)\right| \leq 1+ A_j^r(\tau_{i(n+1)}^r) -   A_j^r(\tau_{in}^r),\quad\forall j\in\mJ_i.\label{eq_d_4}
\end{align}
By \eqref{eq_d_3} and \eqref{eq_d_4},
\begin{equation*}
\left|Q_{ij}^r(\tau_{i(n+1)}^r) - q_{ij}^{*,r}(\tau_{in}^r) \right| \leq 2C_{25}+ A_j^r(\tau_{i(n+1)}^r) -   A_j^r(\tau_{in}^r),\quad\forall j\in\mJ_i.
\end{equation*}
Therefore, the sum in \eqref{eq_d_1} is less than or equal to
\begin{equation}\label{eq_d_5}
\sum_{j\in\mJ_i}\pr\left( A_j^r(\tau_{i(n+1)}^r) -   A_j^r(\tau_{in}^r)  > (C_7-C_6)\ep_2 r -2C_{25},\;\tau_{in}^r \leq r^2T,\;B^r\right).
\end{equation}
Let $b_1:= (C_7-C_6)\ep_2/(4\bar{\lambda})$. Then, the sum in \eqref{eq_d_5} is less than or equal to
\begin{align}
&\sum_{j\in\mJ_i}\pr\left( A_j^r(\tau_{in}^r + b_1r) -  A_j^r(\tau_{in}^r)  > (C_7-C_6)\ep_2 r -2C_{25},\;\tau_{in}^r \leq r^2T,\;B^r\right)\label{eq_d_6}\\
&\hspace{7cm}+ \sum_{j\in\mJ_i}\pr\left( \tau_{i(n+1)}^r -\tau_{in}^r  > b_1r,\;\tau_{in}^r \leq r^2T,\;B^r\right).\label{eq_d_7}
\end{align}
By Lemma \ref{l_tail_2}, if $r\geq r_2$, then the sum in \eqref{eq_d_7} is less than or equal to
\begin{equation}\label{eq_d_8}
C_{26} r^2\e^{-C_{27}r},
\end{equation}
where $C_{26}$ and $C_{27}$ are strictly positive constants independent of $i$, $n$, and $r$. 

There exists an $r_9\in\N_+$ such that if $r\geq r_9$, then 
\begin{equation*}
0.5(C_7-C_6)\ep_2 r-2C_{25} \geq 0.25(C_7-C_6)\ep_2 r.
\end{equation*}
Hence, if $r\geq r_2\vee r_9$, the sum in \eqref{eq_d_6} is less than or equal to
\begin{align}
&\sum_{j\in\mJ_i}\pr\left( A_j^r(\tau_{in}^r + b_1r) -  A_j^r(\tau_{in}^r) -\lambda_j^r b_1r > (C_7-C_6)\ep_2 r -2\bar{\lambda}b_1r -2C_{25},\;\tau_{in}^r \leq r^2T,\;B^r\right)\nonumber\\
&= \sum_{j\in\mJ_i}\pr\left( A_j^r(\tau_{in}^r + b_1r) -  A_j^r(\tau_{in}^r) -\lambda_j^r b_1r > 0.5(C_7-C_6)\ep_2 r-2C_{25},\;\tau_{in}^r \leq r^2T,\;B^r\right)\nonumber\\
&\leq \sum_{j\in\mJ_i}\pr\left( A_j^r(\tau_{in}^r + b_1r) -  A_j^r(\tau_{in}^r) -\lambda_j^r b_1r > 0.25(C_7-C_6)\ep_2 r,\;\tau_{in}^r \leq r^2T,\;B^r\right).\label{eq_d_9}
\end{align}
By Lemma \ref{l_tail_1}, there exists an $r_{10}\in\N_+$ such that $r_{10}\geq r_2\vee r_9$, $r_{10}$ is independent of $i$ and $n$, and if $r\geq r_{10}$, the sum in \eqref{eq_d_9} is less than or equal to 
\begin{equation}\label{eq_d_10}
C_{28}r^2\e^{-C_{29}r},
\end{equation}
where $C_{28}$ and $C_{29}$ are strictly positive constants independent of $i$, $n$, and $r$.

Second, suppose that the $n$th review period is Step 3 of Definition \ref{d_policy}. Without loss of generality, let $m\in\mJ_i^{<,r}(\tau_{in}^r)$ denote the job type that the shared server $i$ does not process in the $n$th review period. Then, the sum in \eqref{eq_d_1} is less than or equal to
\begin{align}
&\sum_{j\in\mJ_i}\pr\bigg(\left|Q_{ij}^r(\tau_{i(n+1)}^r) - q_{ij}^{*,r}(\tau_{in}^r)  \right| > (C_7-C_6)\ep_2 r,\;\tau_{in}^r \leq r^2T,\; \mA_{in}^{(3),r},\;B^r,\nonumber\\
&\hspace{8cm} Q_{im}^r(\tau_{i(n+1)}^r) -\ru{q_{im}^{*,r}(\tau_{in}^r)} \leq C_6\ep_2r\bigg)\label{eq_d_11}\\
&\hspace{1cm}+ \sum_{j\in\mJ_i}\pr\left(Q_{im}^r(\tau_{i(n+1)}^r) -\ru{q_{im}^{*,r}(\tau_{in}^r)} > C_6\ep_2r,\;\tau_{in}^r\leq r^2T,\; \mA_{i(n-1)}^r,\;B^r\right).\label{eq_d_12}
\end{align}
By Lemma \ref{l_tail_3}, there exists an $r_{11}\in\N_+$ such that $r_{11}$ is independent of $n$ and if $r\geq r_{11}$, the sum in \eqref{eq_d_12} is less than or equal to
\begin{equation}\label{eq_d_13}
C_{30}r^{4}\e^{-C_{31}r},
\end{equation}
where $C_{30}$ and $C_{31}$ are strictly positive constants independent of $n$ and $r$. 

Next, let us consider the sum in \eqref{eq_d_11}. By definition of Step 3 (see Definition \ref{d_policy}), if $j\in\mJ_i^{>,r}(\tau_{in}^r)$, then $Q_{ij}^r(\tau_{i(n+1)}^r) =\ru{ q_{ij}^{*,r}(\tau_{in}^r)}$. If $j\in\mJ_i^{\leq,r}(\tau_{in}^r)\backslash\{m\}$, then $Q_{ij}^r(\tau_{i(n+1)}^r) \leq \ru{ q_{ij}^{*,r}(\tau_{in}^r)}$. Therefore, given that $Q_{im}^r(\tau_{i(n+1)}^r) -\ru{q_{im}^{*,r}(\tau_{in}^r)} \leq C_6\ep_2r$, we have
\begin{equation}\label{eq_d_14}
Q_{ij}^r(\tau_{i(n+1)}^r) - q_{ij}^{*,r}(\tau_{in}^r) \leq 1+ C_6\ep_2r,\qquad\forall j\in\mJ_i.
\end{equation}
By \eqref{eq_we}, we have
\begin{equation}\label{eq_d_15}
\sum_{j\in\mJ_i} \frac{Q_{ij}^r(\tau_{i(n+1)}^r)}{\mu_{ij}}= \sum_{j\in\mJ_i} \frac{q_{ij}^{*,r}(\tau_{i(n+1)}^r)}{\mu_{ij}}.
\end{equation}
Hence, in the set $\{\tau_{in}^r \leq r^2T\}\cap\mA_{in}^{(3),r}$, for all $j\in\mJ_i$, we have
\begin{align}
q_{ij}^{*,r}(\tau_{in}^r) -Q_{ij}^r(\tau_{i(n+1)}^r) &\leq q_{ij}^{*,r}(\tau_{i(n+1)}^r) -Q_{ij}^r(\tau_{i(n+1)}^r) + C_6\ep_2r \label{eq_d_16}\\
&=  C_6\ep_2r + \mu_{ij} \sum_{l\in\mJ_i\backslash\{j\}} \frac{Q_{il}^r(\tau_{i(n+1)}^r)-q_{il}^{*,r}(\tau_{i(n+1)}^r)}{\mu_{il}}\label{eq_d_17}\\
&\leq C_6\ep_2r + \mu_{ij} \sum_{l\in\mJ_i\backslash\{j\}} \frac{Q_{il}^r(\tau_{i(n+1)}^r)-q_{il}^{*,r}(\tau_{in}^r)+C_6\ep_2r}{\mu_{il}}\label{eq_d_18}\\
&\leq C_6\ep_2r + \mu_{ij} \sum_{l\in\mJ_i\backslash\{j\}} \frac{1+ 2C_6\ep_2r}{\mu_{il}}\label{eq_d_19}\\
&\leq C_6\ep_2r +\frac{\bar{\mu}}{\underline{\mu}}(J-1)(1+ 2C_6\ep_2r)\nonumber\\
& =\left(1+ 2(J-1)\frac{\bar{\mu}}{\underline{\mu}}\right) C_6\ep_2r + (J-1)\frac{\bar{\mu}}{\underline{\mu}},\label{eq_d_20}
\end{align}
where \eqref{eq_d_16} is by definition of the set $\{\tau_{in}^r \leq r^2T\}\cap\mA_{in}^{(3),r}$, \eqref{eq_d_17} is by \eqref{eq_d_15}, \eqref{eq_d_18} is by definition of the set $\{\tau_{in}^r \leq r^2T\}\cap\mA_{in}^{(3),r}$, and \eqref{eq_d_19} is by \eqref{eq_d_14}.

There exists an $r_{12}\in\N_+$ such that if $r\geq r_{12}$,
\begin{equation*}
C_6\ep_2r \geq  (J-1)\frac{\bar{\mu}}{\underline{\mu}}.
\end{equation*}
By \eqref{eq_d_14} and \eqref{eq_d_20}, if $r\geq r_{12}$,
\begin{equation}\label{eq_d_21}
\left| Q_{ij}^r(\tau_{i(n+1)}^r) - q_{ij}^{*,r}(\tau_{in}^r) \right| \leq \left(2+ 2(J-1)\frac{\bar{\mu}}{\underline{\mu}}\right) C_6\ep_2r,\qquad\forall j\in\mJ_i.
\end{equation}
By \eqref{eq_constants}, we have
\begin{equation}\label{eq_d_22}
C_7-C_6 >\left(2+ 2(J-1)\frac{\bar{\mu}}{\underline{\mu}}\right) C_6.
\end{equation}
Therefore, by \eqref{eq_d_21} and \eqref{eq_d_22}, if $r\geq r_{12}$, the sum in \eqref{eq_d_11} is equal to 0.

Let $r_{13}:= \max\{r_2,r_6,r_7,r_8,r_{10},r_{11},r_{12}\}$. Then $r_{13}$ is independent of $n$. Finally, by \eqref{eq_a_4}, \eqref{eq_b_5}, \eqref{eq_c_3}, \eqref{eq_d_8}, \eqref{eq_d_10}, and \eqref{eq_d_13}, if $r\geq r_{13}$, 
\begin{equation*}
\pr\left(\left(\mA_{in}^r\right)^c \cap \mA_{i(n-1)}^r \cap B^r\right) \leq C_{32} r^{5}\e^{-C_{33}r},
\end{equation*}
where $C_{32}$ and $C_{33}$ are strictly positive constants independent of $n$ and $r$. Therefore, if $r\geq r_{13}$, 
\begin{align*}
\pr\left( \bigcup_{n=1}^{M^r}\left(\mA_{in}^r\right)^c,B^r\right)&\leq \sum_{n=1}^{M^r} \pr\left(\left(\mA_{in}^r\right)^c \cap \mA_{i(n-1)}^r \cap B^r\right)\\
&\leq M^rC_{32} r^{5}\e^{-C_{33}r} \leq \left(2+\left(1+\bar{\mu}\right)JT\right) C_{32} r^{7}\e^{-C_{33}r},
\end{align*}
which converges to $0$ exponentially fast and this completes the proof.

\end{APPENDICES}



\bibliographystyle{informs2014} 
\bibliography{fork_join} 

\newpage
\pagenumbering{arabic}
\numberwithin{theorem}{section}
\numberwithin{lemma}{section}
\numberwithin{assumption}{section}
\numberwithin{remark}{section}
\numberwithin{table}{section}
\numberwithin{figure}{section}
\setcounter{section}{0}
\renewcommand{\thesection}{EC.\arabic{section}}
\setcounter{equation}{0}
\renewcommand{\theequation}{EC.\arabic{equation}}

\begin{center}
\textbf{ELECTRONIC COMPANION}
\end{center}

This electronic companion is associated with the manuscript titled ``Control of Fork-Join Processing Networks with Multiple Job Types and Parallel Shared Resources''. The proofs of the lemmas which are used in the proof of Lemma \ref{l_good_set} are presented. We present some preliminary results in Section \ref{s_pre}. Then, we present the proofs of Lemmas \ref{l_tail_1}, \ref{l_tail_2}, \ref{l_tail_3}, and \ref{l_tail_4} in Sections \ref{l_tail_1_proof}, \ref{l_tail_2_proof}, \ref{l_tail_3_proof}, and \ref{l_tail_4_proof}, respectively.

\section{Preliminary Results}\label{s_pre}

We derive exponentially decaying tail bounds for sum of i.i.d. random variables. Let $\{X_n,n\in\N_+\}$ be a sequence of nonnegative and i.i.d. random variables such that $\E[X_1]=x>0$. Suppose that there exists an $\bar{\al}>0$ such that $\E\left[ \e^{\al X_1}\right]<\infty$ for all $\al\in(-\bar{\al},\bar{\al})$, that is, $X_1$ satisfies the exponential moment assumption (see Assumption \ref{a_moment}). For all $\alpha\in\R$, let 
\begin{equation}\label{eq_legendre_1}
\ell(\al):=\ln \E\left[ \exp\left\{\al\left(X_1-x\right)\right\} \right].
\end{equation}
Then, $\e^{\ell(\alpha)}<\infty$ for all $\al\in(-\bar{\al},\bar{\al})$ by the exponential moment assumption on $X_1$. For $y\geq 0$, let
\begin{equation}\label{eq_legendre_2}
\Lambda^{(1)}(y):=\sup_{\al\in(0,\bar{\al})}\left\{ \al y-\ell(\al)\right\},\qquad\Lambda^{(2)}(y):=\sup_{\al\in(0,\bar{\al})}\left\{ \al y-\ell(-\al)\right\}.
\end{equation}
Then, we have the following result.

\begin{lemma}\label{l_legendre}
Both $\Lambda^{(1)}$ and $\Lambda^{(2)}$ are convex and nondecreasing in $\R_+$, $\Lambda^{(1)}(0)=\Lambda^{(2)}(0)=0$, and $\Lambda^{(1)}(y)>0$ and $\Lambda^{(2)}(y)>0$ for all $y> 0$. 
\end{lemma}

\bProof
First, let us consider $\Lambda^{(1)}$. $\Lambda^{(1)}$ is convex because for any $\theta\in[0,1]$ and $y_1,y_2\in \R_+$,
\begin{align*}
\theta\Lambda^{(1)}(y_1) + (1-\theta)\Lambda^{(1)}(y_2) &= \sup_{\al\in(0,\bar{\al})}\left\{ \theta\al y_1-\theta\ell(\al)\right\} +\sup_{\al\in(0,\bar{\al})}\left\{(1-\theta)\al y_2-(1-\theta)\ell(\al)\right\}\\
&\geq \sup_{\al\in(0,\bar{\al})}\left\{ \al (\theta y_1+(1-\theta)y_2)-\ell(\al)\right\} = \Lambda^{(1)} (\theta y_1+(1-\theta)y_2).
\end{align*}

By Parts (a) and (c) of Lemma 2.2.5 of \citet{dem98}, $\ell$ is convex in $\R$, $\ell$ is differentiable in $(-\bar{\al},\bar{\al})$, and $\ell(0)=\ell'(0)=0$, where $\ell'$ is the derivative of $\ell$. Then, $\ell$ achieves the global minimum at $0$; and since it is convex, $\ell$ is nondecreasing in $\R_+$. Then,
\begin{equation*}
\Lambda^{(1)}(0)=-\inf_{\al\in(0,\bar{\al})}\ell(\al) = \ell(0)=0.
\end{equation*}
Furthermore, for any given $y>0$, there exists an $\al^*\in(0,\bar{\al})$ such that $\Lambda^{(1)}(y)\geq \al^*y-\ell(\al^*)>0$. Therefore, $\Lambda^{(1)}(y)>0$ for all $y>0$. 

For any given $y_1,y_2\in \R_+$ such that $y_1\leq y_2$, because $\Lambda^{(1)}$ is convex and $\Lambda^{(1)}(0)=0$, we have
\begin{equation*}
\Lambda^{(1)}(y_1) \leq \frac{y_1}{y_2}\Lambda^{(1)}(y_2) +\left(1-\frac{y_1}{y_2}\right)\Lambda^{(1)}(0) =  \frac{y_1}{y_2}\Lambda^{(1)}(y_2) \leq \Lambda^{(1)}(y_2).
\end{equation*}
Therefore, $\Lambda^{(1)}$ is nondecreasing in $\R_+$.

The proof for $\Lambda^{(2)}$ follows with exactly the same way, hence we skip it.\eProof

\begin{lemma}\label{l_tail_5}
Let $a$ and $b$ be arbitrary strictly positive constants. There exists an $r_{14}\in\N_+$ such that if $r\geq r_{14}$, then 
\begin{equation*}
\pr \left( \max_{N\in\{1,2,\ldots,\rd{ar}\}} \Bigg| \sum_{n=1}^N (X_n - x) \Bigg|>b r \right) \leq 2\e^{-C_{34} r},
\end{equation*}
where $C_{34}$ is a strictly positive constant independent of $r$.
\end{lemma}

\bProof
We have
\begin{align}\label{eq_tail_1}
&\pr \left( \max_{N\in\{1,2,\ldots,\rd{ar}\}}\Bigg| \sum_{n=1}^{N} (X_n - x) \Bigg|>b r \right) \nonumber\\
&\hspace{1cm} \leq \pr \left(\max_{N\in\{1,2,\ldots,\rd{ar}\}}\sum_{n=1}^{N} (X_n - x) >b r \right) + \pr \left(\max_{N\in\{1,2,\ldots,\rd{ar}\}} \sum_{n=1}^{N} (x-X_n) >b r \right).
\end{align}
Let $\al\in(0,\bar{\al})$ be an arbitrary constant. The first probability in the RHS of \eqref{eq_tail_1} is equal to
\begin{align}
&\pr \left(\max_{N\in\{1,2,\ldots,\rd{ar}\}}\sum_{n=1}^{N} \al(X_n - x) >\al b r \right) \nonumber\\
&\hspace{2cm}= \pr \left(\exp\left \{ \max_{N\in\{1,2,\ldots,\rd{ar}\}}\sum_{n=1}^{N} \al(X_n - x) \right\}> \e^{\al b r} \right)\nonumber\\
&\hspace{2cm}= \pr \left(\max_{N\in\{1,2,\ldots,\rd{ar}\}}\exp\left \{ \sum_{n=1}^{N} \al(X_n - x) \right\}> \e^{\al b r} \right)\nonumber\\
&\hspace{2cm} \leq  \E \left [\exp\left \{ \sum_{n=1}^{\rd{ar}} \al(X_n - x)\right\} \right] \e^{-\al b r} = \E \left [\prod_{n=1}^{\rd{ar}} \exp\left \{\al(X_n - x)\right\} \right] \e^{-\al b r}\label{eq_tail_3}\\
&\hspace{2cm}= \E \left [\exp\left \{\al(X_1 - x) \right\} \right]^{\rd{ar}} \e^{-\al b r} \nonumber\\
&\hspace{2cm}= e^{\ell(\al)\rd{ar}}\e^{-\al b r} = \exp\left\{ -\rd{ar} \left(\al \frac{b r}{\rd{ar}} - \ell(\al) \right)\right\}, \label{eq_tail_4}
\end{align}
where the inequality in \eqref{eq_tail_3} is by Doob's inequality for submartingales (see Theorem 5.4.2 of \citet{dur10}), and the first equality in \eqref{eq_tail_4} is by \eqref{eq_legendre_1}. Similarly, for the second probability in the RHS of \eqref{eq_tail_1}, we can derive that
\begin{equation}\label{eq_tail_5}
\pr \left(\max_{N\in\{1,2,\ldots,\rd{ar}\}} \sum_{n=1}^{N} (x-X_n) >b r \right) \leq \exp\left\{ -\rd{ar} \left(\al \frac{b r}{\rd{ar}} - \ell(-\al) \right)\right\}.
\end{equation}
By \eqref{eq_legendre_2} and because \eqref{eq_tail_4} and \eqref{eq_tail_5} hold for all $\al\in(0,\bar{\al})$, the RHS of \eqref{eq_tail_1} is less than or equal to
\begin{equation}\label{eq_tail_6}
\exp\left\{ -\rd{ar} \Lambda^{(1)} \left( \frac{b r}{\rd{ar}}\right)\right\} + \exp\left\{ -\rd{ar} \Lambda^{(2)} \left( \frac{b r}{\rd{ar}}\right)\right\}.
\end{equation}
There exists $r_{14}\in\N_+$ such that if $r \geq r_{14}$, we have $b r/\rd{ar}\geq 0.5b/a>0$ and $\rd{ar} \geq 0.5 a r$. By Lemma \ref{l_legendre}, $\Lambda^{(i)}(y) \geq \Lambda^{(i)}(0.5b/a)>0$ for all $y\geq 0.5b/a$ and $i\in\{1,2\}$. Therefore, the sum in \eqref{eq_tail_6} converges to 0 with exponential rate. To complete the proof, let 
\begin{equation*}
C_{34}:= \frac{a}{2} \left(\Lambda^{(1)}\left(\frac{b}{2a}\right)\wedge \Lambda^{(2)}\left(\frac{b}{2a}\right)\right).
\end{equation*}
\eProof

\section{Proof of Lemma \ref{l_tail_1}}\label{l_tail_1_proof}

We will only prove \eqref{eq_t1_2}. The proofs of \eqref{eq_t1_1} and \eqref{eq_t1_3} are very similar. Fix an arbitrary $i,l\in\mI$ and $j\in\mJ_i$. Let
\begin{equation*}
\eta^r:= \inf\left\{ t\in\R_+:  \left| S_{ij}(T_{ij}^r(\tau_{ln}^r) + t) -  S_{ij}(T_{ij}^r(\tau_{ln}^r)) -\mu_{ij}t \right| >br\right\},
\end{equation*}
where $\inf\{\emptyset\}=\infty$ for completeness. Let $B_{ln}^r:=\{\tau_{ln}^r \leq r^2T\}\cap B^r$. Then,
\begin{equation}\label{eq_lt1_1}
\pr\left( \sup_{0\leq t\leq ar} \left| S_{ij}(T_{ij}^r(\tau_{ln}^r) + t) -  S_{ij}(T_{ij}^r(\tau_{ln}^r)) -\mu_{ij}t \right| >br,\;B_{ln}^r\right) = \pr\left( \eta^r \leq ar,\;B_{ln}^r\right).
\end{equation}
Let
\begin{subequations}\label{eq_lt1_2}
\begin{align}
&V_{ij}^{(1),r} (Y):= \sum_{y= S_{ij}(T_{ij}^r(\tau_{ln}^r))+1}^{S_{ij}(T_{ij}^r(\tau_{ln}^r))+Y} v_{ijy},\qquad\forall Y\in\N_+,\\
&V_{ij}^{(2),r} (Y):= \sum_{y= S_{ij}(T_{ij}^r(\tau_{ln}^r))+2}^{S_{ij}(T_{ij}^r(\tau_{ln}^r))+Y} v_{ijy},\qquad\forall Y\in\{2,3,\dots\},\\
&V_{ij}^{(1),r} (-Y):=0,\quad\forall Y\in\N,\qquad V_{ij}^{(2),r} (-Y):=0,\quad \forall Y\in\N\cup\{-1\}.
\end{align}
\end{subequations}
Then, for all $t\in\R_+$,
\begin{align*}
&\left\{ \left| S_{ij}(T_{ij}^r(\tau_{ln}^r) + t) -  S_{ij}(T_{ij}^r(\tau_{ln}^r)) -\mu_{ij}t \right| >br \right\} \\
&\hspace{2cm} = \left\{ S_{ij}(T_{ij}^r(\tau_{ln}^r) + t) -  S_{ij}(T_{ij}^r(\tau_{ln}^r)) -\mu_{ij}t >br \right\}\\
&\hspace{6cm}\cup \left\{ S_{ij}(T_{ij}^r(\tau_{ln}^r) + t) -  S_{ij}(T_{ij}^r(\tau_{ln}^r)) -\mu_{ij}t < -br \right\} \\
&\hspace{2cm} \subseteq \left\{ S_{ij}(T_{ij}^r(\tau_{ln}^r) + t) -  S_{ij}(T_{ij}^r(\tau_{ln}^r))> \rd{br + \mu_{ij}t} \right\}\\
&\hspace{6cm}\cup \left\{ S_{ij}(T_{ij}^r(\tau_{ln}^r) + t) -  S_{ij}(T_{ij}^r(\tau_{ln}^r)) < \ru{-br +\mu_{ij}t} \right\} \\
&\hspace{2cm} \subseteq   \left\{ V_{ij}^{(2),r} \left(\rd{br + \mu_{ij}t}\right)<t \right\} \cup  \left\{ V_{ij}^{(1),r} \left(\ru{-br +\mu_{ij}t}\right)>t \right\}.
\end{align*}
Let us define
\begin{align*}
&\eta^r_1:= \inf\left\{ t\in\R_+: V_{ij}^{(2),r} \left(\rd{br + \mu_{ij}t}\right)<t\right\},\\
&\eta^r_2:= \inf\left\{ t\in\R_+: V_{ij}^{(1),r} \left(\ru{-br +\mu_{ij}t}\right)>t\right\}.
\end{align*}
Then $\eta^r\geq \eta^r_1\wedge \eta^r_2$. Hence,
\begin{equation}\label{eq_lt1_3}
\pr\left( \eta^r \leq ar,\;B_{ln}^r\right) \leq \pr\left( \eta^r_1 \leq ar,\;B_{ln}^r\right) + \pr\left( \eta^r_2 \leq ar,\;B_{ln}^r\right).
\end{equation}

First,
\begin{align}
&\pr\left( \eta^r_1 \leq ar,\;B_{ln}^r\right) =\pr\left( \inf_{0\leq t\leq ar} \left\{V_{ij}^{(2),r} \left(\rd{br + \mu_{ij}t}\right) - t\right\} <0,\;B_{ln}^r\right) \nonumber\\
&\hspace{0cm}\leq  \pr\left( \min_{y\in\left\{\rd{br},\ldots,\rd{(b +\mu_{ij}a)r}\right\}} \left\{V_{ij}^{(2),r} (y) - \frac{y+1-br}{\mu_{ij}}\right\} <0,\;B_{ln}^r\right) \nonumber \\
&\hspace{0cm}=  \pr\left( \min_{y\in\left\{\rd{br},\ldots,\rd{(b +\mu_{ij}a)r}\right\}} \left\{V_{ij}^{(2),r} (y) - \frac{y-1}{\mu_{ij}}\right\} <\frac{2-br}{\mu_{ij}},\;B_{ln}^r\right) \nonumber \\
&\hspace{0cm}=  \pr\left( \max_{y\in\left\{\rd{br},\ldots,\rd{(b +\mu_{ij}a)r}\right\}} \left\{\frac{y-1}{\mu_{ij}}-V_{ij}^{(2),r} (y)\right\} >\frac{br-2}{\mu_{ij}},\;B_{ln}^r\right) \nonumber \\
&\hspace{0cm}\leq \sum_{y_1=0}^{\rd{C_8r^2T} - 1} \pr\left( \max_{y\in\left\{\rd{br},\ldots,\rd{(b +\mu_{ij}a)r}\right\}} \left\{\frac{y-1}{\mu_{ij}}-V_{ij}^{(2),r} (y)\right\} >\frac{br-2}{\mu_{ij}},\;S_{ij}(T_{ij}^r(\tau_{ln}^r))=y_1\right) \nonumber\\
&\hspace{0cm}\leq C_8r^2T \pr\left( \max_{y\in\left\{\rd{br},\ldots,\rd{(b +\mu_{ij}a)r}\right\}} \left\{\frac{y-1}{\mu_{ij}}-V_{ij} (y-1)\right\} >\frac{br-2}{\mu_{ij}}\right) \nonumber\\
&\hspace{0cm}\leq C_8r^2T \pr\left( \max_{y\in\left\{1,\ldots,\rd{(b +\mu_{ij}a)r}\right\}} \left\{\frac{y}{\mu_{ij}}-V_{ij} (y)\right\} >\frac{br-2}{\mu_{ij}}\right). \label{eq_lt1_4}
\end{align}

Second,
\begin{align}
&\pr\left( \eta^r_2 \leq ar,\;B_{ln}^r\right) =\pr\left( \sup_{0\leq t\leq ar} \left\{V_{ij}^{(1),r} \left(\ru{-br +\mu_{ij}t}\right) - t\right\} >0,\;B_{ln}^r\right) \nonumber\\
&\hspace{0cm}\leq  \pr\left( \max_{y\in\left\{1,\ldots,\left\lceil(b+a\mu_{ij})r\right\rceil\right\}} \left\{V_{ij}^{(1),r} (y) - \frac{y+br-1}{\mu_{ij}}\right\} >0,\;B_{ln}^r\right) \nonumber \\
&\hspace{0cm}= \pr\left( \max_{y\in\left\{1,\ldots,\left\lceil(b+a\mu_{ij})r\right\rceil\right\}} \left\{V_{ij}^{(1),r} (y) - \frac{y}{\mu_{ij}}\right\} >\frac{br-1}{\mu_{ij}},\;B_{ln}^r\right) \nonumber \\
&\hspace{0cm}\leq \sum_{y_1=0}^{\rd{C_8r^2T} - 1} \pr\left(\max_{y\in\left\{1,\ldots,\left\lceil(b+a\mu_{ij})r\right\rceil\right\}} \left\{V_{ij}^{(1),r} (y) - \frac{y}{\mu_{ij}}\right\} >\frac{br-1}{\mu_{ij}},\;S_{ij}(T_{ij}^r(\tau_{ln}^r))=y_1\right) \nonumber\\
&\hspace{0cm}\leq C_8r^2T \pr\left(\max_{y\in\left\{1,\ldots,\left\lceil(b+a\mu_{ij})r\right\rceil\right\}} \left\{V_{ij} (y) - \frac{y}{\mu_{ij}}\right\} >\frac{br-1}{\mu_{ij}}\right).\label{eq_lt1_5}
\end{align}

By Lemma \ref{l_tail_5}, there exists an $r_{ij}\in\N_+$ such that if $r\geq r_{ij}$, the sum of the terms in \eqref{eq_lt1_4} and \eqref{eq_lt1_5} is less than or equal to 
\begin{equation}\label{eq_lt1_6}
4C_8T r^2\e^{-C_{ij}r},
\end{equation}
where $C_{ij}$ is a strictly positive constant independent of $l$, $n$, and $r$. Finally, \eqref{eq_t1_2} follows from \eqref{eq_lt1_1}, \eqref{eq_lt1_3}, \eqref{eq_lt1_6}, and letting $r_1:= \max_{j\in\mJ,i\in\mI_j} r_{ij}$, $C_9:=4C_8T$, and $C_{10}:=\min_{j\in\mJ,i\in\mI_j} C_{ij}$.

\section{Proof of Lemma \ref{l_tail_2}}\label{l_tail_2_proof}

Recall that Step 2 lasts at most as much as the sum of a residual inter-arrival time and a service time. Hence,
\begin{align}
&\pr\left(\tau_{i(n+1)}^r - \tau_{in}^r>b r,\; \tau_{in}^r\leq r^2T,\;B^r\right) \nonumber\\
&\hspace{2cm}\leq \pr\left( \max_{l\in\left\{1,2,\ldots,\rd{C_8r^2T}\right\}}\max_{j\in\mJ_i} \left(u_{jl}^r+v_{ijl}\right) > b r \right)\nonumber\\
&\hspace{2cm}\leq \sum_{l=1}^{\rd{C_8r^2T}} \sum_{j\in\mJ_i}\pr\left(u_{jl}^r+v_{ijl}> b r \right)\nonumber\\
&\hspace{2cm}= C_8r^2T \sum_{j\in\mJ_i}\pr\left(u_{j1}^r+v_{ij1}> b r \right) =C_8r^2T \sum_{j\in\mJ_i}\pr\left(\frac{\bar{u}_{j1}}{\lambda_j^r}+v_{ij1}> b r \right)\nonumber\\
&\hspace{2cm}\leq C_8r^2T  \sum_{j\in\mJ_i} \left(\pr\left(\bar{u}_{j1}> 0.25\lambda_jb r \right)+ \pr\left( v_{ij1}> 0.5b r \right)\right)\label{eq_l_tail_2_1}\\
&\hspace{2cm}\leq C_8r^2T  \sum_{j\in\mJ_i} \left( \pr\left(\e^{0.5\bar{\al}\bar{u}_{j1}}> \e^{0.1\bar{\al}\lambda_jb r} \right)+ \pr\left( \e^{0.5\bar{\al}v_{ij1}}> \e^{0.25\bar{\al}b r} \right)\right)\nonumber\\
&\hspace{2cm}\leq C_8r^2T \sum_{j\in\mJ_i} \left( \E\left[\e^{0.5\bar{\al}\bar{u}_{j1}}\right]  \e^{-0.1\bar{\al}\lambda_jb r} +  \E\left[\e^{0.5\bar{\al}v_{ij1}}\right] \e^{-0.25\bar{\al}b r} \right)\label{eq_l_tail_2_2}\\
&\hspace{2cm}\leq C_{11}r^2\e^{-C_{12} r}, \label{eq_l_tail_2_3}
\end{align}
where \eqref{eq_l_tail_2_1} is by the fact that $r\geq r_2$, \eqref{eq_l_tail_2_2} is by Markov's inequality, and
\begin{equation*}
C_{11}:=C_8T J\max_{i\in\mI,j\in\mJ_i} \left\{\E\left[ \e^{0.5\bar{\al}\bar{u}_{j1}}\right]+ \E\left[ \e^{0.5\bar{\al}v_{ij1}}\right]\right\},\qquad C_{12}:= \bar{\al}b\left(0.25\wedge \left(0.1\underline{\lambda}\right)\right).
\end{equation*}
Notice that both $C_{11}$ and $C_{12}$ are strictly positive and finite constants by Assumption \ref{a_moment} and independent of $i$, $n$, and $r$.

\section{Proof of Lemma \ref{l_tail_3}}\label{l_tail_3_proof}

By Step 3 of Definition \ref{d_policy}, we have
\begin{equation*}
\tau_{i(n+1)}^r - \tau_{in}^r =\inf\Bigg\{t\geq 0: \sum_{j\in\mJ_i\backslash\{m\}} \sum_{x=1}^{\left(Q_{ij}^r(\tau_{in}^r)-\lceil q_{ij}^{*,r}(\tau_{in}^r)\rceil+A_j^r(\tau_{in}^r+t)-A_j^r(\tau_{in}^r)\right)^+} v_{ij(S_{ij}(T_{ij}^r(\tau_{in}^r))+x)} - t\leq 0 \Bigg\}.
\end{equation*}

For notational convenience, let
\begin{align*}
&\Delta_{ij}^{n,r}:=Q_{ij}^r(\tau_{in}^r) -\ru{q_{ij}^{*,r}(\tau_{in}^r)},&& L_{ij}^{n,r} := \frac{q_{ij}^{*,r}(\tau_{in}^r)-Q_{ij}^r(\tau_{in}^r)}{\lambda_j^r},\\
&B_{in}^{r}:=\left\{\max_{j\in\mJ_i}|\Delta_{ij}^{n,r}| \leq ar,\;\tau_{in}^r\leq r^2T,\;B^r\right\},&& \ot{V}_{ijn}^r(X):=\sum_{x=1}^{X} v_{ij(S_{ij}(T_{ij}^r(\tau_{in}^r))+x)},
\end{align*}
for all $X\in\N_+$ and $j\in\mJ$ such that $\sum_{x=y}^{z}v_{ijx}:=0$ for all $y>z$. Then, the probability in the LHS of \eqref{eq_l_tail_3_0} is equal to
\begin{align}
&\pr\left( \tau_{i(n+1)}^r - \tau_{in}^r > \frac{-\Delta_{im}^{n,r}}{\lambda_m}+ br,\; B_{in}^{r} \right)\label{eq_l_tail_3_00}\\
&\hspace{2cm}+\pr\left( Q_{im}^r(\tau_{i(n+1)}^r) -\ru{q_{im}^{*,r}(\tau_{in}^r)} > 2b\lambda_m r,\;\tau_{i(n+1)}^r - \tau_{in}^r \leq \frac{-\Delta_{im}^{n,r}}{\lambda_m}+ br,\; B_{in}^{r} \right)\label{eq_l_tail_3_01}.
\end{align}

By definition of Step 3 (see Definition \ref{d_policy}) and \eqref{eq_we}, $\Delta_{im}^{n,r}<0$. Let
\begin{equation*}
B_{inm}^{(1),r}(y_1) := B_{in}^{r} \cap \left\{\Delta_{im}^{n,r}=-y_1\right\},\quad\forall y_1\in\{1,2,\ldots,\rd{ar}\}.
\end{equation*}

\textbf{\textit{The probability in \eqref{eq_l_tail_3_00}}} It is equal to
\begin{align}
&\sum_{y_1=1}^{\rd{ar}}\pr\left( \tau_{i(n+1)}^r - \tau_{in}^r > \frac{y_1}{\lambda_m}+ br,\; B_{inm}^{(1),r}(y_1)\right)\nonumber\\
&=\sum_{y_1=1}^{\rd{ar}}\pr\Bigg(\inf_{t\geq 0} \Bigg\{ \sum_{j\in\mJ_i\backslash\{m\}} \sum_{x=1}^{\left(Q_{ij}^r(\tau_{in}^r)-\lceil q_{ij}^{*,r}(\tau_{in}^r) \rceil+A_j^r(\tau_{in}^r+t)-A_j^r(\tau_{in}^r)\right)^+} v_{ij(S_{ij}(T_{ij}^r(\tau_{in}^r))+x)} - t\leq 0 \Bigg\} \nonumber\\
& \hspace{12cm}> \frac{y_1}{\lambda_m}+ b r,\;B_{inm}^{(1),r}(y_1)\Bigg) \nonumber\\
&=\sum_{y_1=1}^{\rd{ar}}\pr\Bigg(\inf_{t\geq 0} \Bigg\{ \sum_{j\in\mJ_i\backslash\{m\}} \ot{V}_{ijn}^r\left(\Delta_{ij}^{n,r}+A_j^r(\tau_{in}^r+t)-A_j^r(\tau_{in}^r)\right) - t\leq 0 \Bigg\} > \frac{y_1}{\lambda_m}+ b r,\;B_{inm}^{(1),r}(y_1)\Bigg)\nonumber\\
&=\sum_{y_1=1}^{\rd{ar}}\pr\Bigg(\inf_{0\leq t\leq \frac{y_1}{\lambda_m}+ b r} \Bigg\{ \sum_{j\in\mJ_i\backslash\{m\}} \ot{V}_{ijn}^r\left(\Delta_{ij}^{n,r}+A_j^r(\tau_{in}^r+t)-A_j^r(\tau_{in}^r)\right) - t \Bigg\}>0,\;B_{inm}^{(1),r}(y_1)\Bigg)\nonumber\\
&\leq\sum_{y_1=1}^{\rd{ar}}\pr\Bigg( \sum_{j\in\mJ_i\backslash\{m\}} \ot{V}_{ijn}^r\left(\Delta_{ij}^{n,r}+A_j^r\left(\tau_{in}^r+\frac{y_1}{\lambda_m}+ b r\right)-A_j^r(\tau_{in}^r)\right)>\frac{y_1}{\lambda_m}+ b r,\;B_{inm}^{(1),r}(y_1)\Bigg).\label{eq_l_tail_3_1}
\end{align}

By Assumption \ref{a_ht} Part 1, there exists an $r_{15}\in\N_+$ such that if $r\geq r_{15}$, we have
\begin{equation*}
\frac{1}{\lambda_m^r} -\frac{1}{\lambda_m} \leq \frac{br}{2\rd{ar}}.
\end{equation*}
Then, if $r\geq r_{15}$, in the set $B_{inm}^{(1),r}(y_1)$, we have
\begin{equation*}
L_{im}^{n,r} + 0.5br \leq \frac{\ru{q_{im}^{*,r}(\tau_{in}^r)}-Q_{im}^r(\tau_{in}^r)}{\lambda_m^r} + 0.5b r \leq \frac{\ru{q_{im}^{*,r}(\tau_{in}^r)}-Q_{im}^r(\tau_{in}^r)}{\lambda_m}+ b r=\frac{y_1}{\lambda_m}+ b r.
\end{equation*}
Hence, by Lemma \ref{l_review}, if $r\geq r_2\vee r_{15}$, for all $y_1\in\{1,2,\ldots,\rd{ar}\}$, in $B_{inm}^{(1),r}(y_1)$, we have
\begin{align}
\frac{y_1}{\lambda_m}+ b r=& \sum_{j\in\mJ} \frac{\left(\lambda_j^r \left(\frac{y_1}{\lambda_m}+ b r\right)  - q_{ij}^{*,r}(\tau_{in}^r)+Q_{ij}^r(\tau_{in}^r)\right)^+}{\mu_{ij}}\nonumber\\
\geq& \sum_{j\in\mJ_i\backslash\{m\}} \frac{\left(\lambda_j^r\left(\frac{y_1}{\lambda_m}+ b r\right)  - q_{ij}^{*,r}(\tau_{in}^r)+Q_{ij}^r(\tau_{in}^r)\right)^+}{\mu_{ij}} \nonumber\\
&\hspace{4cm}+  \frac{\left(\lambda_m^r \left(L_{im}^{n,r} + 0.5br\right)  - q_{im}^{*,r}(\tau_{in}^r)+Q_{im}^r(\tau_{in}^r)\right)^+}{\mu_{im}} \nonumber\\
\geq& \sum_{j\in\mJ_i\backslash\{m\}} \frac{\left(\lambda_j^r \left(\frac{y_1}{\lambda_m}+ b r\right)  - \ru{q_{ij}^{*,r}(\tau_{in}^r)}+Q_{ij}^r(\tau_{in}^r)\right)^+}{\mu_{ij}} + \frac{\lambda_m br}{4\mu_{im}}.\label{eq_l_tail_3_2}
\end{align}
Therefore, by \eqref{eq_l_tail_3_2}, if $r\geq r_2\vee r_{15}$, the probability in \eqref{eq_l_tail_3_1} is less than or equal to
\begin{align}
 &\sum_{y_1=1}^{\rd{ar}}\sum_{j\in\mJ_i\backslash\{m\}}  \pr\Bigg( \ot{V}_{ijn}^r\left(\Delta_{ij}^{n,r}+A_j^r\left(\tau_{in}^r+\frac{y_1}{\lambda_m}+ b r\right)-A_j^r(\tau_{in}^r)\right) - \frac{\left(\lambda_j^r\left(\frac{y_1}{\lambda_m}+ b r\right) +\Delta_{ij}^{n,r}\right)^+}{\mu_{ij}}\nonumber\\
 &\hspace{11cm}  > \frac{\lambda_m b r}{4J\mu_{im}}, \;B_{inm}^{(1),r}(y_1)\Bigg).\label{eq_l_tail_3_3}
\end{align}

Let $c>0$ be an arbitrary constant such that
\begin{equation*}
c <\left[\frac{\underline{\lambda}b}{8J\bar{\mu}} \left( \frac{1.25\bar{\lambda}}{\underline{\mu}} \left( \frac{a}{\underline{\lambda}}+b \right)\right)^{-1} \right] \wedge \frac{4b\underline{\lambda}}{a+b\underline{\lambda}}.
\end{equation*}
By Assumption \ref{a_ht} Part 1, there exists an $r_{16}\in\N_+$ such that if $r\geq r_{16}$, for all $j\in\mJ_i$, we have
\begin{align*}
&(1-0.25c)\lambda_j < \lambda_j^r < (1+0.25c)\lambda_j,\\
& c\lambda_jbr >8,\\
&\frac{1}{\mu_{ij}}\left((1+c)\lambda_j - \lambda_j^r\right)\left(\frac{\rd{ar}}{\lambda_m}+ b r\right)< \frac{\lambda_m b r}{8J\mu_{im}}.
\end{align*}
Hence, if $r\geq r_{16}$, for all $j\in\mJ_i$, $y_1\in\{1,\ldots,\rd{ar}\}$, and $\omega\in\Omega$, we have,
\begin{align}
& \bigg\lfloor (1+c)\lambda_j\left(\frac{y_1}{\lambda_m}+ b r\right) \bigg\rfloor -  \lambda_j^r\left(\frac{y_1}{\lambda_m}+ b r\right) -1 \geq 0.75c\lambda_jbr -2 > 0.5c\lambda_jbr,\label{eq_l_tail_3_4}\\
&\frac{1}{\mu_{ij}}\left(\left( \bigg\lfloor (1+c)\lambda_j\left(\frac{y_1}{\lambda_m}+ b r\right) \bigg\rfloor +\Delta_{ij}^{n,r}(\omega)\right)^+ -\left(\lambda_j^r\left(\frac{y_1}{\lambda_m}+ b r\right) +\Delta_{ij}^{n,r}(\omega)\right)^+\right)\nonumber\\
&\hspace{7cm} \leq \frac{1}{\mu_{ij}}\left((1+c)\lambda_j - \lambda_j^r\right)\left(\frac{y_1}{\lambda_m}+ b r\right)< \frac{\lambda_m b r}{8J\mu_{im}}.\label{eq_l_tail_3_5}
\end{align}
Next, let us define the set
\begin{equation}\label{eq_l_tail_3_6}
B_{in}^{(2),r}(y_1):= \bigcap_{j\in\mJ_i}\left\{ A_j^r\left(\tau_{in}^r+\frac{y_1}{\lambda_m}+ b r\right)-A_j^r(\tau_{in}^r) \leq  \bigg\lfloor (1+c)\lambda_j\left(\frac{y_1}{\lambda_m}+ b r\right) \bigg\rfloor \right\}
\end{equation}
for all $y_1\in\{1,2,\ldots,\rd{ar}\}$. If $r\geq r_{16}$,
\begin{align}
&\sum_{y_1=1}^{\rd{ar}}\pr \left( \left( B_{in}^{(2),r}(y_1)\right)^c \cap B_{inm}^{(1),r}(y_1)\right) \nonumber\\
&\hspace{1cm} \leq \sum_{j\in\mJ_i}\sum_{y_1=1}^{\rd{ar}}\pr \Bigg( A_j^r\left(\tau_{in}^r+\frac{y_1}{\lambda_m}+ b r\right)-A_j^r(\tau_{in}^r) >  \bigg\lfloor (1+c)\lambda_j\left(\frac{y_1}{\lambda_m}+ b r\right) \bigg\rfloor,\nonumber\\
&\hspace{9cm} \tau_{in}^r\leq r^2T,\;A_j^r(r^2T)\leq \rd{C_8r^2T}-1\Bigg) \nonumber\\
&\hspace{1cm} \leq \sum_{j\in\mJ_i}\sum_{y_1=1}^{\rd{ar}}\pr \Bigg( \sum_{x=2}^{\big\lfloor (1+c)\lambda_j\left(\frac{y_1}{\lambda_m}+ b r\right) \big\rfloor} \bar{u}_{j(A_j^r(\tau_{in}^r)+x)} < \lambda_j^r\left(\frac{y_1}{\lambda_m}+ b r\right),\nonumber\\
&\hspace{9cm} \tau_{in}^r\leq r^2T,\;A_j^r(r^2T)\leq \rd{C_8r^2T}-1\Bigg) \nonumber\\
&\hspace{1cm} \leq \sum_{j\in\mJ_i}\sum_{y_1=1}^{\rd{ar}} \sum_{y_2=0}^{\rd{C_8r^2T}-1}\pr \Bigg( \sum_{x=2}^{\big\lfloor (1+c)\lambda_j\left(\frac{y_1}{\lambda_m}+ b r\right) \big\rfloor} \bar{u}_{j(y_2+x)} < \lambda_j^r\left(\frac{y_1}{\lambda_m}+ b r\right),\;A_j^r(\tau_{in}^r)=y_2\Bigg) \nonumber\\
&\hspace{1cm} \leq C_8r^2T\sum_{j\in\mJ_i}\sum_{y_1=1}^{\rd{ar}} \pr \Bigg( \sum_{x=1}^{\big\lfloor (1+c)\lambda_j\left(\frac{y_1}{\lambda_m}+ b r\right) \big\rfloor -1} (\bar{u}_{jx}-1)\nonumber\\
&\hspace{7cm} < \lambda_j^r\left(\frac{y_1}{\lambda_m}+ b r\right) -\bigg\lfloor (1+c)\lambda_j\left(\frac{y_1}{\lambda_m}+ b r\right) \bigg\rfloor +1\Bigg) \nonumber\\
&\hspace{1cm} \leq C_8r^2T\sum_{j\in\mJ_i}\sum_{y_1=1}^{\rd{ar}} \pr \Bigg( \Bigg| \sum_{x=1}^{\big\lfloor (1+c)\lambda_j\left(\frac{y_1}{\lambda_m}+ b r\right) \big\rfloor -1} (\bar{u}_{jx}-1) \Bigg| >0.5c\lambda_jbr \Bigg) \label{eq_l_tail_3_7}\\
&\hspace{1cm} \leq C_8ar^{3}T\sum_{j\in\mJ_i}\pr \Bigg(\max_{y_3\in\left\{1,2,\ldots,\big\lfloor (1+c)\lambda_j\left(\frac{\rd{ar}}{\lambda_m}+ b r\right)\big\rfloor\right\}} \Bigg| \sum_{x=1}^{y_3} (\bar{u}_{jx}-1) \Bigg| >0.5c\lambda_jbr \Bigg) \nonumber\\
&\hspace{1cm} \leq 2C_8JTar^{3}\e^{-C_{35}r} \label{eq_l_tail_3_8},
\end{align}
where $C_{35}$ is a strictly positive constant independent of $n$ and $r$, \eqref{eq_l_tail_3_7} is by \eqref{eq_l_tail_3_4}, and \eqref{eq_l_tail_3_8} is by Lemma \ref{l_tail_5} and holds for all $r\geq r_{17}$ such that $r_{17}\in\N_+$ is a constant independent of $n$.

By \eqref{eq_l_tail_3_5} and \eqref{eq_l_tail_3_6}, if $r\geq r_{16}$, the sum in \eqref{eq_l_tail_3_3} is less than or equal to
\begin{align}
 &\sum_{y_1=1}^{\rd{ar}}\sum_{j\in\mJ_i\backslash\{m\}}  \pr\Bigg( \ot{V}_{ijn}^r\left(\Delta_{ij}^{n,r} +  \bigg\lfloor (1+c)\lambda_j\left(\frac{y_1}{\lambda_m}+ b r\right) \bigg\rfloor \right) - \frac{\left(\Delta_{ij}^{n,r}+ \Big\lfloor (1+c)\lambda_j\left(\frac{y_1}{\lambda_m}+ b r\right) \Big\rfloor \right)^+}{\mu_{ij}}\nonumber\\
 &\hspace{2cm}  > \frac{\lambda_m b r}{8J\mu_{im}}, \;B_{inm}^{(1),r}(y_1),\;B_{in}^{(2),r}(y_1)\Bigg) +  \sum_{y_1=1}^{\rd{ar}} J\pr \left( \left( B_{in}^{(2),r}(y_1)\right)^c \cap B_{inm}^{(1),r}(y_1)\right) .\label{eq_l_tail_3_9}
\end{align}

The first sum in \eqref{eq_l_tail_3_9} is less than or equal to
\begin{align}
 &\sum_{j\in\mJ_i\backslash\{m\}} \sum_{y_1=1}^{\rd{ar}}\sum_{y_2=0}^{\rd{C_8r^2T}-1} \sum_{y_3=-\rd{ar}}^{\rd{ar}} \pr\Bigg( \ot{V}_{ijn}^r\left(y_3 +  \bigg\lfloor (1+c)\lambda_j\left(\frac{y_1}{\lambda_m}+ b r\right) \bigg\rfloor \right)\nonumber\\
 &\hspace{2.5cm}  - \frac{\left(y_3+ \Big\lfloor (1+c)\lambda_j\left(\frac{y_1}{\lambda_m}+ b r\right) \Big\rfloor \right)^+}{\mu_{ij}}> \frac{\lambda_m b r}{8J\mu_{im}},\; S_{ij}(T_{ij}^r(\tau_{in}^r))=y_2,\;\Delta_{ij}^{n,r}=y_3\Bigg) \nonumber\\
&\leq\sum_{j\in\mJ_i\backslash\{m\}} \sum_{y_1=1}^{\rd{ar}}\sum_{y_2=0}^{\rd{C_8r^2T}-1} \sum_{y_3=-\rd{ar}}^{\rd{ar}} \pr\left(\sum_{x=1}^{y_3 +  \big\lfloor (1+c)\lambda_j\left(\frac{y_1}{\lambda_m}+ b r\right) \big\rfloor  } \left(v_{ij(y_2+x)} -\frac{1}{\mu_{ij}} \right) > \frac{\lambda_m b r}{8J\mu_{im}}\right)\nonumber\\
&\leq C_8r^2T\sum_{j\in\mJ_i\backslash\{m\}} \sum_{y_1=1}^{\rd{ar}}\sum_{y_3=-\rd{ar}}^{\rd{ar}} \pr\left(\sum_{x=1}^{y_3 +  \big\lfloor (1+c)\lambda_j\left(\frac{y_1}{\lambda_m}+ b r\right) \big\rfloor } \left(v_{ijx} -\frac{1}{\mu_{ij}} \right) > \frac{\lambda_m b r}{8J\mu_{im}}\right)\nonumber\\
&\leq C_8r^2T\left(2ar+1\right)^2\sum_{j\in\mJ_i\backslash\{m\}}  \pr\left( \max_{y_4\in\left\{1,2,\ldots,\rd{ar}+ \big\lfloor (1+c)\lambda_j\left(\frac{\rd{ar}}{\lambda_m}+ b r\right) \big\rfloor\right\}}\sum_{x=1}^{y_4} \left(v_{ijx} -\frac{1}{\mu_{ij}} \right) > \frac{\lambda_m b r}{8J\mu_{im}}\right)\nonumber\\
&\leq C_8r^2T\left(2ar+1\right)^2 J2\e^{-C_{36}r},\label{eq_l_tail_3_10}
\end{align}
where \eqref{eq_l_tail_3_10} is by Lemma \ref{l_tail_5}, $C_{36}$ is a strictly positive constant independent of $n$ and $r$, and \eqref{eq_l_tail_3_10} holds if $r\geq r_{18}$ for some $r_{18}\in\N_+$ such that $r_{18}$ is a constant independent of $n$.

Therefore, by \eqref{eq_l_tail_3_8}, \eqref{eq_l_tail_3_9}, and \eqref{eq_l_tail_3_10}, if $r\geq \max\{r_{15},r_{16},r_{17},r_{18}\}$, the probability in \eqref{eq_l_tail_3_00} is less than or equal to 
\begin{equation}\label{eq_l_tail_3_11}
2C_8J^2T\left(4a^2+5a+1\right)r^{4} \e^{-\left( C_{35}\wedge C_{36}\right)r}.
\end{equation}

\textbf{\textit{The probability in \eqref{eq_l_tail_3_01}}} By definition of Step 3 (see Definition \ref{d_policy}),
\begin{equation*}
Q_{im}^r(\tau_{i(n+1)}^r) =  Q_{im}^r(\tau_{in}^r) + A_m^r(\tau_{i(n+1)}^r) - A_m^r(\tau_{in}^r).
\end{equation*}
Hence, 
\begin{equation*}
Q_{im}^r(\tau_{i(n+1)}^r) -\ru{q_{im}^{*,r}(\tau_{in}^r)} = A_m^r(\tau_{i(n+1)}^r) - A_m^r(\tau_{in}^r) + \Delta_{im}^{n,r}.
\end{equation*}
Therefore, the probability in \eqref{eq_l_tail_3_01} is equal to 
\begin{align}
&\pr\left( A_m^r(\tau_{i(n+1)}^r) - A_m^r(\tau_{in}^r) > -\Delta_{im}^{n,r} + 2b\lambda_mr,\;\tau_{i(n+1)}^r - \tau_{in}^r \leq \frac{-\Delta_{im}^{n,r}}{\lambda_m}+ br,\; B_{in}^{r} \right)\nonumber \\
&\hspace{1cm} \leq \pr\left( A_m^r\left(\tau_{in}^r + \frac{-\Delta_{im}^{n,r}}{\lambda_m}+ br\right) - A_m^r(\tau_{in}^r) > -\Delta_{im}^{n,r} + 2b\lambda_mr,\; B_{in}^{r} \right)\nonumber \\
&\hspace{1cm} =\sum_{y_1=1}^{\rd{ar}}\pr\left( A_m^r\left(\tau_{in}^r + \frac{y_1}{\lambda_m}+ br\right) - A_m^r(\tau_{in}^r) > y_1 + 2b\lambda_mr,\; B_{inm}^{(1),r}(y_1)\right) \label{eq_l_tail_3_12}.
\end{align}
Similar to how we derive the bound in \eqref{eq_l_tail_3_8}, we can prove that there exists an $r_{19}\in\N_+$ independent of $n$ such that if $r\geq r_{19}$, the sum in \eqref{eq_l_tail_3_12} is less than or equal to 
\begin{equation}\label{eq_l_tail_3_13}
2C_8Tar^{3}\e^{-C_{37}r},
\end{equation}
where $C_{37}$ is a strictly positive constant independent of $n$ and $r$.

Finally, let $r_3:=\max\{r_{15},r_{16},r_{17},r_{18},r_{19}\}$. Then, $r_3$ is independent of $n$ and $\omega$. By \eqref{eq_l_tail_3_11} and \eqref{eq_l_tail_3_13}, if $r\geq r_3$, the probability in the LHS of \eqref{eq_l_tail_3_0} is less than or equal to $C_{13}r^{4}\e^{-C_{14}r}$, where $C_{13}:= 2C_8J^2T(4a^2+6a+1)$ and $C_{14}:= C_{35}\wedge C_{36}\wedge C_{37}$.

\section{Proof of Lemma \ref{l_tail_4}}\label{l_tail_4_proof}

Let us fix arbitrary $a>0$ and $b>0$. The probability in \eqref{eq_t4_0} is less than or equal to
\begin{align}
&\sum_{i\in\mI} \pr\left(\sup_{0\leq t\leq ar } \left|W_i^{r}(\tau_{ln}^r+t)- W_i^{r}(\tau_{ln}^r) \right| >b r,\;\tau_{ln}^r \leq r^2T,\;B^r\right) \label{eq_t4_1}\\
&\hspace{3cm}+\sum_{j\in\mJ}\sum_{k\in\mK_j} \pr\left(\sup_{0\leq t\leq ar } \left|Q_k^{r}(\tau_{ln}^r+ t)- Q_k^{r}(\tau_{ln}^r) \right| >b r,\;\tau_{ln}^r \leq r^2T,\;B^r\right). \label{eq_t4_2}
\end{align}

Let us focus on the sum in \eqref{eq_t4_1} first. For all $j\in\mJ$, $i\in\mI_j$, $l\in\mI$, $n\in\{1,2,\ldots,M^r\}$, $r\in\N_+$ and $t\in\R_+$, let us define the shifted processes
\begin{align*}
T_{ij}^{l,n,r}(t) &:= T_{ij}^{r}(\tau_{ln}^r + t) -  T_{ij}^{r}(\tau_{ln}^r) \\
I_{i}^{l,n,r}(t) &:= I_{i}^{r}(\tau_{ln}^r + t) -  I_{i}^{r}(\tau_{ln}^r) \\
A_{j}^{l,n,r}(t) &:= A_{j}^{r}(\tau_{ln}^r + t) -  A_{j}^{r}(\tau_{ln}^r) -\lambda_j^r t \\
S_{ij}^{l,n,r}(t) &:= S_{ij}(T_{ij}^r(\tau_{ln}^r + t)) -  S_{ij}(T_{ij}^r(\tau_{ln}^r)) -\mu_{ij} T_{ij}^{l,n,r}(t) \\
Q_{ij}^{l,n,r}(t) &:= Q_{ij}^{r}(\tau_{ln}^r + t) \\
X_{ij}^{l,n,r}(t) &:= Q_{ij}^{r}(\tau_{ln}^r) + A_{j}^{l,n,r}(t)- S_{ij}^{l,n,r}(t)\\
W_{i}^{l,n,r}(t) &:= W_{i}^{r}(\tau_{ln}^r + t) = \sum_{j\in\mJ_i} \frac{Q_{ij}^{l,n,r}(t) }{\mu_{ij}}\\
\rho_i^r &:=  \sum_{j\in\mJ_i} \frac{\lambda_j^r}{\mu_{ij}},
\end{align*}
By some algebra, for all $i,l\in\mI$, $n\in\{1,2,\ldots,M^r\}$, $r\in\N_+$ and $t\in\R_+$, we have 
\begin{align}
&W_{i}^{l,n,r}(t)  = \sum_{j\in\mJ_i} \frac{X_{ij}^{l,n,r}(t) }{\mu_{ij}} + \left(\rho_i^r -1\right)t +I_{i}^{l,n,r}(t),\nonumber\\
&\left(W_{i}^{l,n,r},I_{i}^{l,n,r}\right) = \left(\Phi,\Psi \right) \left( \sum_{j\in\mJ_i} \frac{X_{ij}^{l,n,r}}{\mu_{ij}} + \left(\rho_i^r -1\right)e \right).\label{eq_t4_3}
\end{align}

Then, by \eqref{eq_t4_3},
\begin{align}
&\sup_{0\leq t\leq ar } \left|W_i^{r}(\tau_{ln}^r+t)- W_i^{r}(\tau_{ln}^r) \right| =\sup_{0\leq t\leq ar } \left|W_{i}^{l,n,r}(t) -W_i^{r}(\tau_{ln}^r) \right| \nonumber\\
&\hspace{0.5cm} =\sup_{0\leq t\leq ar } \Bigg| \sum_{j\in\mJ_i} \frac{X_{ij}^{l,n,r}(t) }{\mu_{ij}} +  \left(\rho_i^r -1\right)t + \sup_{0\leq s \leq t} \Bigg( -\sum_{j\in\mJ_i} \frac{X_{ij}^{l,n,r}(s) }{\mu_{ij}} -  \left(\rho_i^r -1\right)s\Bigg)^+ - W_i^{r}(\tau_{ln}^r)\Bigg| \nonumber\\
&\hspace{0.5cm} =\sup_{0\leq t\leq ar } \Bigg| \sum_{j\in\mJ_i} \frac{A_{j}^{l,n,r}(t)- S_{ij}^{l,n,r}(t)}{\mu_{ij}} + \left(\rho_i^r -1\right)t  \nonumber\\
&\hspace{5cm} + \sup_{0\leq s \leq t} \Bigg( \sum_{j\in\mJ_i} \frac{-Q_{ij}^{r}(\tau_{ln}^r) - A_{j}^{l,n,r}(s)+ S_{ij}^{l,n,r}(s) }{\mu_{ij}} -  \left(\rho_i^r -1\right)s\Bigg)^+ \Bigg| \nonumber \\
&\hspace{0.5cm} \leq  \sup_{0\leq t\leq ar } \Bigg| \left(\rho_i^r -1\right)t + \sup_{0\leq s \leq t} \Bigg( \sum_{j\in\mJ_i} \frac{-Q_{ij}^{r}(\tau_{ln}^r)}{\mu_{ij}} -  \left(\rho_i^r -1\right)s\Bigg)^+ \Bigg| \label{eq_t4_4}\\
&\hspace{8cm} + 2\sum_{j\in\mJ_i} \frac{1} {\mu_{ij}}  \left(\big\Vert A_{j}^{l,n,r}\big\Vert_{ar} + \big\Vert S_{ij}^{l,n,r}\big\Vert_{ar}\right). \label{eq_t4_5}
\end{align}
If $\rho_i^r\geq 1$, then the term in \eqref{eq_t4_4} is equal to
\begin{equation}\label{eq_t4_6}
 \left(\rho_i^r -1\right)ar.
\end{equation}
If $\rho_i^r < 1$, then the term in \eqref{eq_t4_4} is equal to
\begin{align}
&\sup_{0\leq t\leq ar } \Bigg| - \left(1-\rho_i^r \right)t +  \Bigg(\left(1-\rho_i^r\right)t  -\sum_{j\in\mJ_i} \frac{Q_{ij}^{r}(\tau_{ln}^r)}{\mu_{ij}}\Bigg)^+ \Bigg| \nonumber\\
&\hspace{5cm}= \sup_{0\leq t\leq ar } \Bigg|  \left(\left(1-\rho_i^r\right)t \right) \wedge \sum_{j\in\mJ_i} \frac{Q_{ij}^{r}(\tau_{ln}^r)}{\mu_{ij}} \Bigg| \leq \left(1-\rho_i^r\right)ar. \label{eq_t4_7}
\end{align}
Therefore, by \eqref{eq_t4_6} and \eqref{eq_t4_7}, the sum of the terms in \eqref{eq_t4_4} and \eqref{eq_t4_5} is less than or equal to
\begin{equation*}
 \left|\rho_i^r -1\right|ar + 2\sum_{j\in\mJ_i} \frac{1} {\mu_{ij}}  \left(\big\Vert A_{j}^{l,n,r}\big\Vert_{ar} + \big\Vert S_{ij}^{l,n,r}\big\Vert_{ar}\right).
\end{equation*}
Therefore, the sum in \eqref{eq_t4_1} is less than or equal to
\begin{equation}\label{eq_t4_8}
\sum_{i\in\mI} \pr\left( \left|\rho_i^r -1\right|ar + 2\sum_{j\in\mJ_i} \frac{1} {\mu_{ij}}  \left(\big\Vert A_{j}^{l,n,r}\big\Vert_{ar} + \big\Vert S_{ij}^{l,n,r}\big\Vert_{ar}\right) >b r,\;\tau_{ln}^r \leq r^2T,\;B^r\right).
\end{equation}
By Assumption \ref{a_ht} Parts 1 and 2, there exists an $r_{20}\in\N_+$ such that if $r\geq r_{20}$,
\begin{equation}\label{eq_t4_8_0}
\left|\rho_i^r -1\right| \leq \frac{b}{2a},\qquad \forall i\in\mI.
\end{equation}
Therefore, by \eqref{eq_t4_8_0}, if $r\geq r_{20}$, the sum in \eqref{eq_t4_8} is less than or equal to
\begin{align}
& \sum_{i\in\mI} \pr\left( \sum_{j\in\mJ_i} \frac{1} {\mu_{ij}}  \left(\big\Vert A_{j}^{l,n,r}\big\Vert_{ar} + \big\Vert S_{ij}^{l,n,r}\big\Vert_{ar}\right) >0.25 b r,\;\tau_{ln}^r \leq r^2T,\;B^r\right)\nonumber\\
&\hspace{2cm}\leq \sum_{i\in\mI}\sum_{j\in\mJ_i}  \pr\left( \left(\big\Vert A_{j}^{l,n,r}\big\Vert_{ar} + \big\Vert S_{ij}^{l,n,r}\big\Vert_{ar}\right) >\frac{b\underline{\mu}}{4J} r,\;\tau_{ln}^r \leq r^2T,\;B^r\right)\nonumber\\
&\hspace{2cm}\leq \sum_{i\in\mI}\sum_{j\in\mJ_i}  \pr\left(\big\Vert A_{j}^{l,n,r}\big\Vert_{ar} >\frac{b\underline{\mu}}{8J} r,\;\tau_{ln}^r \leq r^2T,\;B^r\right)\label{eq_t4_9}\\
&\hspace{4cm}+ \sum_{i\in\mI}\sum_{j\in\mJ_i}  \pr\left( \big\Vert S_{ij}^{l,n,r}\big\Vert_{ar}>\frac{b\underline{\mu}}{8J} r,\;\tau_{ln}^r \leq r^2T,\;B^r\right)\label{eq_t4_10}.
\end{align}

First, let us consider the sum in \eqref{eq_t4_9}, which is equal to
\begin{equation}\label{eq_t4_11}
\sum_{i\in\mI}\sum_{j\in\mJ_i}  \pr\left(\sup_{0\leq t\leq ar} \left| A_{j}^{r}(\tau_{ln}^r + t) -  A_{j}^{r}(\tau_{ln}^r) -\lambda_j^r t\right| >\frac{b\underline{\mu}}{8J} r,\;\tau_{ln}^r \leq r^2T,\;B^r\right).
\end{equation}
By Lemma \ref{l_tail_1}, there exists an $r_{21}\in\N_+$ such that $r_{21}$ is independent of $l$ and $n$ and if $r\geq r_{21}$, then the sum in \eqref{eq_t4_11} is less than or equal to
\begin{equation}\label{eq_t4_12}
C_{38}r^2\e^{-C_{39}r},
\end{equation}
where $C_{38}$ and $C_{39}$ are strictly positive constants independent of $l$, $n$, and $r$.

Second, let us consider the sum in \eqref{eq_t4_10}. By definition, we have
\begin{equation}\label{eq_t4_13}
\big\Vert S_{ij}^{l,n,r}\big\Vert_{ar} = \sup_{0\leq t\leq ar} \left| S_{ij}(T_{ij}^r(\tau_{ln}^r + t)) -  S_{ij}(T_{ij}^r(\tau_{ln}^r)) -\mu_{ij} \left(T_{ij}^{r}(\tau_{ln}^r + t) - T_{ij}^{r}(\tau_{ln}^r)\right)\right|.
\end{equation}

For all $j\in\mJ$, $i\in\mI_j$, $l\in\mI$, $n\in\{1,2,\ldots,M^r\}$, $r\in\N_+$, and $t\in[0,ar]$, because $0\leq T_{ij}^{r}(\tau_{ln}^r + t) - T_{ij}^{r}(\tau_{ln}^r)\leq t$, there exists $f_{ij}^{l,n,r}(t)\in[0,t]$ such that $T_{ij}^{r}(\tau_{ln}^r + t)=T_{ij}^{r}(\tau_{ln}^r)+f_{ij}^{l,n,r}(t)$. Then, for all $t\in[0,ar]$,
\begin{align*}
&\left| S_{ij}(T_{ij}^r(\tau_{ln}^r + t)) -  S_{ij}(T_{ij}^r(\tau_{ln}^r)) -\mu_{ij} \left(T_{ij}^{r}(\tau_{ln}^r + t) - T_{ij}^{r}(\tau_{ln}^r)\right)\right|\\
&\hspace{2cm} = \left| S_{ij}(T_{ij}^r(\tau_{ln}^r )+f_{ij}^{l,n,r}(t)) -  S_{ij}(T_{ij}^r(\tau_{ln}^r)) -\mu_{ij} f_{ij}^{l,n,r}(t)\right|\\
&\hspace{2cm} \leq \sup_{0\leq s\leq t}\left| S_{ij}(T_{ij}^r(\tau_{ln}^r )+s) -  S_{ij}(T_{ij}^r(\tau_{ln}^r)) -\mu_{ij}s \right|\\
&\hspace{2cm} \leq \sup_{0\leq s\leq ar}\left| S_{ij}(T_{ij}^r(\tau_{ln}^r )+s) -  S_{ij}(T_{ij}^r(\tau_{ln}^r)) -\mu_{ij}s \right|.
\end{align*}
By \eqref{eq_t4_13} and the fact that the last inequality above holds uniformly for all $t\in[0,ar]$, we have
\begin{equation*}
\big\Vert S_{ij}^{l,n,r}\big\Vert_{ar} \leq \sup_{0\leq t\leq ar} \left| S_{ij}(T_{ij}^r(\tau_{ln}^r )+t) -  S_{ij}(T_{ij}^r(\tau_{ln}^r)) -\mu_{ij}t \right|.
\end{equation*}
Therefore, the sum in \eqref{eq_t4_10} is less than or equal to
\begin{equation}\label{eq_t4_14}
\sum_{i\in\mI}\sum_{j\in\mJ_i}  \pr\left(\sup_{0\leq t\leq ar} \left| S_{ij}(T_{ij}^r(\tau_{ln}^r )+t) -  S_{ij}(T_{ij}^r(\tau_{ln}^r)) -\mu_{ij}t \right| >\frac{b\underline{\mu}}{8J} r,\;\tau_{ln}^r \leq r^2T,\;B^r\right).
\end{equation}
By Lemma \ref{l_tail_1}, there exists an $r_{22}\in\N_+$ such that $r_{22}$ is independent of $l$ and $n$ and if $r\geq r_{22}$, then the sum in \eqref{eq_t4_14} is less than or equal to
\begin{equation}\label{eq_t4_15}
C_{40}r^2\e^{-C_{41}r},
\end{equation}
where $C_{40}$ and $C_{41}$ are strictly positive constants independent of $l$, $n$, and $r$.

Consequently, by \eqref{eq_t4_12} and \eqref{eq_t4_15}, if $r\geq r_{20}\vee r_{21}\vee r_{22}$, then the sum in \eqref{eq_t4_1} is less than or equal to
\begin{equation}\label{eq_t4_16}
\left(C_{38}+C_{40}\right)r^2\e^{-\left(C_{39}\wedge C_{41}\right)r}.
\end{equation}

Next, let us consider the sum in \eqref{eq_t4_2}. Let $\ot{\mK}_j^L:=\{k\in\mK_j:\lambda_j<\mu_k\}$ and $\ot{\mK}_j^H:=\{k\in\mK_j:\lambda_j=\mu_k\}$ for all $j\in\mJ$. Then, $\ot{\mK}_j^L\subset \mK_j^L$ and $\ot{\mK}_j^L \cup \ot{\mK}_j^H= \mK_j$ for all $j\in\mJ$ by Assumption \ref{a_ht} Parts 1 and 4. $\ot{\mK}_j^L$ denotes the set of dedicated servers associated with the job type $j$ that are in light traffic and whose corresponding limiting arrival rate is strictly less than its service rate. Then, the sum in \eqref{eq_t4_2} is equal to
\begin{align}
&\sum_{j\in\mJ}\sum_{k\in\ot{\mK}_j^H} \pr\left(\sup_{0\leq t\leq ar } \left|Q_k^{r}(\tau_{ln}^r+ t)- Q_k^{r}(\tau_{ln}^r) \right| >b r,\;\tau_{ln}^r \leq r^2T,\;B^r\right) \label{eq_t4_17}\\
&\hspace{2cm} +\sum_{j\in\mJ}\sum_{k\in\ot{\mK}_j^L} \pr\left(\sup_{0\leq t\leq ar } \left|Q_k^{r}(\tau_{ln}^r+ t)- Q_k^{r}(\tau_{ln}^r) \right| >b r,\;\tau_{ln}^r \leq r^2T,\;B^r\right). \label{eq_t4_18}
\end{align}
Similar to how we derive \eqref{eq_t4_16}, we can prove that there exists an $r_{23}\in\N_+$ such that $r_{23}$ is independent of $l$ and $n$ and if $r\geq r_{23}$, then the sum in \eqref{eq_t4_17} is less than or equal to
\begin{equation}\label{eq_t4_19}
C_{42}r^2\e^{-C_{43}r},
\end{equation}
where $C_{42}$ and $C_{43}$ are strictly positive constants independent of $l$, $n$, and $r$. 

However, we cannot use the same technique to derive an exponential tail bound for the sum in \eqref{eq_t4_18}. Because $\lambda_j^r\rightarrow \lambda_j < \mu_k$ for all $j\in\mJ$ and $k\in\ot{\mK}_j^L$, the inequality in \eqref{eq_t4_8_0} with $\rho_i^r$ replaced with $\lambda_j^r/\mu_k$ may not hold for the dedicated server $k\in\ot{\mK}_j^L$. Therefore, the term in the RHS of \eqref{eq_t4_7} becomes a very loose bound. Intuitively, if $k\in\ot{\mK}_j^L$ and $Q_k^r(\tau_{ln}^r)$ is too large, the dedicated server $k\in\ot{\mK}_j^L$ can process many jobs within $(\tau_{ln}^r,\tau_{ln}^r + ar)$ and so we can have $Q_k^r(\tau_{ln}^r) - Q_k^r(\tau_{ln}^r + ar)> br$. Therefore, we need show that $Q_k^r(\tau_{ln}^r)$ can never be too large for all $j\in\mJ$ and $k\in\ot{\mK}_j^L$. In fact, the sum in \eqref{eq_t4_18} is less than or equal to
\begin{equation}\label{eq_t4_20}
\sum_{j\in\mJ}\sum_{k\in\ot{\mK}_j^L} \pr\left(\sup_{0\leq t\leq r^2T + ar } Q_k^{r}(t) >b r\right). 
\end{equation}
By Proposition 5 of \citet{ozk19}, there exists an $r_{24}\in\N_+$ such that if $r\geq r_{24}$, the sum in \eqref{eq_t4_20} is less than or equal to
\begin{equation}\label{eq_t4_21}
C_{44}r^5\e^{-C_{45}r},
\end{equation}
where $C_{44}$ and $C_{45}$ are strictly positive constants independent of $l$, $n$, and $r$.

By \eqref{eq_t4_16}, \eqref{eq_t4_19}, and \eqref{eq_t4_21}, Lemma \ref{l_tail_4} follows from letting 
\begin{equation*}
r_4:= r_{20}\vee r_{21}\vee r_{22} \vee r_{23} \vee r_{24},\quad C_{15}:= C_{38} + C_{40} + C_{42} + C_{44},\quad C_{16}:= C_{39} \wedge C_{41} \wedge C_{43} \wedge C_{45}.
\end{equation*}

\end{document}